\documentclass[11pt,oneside, final]{amsart}
\copyrightinfo{0}{Iranian Mathematical Society}
\pagespan{1}{\pageref*{LastPage}}
\usepackage{etoolbox,lastpage}
\commby{}
\usepackage{amsmath,amsthm,amscd,amsfonts,amssymb,enumerate}
\usepackage{graphicx}
\usepackage{color}
\usepackage[colorlinks]{hyperref}
\usepackage{esint}
\usepackage{empheq}

\makeatletter \@addtoreset{equation}{section}

\newtheorem{theorem}{Theorem}[section]
\newtheorem{proposition}[theorem]{Proposition}
\newtheorem{lemma}[theorem]{Lemma}
\newtheorem{corollary}[theorem]{Corollary}

\theoremstyle{definition}
\newtheorem{definition}[theorem]{Definition}

\newtheorem{remark}[theorem]{Remark}

\theoremstyle{remark}
\numberwithin{equation}{section}
\DeclareMathOperator{\esssup}{ess sup}
\DeclareMathOperator{\essinf}{ess inf}

\begin{document}


\title[Horizontal inverse mean curvature flow in the Heisenberg group]{Horizontal inverse mean curvature flow in the Heisenberg group$^{**}$}

\author[J. Cui]{Jingshi Cui}

\author[P. Zhao]{Peibiao Zhao$^*$}

\thanks{$^*$Corresponding author}
%
\begin{abstract}
    Huisken and Ilmanen [J. Differential Geom., 2001] created the theory of weak solutions for inverse mean curvature flows (IMCF) of hypersurfaces on Riemannian manifolds, and proved successfully a Riemannian version of the Penrose inequality.

    The present paper investigates and constructs a sub-Riemannian version of the theory of weak solutions for inverse mean curvature flows of surfaces in the first Heisenberg group $\mathbb{H}^{1}$. The level set formulation of the IMCF in $\mathbb{H}^{1}$ is given by
    \begin{equation}\label{0.1}
        \begin{cases}
            \operatorname{div}_{0}\left(\frac{\nabla_{0} u}{|\nabla_{0} u|}\right) =|\nabla_{0} u|,  \qquad  &\text{ in } \Omega,  \\
            u = 0,  \qquad  &\text{ on } \partial \Omega,
        \end{cases}
    \end{equation}
    where $\Omega \subset \mathbb{H}^{1}$ is an open set with smooth boundary, and $\Omega^{c} = \mathbb{H}^{1}\setminus \Omega = \{ u \leq 0\}$ is bounded.  Let $w_{p} = \exp \left( \frac{u_{p}}{1-p}\right)$ and  $w_{p}$ satisfies
    \begin{equation}\label{0.3}
        \begin{cases}
            \operatorname{div}_{0}\left(|\nabla_{0} w_{p}|^{p-2} \nabla_{0} w_{p}\right) = 0,  \qquad  &\text{ in } \Omega,  \\
            w_{p}=1,  \qquad  &\text{ on } \partial \Omega.
        \end{cases}
    \end{equation}
    Following the argument by Moser [J. Eur. Math. Soc., 2007], the key ingredient  in proving the existence of weak solutions to (\ref{0.1}) is to establish a uniform interior estimate for $|\nabla_{0} u_{p}|$. However, due to the lack of boundary continuity of  $|\nabla_{0} u_{p}| = (p-1)\frac{|\nabla_{0} w_{p}|}{w_{p}} \in C^{0,\beta}(\Omega)$ ($0< \beta <1, p>1$) by Zhong and Mukherjee [Anal. PDE, 2021], the standard method in [R. Moser, J. Eur. Math. Soc., 2007] cannot be applied to obtain a uniform interior estimate for $|\nabla_{0} u_{p}|$. Fortunately, the present paper discovers two refined inequalities: Harnack inequality and  Lipschitz estimate for $w_p$, which allow one to obtain interior estimates for $|\nabla_{0} u_{p}|$ independent of $p$. By further combining them with Arzel$\grave{\rm a} $-Ascoli theorem, the weak solution of (\ref{0.1}) can then be generated as the limit of $u_{p}$ as $p \to 1$, where $w_{p} = \exp \left( \frac{u_{p}}{1-p}\right)$ and $w_{p}$ is of solutions to (\ref{0.3}).
    
    As an important application of the IMCF in $\mathbb{H}^{1}$, a positive answer to an open problem posed in [F. Montefalcon, Ann. Mat. Pura Appl. (4), 2014]:Heintze-Karcher inequality in $\mathbb{H}^{1}$ is provided.
\end{abstract}

\maketitle
%

\section{Introduction}
Roughly speaking, the study of curvature flows, that is, the study of the evolution of surfaces depending on their curvature has a long history, of which the inverse mean curvature flow plays an important role. Given a closed manifold $M^{n}$, we say that a smooth family of immersions $X: M^{n} \times[0, T_{0}] \to \mathbb{R}^{n+1}$ is a classical solution of inverse mean curvature flow (IMCF) if
\begin{align}\label{1.1}
    \frac{\partial X}{\partial t}=\frac{\nu}{H},
\end{align}
where $H >0$ and $\nu$ are the mean curvature and the outward unit normal vector of the hypersurface $M_{t}:=X(M^{n}, t)$ at a point $X$, respectively.

Inverse mean curvature flow has many interesting and important applications in geometric analysis and general relativity.  In geometric analysis, it has been employed to address various important problems. For example, it has been used to prove the Poincar$\acute{\rm e}$ conjecture for manifolds with Yamabe invariant greater than those of $\mathbb{P} \mathbb{R}^{3}$ \cite{BN2004}, to establish the Minkowski inequality in the anti-de Sitter-Schwarzschild manifold \cite{BHW2016}, and to derive Alexandrov-Fenchel type inequalities in space forms \cite{DG2016,GP2017,MS2016, WX2015}. Other inequalities have also been established using this flow \cite{GWWX2015, LWX2014, LW2017}. In general relativity, Huisken and Ilmanen \cite{HI2001, HI2008} established a Riemannian version of the Penrose inequality using this flow, while Bray \cite{B2001} extended this result to a more general form of the Penrose inequality in Riemannian manifolds. These advancements mark significant progress toward resolving the well-known Penrose conjecture in general relativity.

It is well known that the inverse mean curvature flow in Euclidean space $\mathbb{R}^{n+1}$ does not produce singularities in finite time, provided that the initial hypersurface is a topological sphere, even if it is not convex. In fact, Gerhardt \cite{G1990} and Urbas \cite{U1990} proved that if the initial hypersurface is star-shaped with strictly positive mean curvature, then the classical solutions to (\ref{1.1}) remain smooth for all $t \geq 0$ and converge homothetically to round spheres after exponential rescaling as $t \to \infty$.

Another important observation is that there exist curvature flows that develop singularities within a finite time for general initial hypersurfaces. For instance, Colding and Minicozzi \cite{CM2012, CM2016} proved that singularities are unavoidable in mean curvature flow. Similarly, Harvie \cite{H2019} demonstrated that for non-topological spherical initial hypersurfaces, the inverse mean curvature flow typically develops intersections or singularities in finite time, leading to the breakdown of classical motions. This phenomenon is closely related to changes in the topology of the evolving surface, such as the formation of self-intersections or the merging of initially disjoint surface components. To address such challenges, Huisken and Ilmanen \cite{HI2001} introduced a weak solution for IMCF in Riemannian manifolds and established its existence, uniqueness, and connectedness. To the best of our knowledge, the study of IMCF in degenerate Riemannian spaces remains an open problem.

Sub-Riemannian spaces are degenerate Riemannian spaces in which the Riemannian inner product is defined only on subspaces of the tangent space. The simplest model of sub-Riemannian spaces is the first Heisenberg group $\mathbb{H}^{1}$, which is not Euclidean on any scale. In a sense, the Heisenberg group plays a role in sub-Riemannian manifolds analogous to that of Euclidean space in Riemannian manifolds. A natural idea is to study the inverse mean curvature flow in the first Heisenberg group. There are two folds for choosing the Heisenberg group as models: it is possible to define intrinsic derivatives and hence introduce the notion of horizontal mean curvature; it has a rich structure of dilations and translations, and some basic functional theory and partial differential equation theory have been established on the Heisenberg group.

In this paper we investigate IMCF in the first Heisenberg group, referred to as horizontal inverse mean curvature flow (HIMCF), which is obtained by replacing the geometric elements of the inverse mean curvature flow (\ref{1.1}) by the corresponding horizontal quantities in $\mathbb{H}^{1}$. We study the existence and uniqueness of weak solutions of the HIMCF and explore the geometric properties and applications.

Let $x : M \times [0,T_{0}) \to \mathbb{H}^{1}$ be a family of smooth embeddings, $M \subset \mathbb{H}^{1}$, $M_{0} = x(M,0)$ and $M_{t}= x(M,t)$. We say that $M_{t}\setminus \Sigma (M_{t})$ (where $\Sigma (M_{t})$ denotes the set of characteristic points of $M_{t}$) is a classical solution of the horizontal inverse mean curvature flow with initial surface $M_{0}$ if and only if, for any point $x\in M_{t} \setminus \Sigma (M_{t})$ it satisfies
\begin{align}\label{1.2}
    \left\langle \frac{\partial x}{\partial t}, \nu \right\rangle = H_{0}^{-1} \left\langle \nu_{0}, \nu \right\rangle ,
\end{align}
where $H_{0} > 0$ is the horizontal mean curvature, $\nu_{0}$ and $\nu$ are the horizontal unit outward normal vector and the unit outward normal vector of $M_{t}$, respectively. By Chow's connectivity theorem in \cite{M1996}, any point $x \in \mathbb{H}^{1}$ can be joined to the origin $o=(0,0,0)$ by a horizontal curve $\gamma(t) : [0,T_{0}] \to \mathbb{H}^{1}$. Since the evolution equation (\ref{1.2}) holds at every point $x \in M_{t} \setminus \Sigma (M_{t})$, it is possible to verify that any smooth horizontal curve $\gamma(t) \in M_{t} \setminus \Sigma (M_{t})$ satisfies
\begin{align}\label{1.3}
    \frac{\partial }{\partial t} \gamma = H_{0}^{-1} \nu_{0} .
\end{align}
In order to simplify the problem formulation, the system of equations can be transformed into a scalar PDE by using the parametrization method of differential geometry. The level set approach (see e.g. \cite{G2002, I1992, OS1988}) not only can describe the behaviour of singularities on $M_{t}$,  but also is a powerful tool for constructing weak solutions that allow singularities.

Given an initial surface $M_{0}$ which is the boundary of a bounded open set $E_{0}$, there exists a smooth function $u(x): \mathbb{H}^{1}  \to \mathbb{R}$ such that
\begin{align*}
    E_{0}= \{x \in \mathbb{H}^{1}: u(x) < 0 \},\qquad  M_{0} =  \partial E_{0},
\end{align*}
and the evolving surface $M_{t}$ is parametrized by
\begin{align*}
    E_{t}= \{x \in \mathbb{H}^{1}: u(x) < t \},\qquad  M_{t} =  \partial E_{t},
\end{align*}
where $M_{t}$ is allowed to contain singularities. Let $\Omega = \mathbb{H}^{1} \setminus E_{0}$ and $\partial \Omega = \partial E_{0}$. Then, the geometric evolution equation (\ref{1.2}) can be expressed as the subelliptic Dirichlet problem (\ref{0.1}):
\begin{equation*}
    \begin{cases}
        \operatorname{div}_{0}\left(\frac{\nabla_{0} u}{|\nabla_{0} u|}\right) =|\nabla_{0} u|,  \qquad  &\text{ in } \Omega,  \\
        u = 0,  \qquad  &\text{ on } \partial \Omega.
    \end{cases}
\end{equation*}
Notice that the equation (\ref{0.1}) is degenerate, which makes the analysis of the existence and uniqueness of solutions more delicate. Furthermore, (\ref{1.2}), (\ref{1.3}) and (\ref{0.1}) are not sufficient to describe the evolution perfectly, as they are not well defined in $\Sigma(M_{t})$. Due to the non-zero dimensions and the complexity of the structure of $\Sigma(M_{t})$, many of the tactics that prevailed in Euclidean spaces no longer worked. For instance, the classical elliptic regularization methods cannot be directly applied, and establishing a priori estimates for solutions of (\ref{0.1}) remains a challenging problem. To overcome the difficulties arising from the presence of characteristic points in the Heisenberg group, we adopt the notion of weak solutions to (\ref{1.2}) and introduce the generalized horizontal inverse mean curvature flow.

Since the left-hand side of (\ref{0.1}) is in divergence form, it is usual to consider the definition of a weak solution in the distributional sense (by integrating by parts of derivatives over a smooth test function) as follows
\begin{align*}
    \int_{\Omega} |\nabla_{0} u|^{-1} \left\langle \nabla_{0} u, \nabla_{0} \varphi \right\rangle dx = - \int_{\Omega} |\nabla_{0} u| \varphi  dx , \qquad \text{ for all } \varphi \in HW^{1,1}_{0} (\Omega).
\end{align*}
However, this way of defining weak solutions fails, the above formula does not make sense at the characteristic points where $|\nabla_{0} u| = 0$. In the spirit of \cite{HZ1994,LT1978,V1998}, we attempt to construct a suitable functional in the horizontal sense such that a weak solution of (\ref{1.2}) is defined as a minimizer of the functional
\begin{align}\label{1.5}
    J_{u}^{K}(v):=\int_{K} |\nabla_{0} v|+ v|\nabla_{0} u| dx ,
\end{align}
where $K$ is a compact set containing $\{v\neq u\}$. It can be verified that the minimizer of the functional (\ref{1.5}) satisfies the evolution equation (\ref{0.1}). In order to describe the behaviour of $M_{t}$ at singular points, it is necessary to introduce another notion of weak solution, which depends on the level set $F_{t} = \{ v < t \}$, as follows
\begin{align}\label{1.6}
    J_{u}^{K}(F_{t})=P_{\mathbb{H}}(F_{t},K) - \int_{F_{t} \cap K} |\nabla_{0} u|dx .
\end{align}
It can be shown that the  two forms of the  weak solution of (\ref{1.2}) mentioned above are equivalent (see Lemma \ref{lem3.3}).

Although it is possible to define a generalized horizontal inverse mean curvature flow (cf. Definition \ref{def3.5}) via weak solutions, in which the evolving surface may contain characteristic points, there is one more problem that needs to be solved: the evolving surface may develop singularities. For example, if $E_{0}$ is the union of two disjoint Kor$\acute{ \rm a}$nyi balls, then after a sufficiently long time, the horizontal mean curvature may tend to zero at some points, and singularities may develop. Similar to the IMCF in Euclidean space, we also adopt the concept of jumps and introduce appropriate horizontal jumps to address the singularity in the Heisenberg group. The horizontal jump procedure ensures that the HIMCF can pass through the singularity and continue its evolution thereafter. More precisely, the horizontal jump procedure removes the region in $M_{t}$ where $H_{0}=0$ and replaces it with a more regular region. Possibly, this modification combines the separated components of $E_{t}=\{u < t\}$ into a connected set. The horizontal minimizing hull $E_{t}^{'}$ is introduced to heuristically describe the surface after performing the horizontal jump.

Now we can investigate the existence and uniqueness of the weak solution to HIMCF. Using an approach similar to that in \cite{HI2001}, we establish a comparison principle for $E_{t}$, which implies the uniqueness of the weak solution to (\ref{0.1}). However, it is challenging to directly extend the elliptic regularization method in \cite{HI2001} to prove the existence of weak solutions in the Heisenberg group. Based on the approximation method in \cite{M2007}, a weak solution of (\ref{0.1}) can be obtained as the limit of weak solutions $u_{p}$ to (\ref{5.2a})-(\ref{5.2b}) as $p \to 1$. The key ingredient of this method is to establish a horizontal gradient bound for $u_{p}$ that is uniform in $p$, which is obtained by establishing a refined Harnack inequality and a refined Lipschitz estimate for the weak solution $w_{p}$ to (\ref{0.3}). Consequently, there exists a subsequence $p_{i} \to 1$ such that $u_{p_{i}} \to u$ locally uniformly in $\Omega$, where $u \in HW^{1,\infty}_{loc}(\Omega)$. With essentially the same argument as that in \cite{HI2001}, it follows that $u$ is a weak solution of (\ref{0.1}).
\begin{theorem}\label{thm1.1}
    Let $E_{0}$ be a bounded open set in $\mathbb{H}^{1}$ with smooth boundary $M_{0} = \partial E_{0}$. Then the initial value problem (\ref{0.1}) has a unique locally Lipschitz weak solution $u(x)$, i.e., the generalized HIMCF with the initial surface $M_{0}$ exists for all time and is unique.
\end{theorem}
\begin{remark}
    In \cite{PV2024}, the authors establish the existence of weak solutions to (\ref{0.1}) by combining Moser's method \cite{M2007} with a Riemannian approximation method, provided that $\Omega$ satisfies an exterior uniform gauge-ball condition.
\end{remark}

In particular, if the initial surface is the unit Kor$\acute{\rm a} $nyi sphere, then the evolving surface $M_{t}$ along the HIMCF can be expressed explicitly as
\begin{align*}
    M_{t} = \left\{x \in \mathbb{H}^{1}, \left(x_{1}^{2}+x_{2}^{2}\right)^{2}+16 x_{3}^{2} = e^{\frac{4}{3} t}\right\}.
\end{align*}

In the following, we analyze the geometrical properties of the generalized horizontal inverse mean curvature flow. Firstly, we investigate the regularity of the weak solution and show that it coincides with the classical solution for a short time under suitable initial conditions.

\begin{corollary}\label{cor1.2}
    Let $E_{0} \subset \mathbb{H}^{1}$ be a bounded open set with smooth noncharacteristic boundary $M_{0} = \partial E_{0}$. Assume that $H_{0}\big|_{M_{0}} > 0$ and $E_{0}$ is a horizontal minimizing hull. Then the weak solution $\{E_{t}\}_{0\leq t < \infty}$ of the initial value problem (\ref{0.1}) coincides with the smooth, classical solution for a short time, provided that $E_{t}$ remains precompact.
\end{corollary}

Secondly, we show that the generalized horizontal inverse mean curvature flow is an expanding flow, characterized by the exponential expansion of the horizontal perimeter $P_{\mathbb{H}}(E_{t})$.

\begin{corollary}\label{cor1.3}
    Let $E_{0} \subset \mathbb{H}^{1}$ be a bounded open set with smooth boundary. Suppose that $\left\{E_{t} \right\}_{t>0}$ is the weak solution of (\ref{0.1}) with initial value $E_{0}$, and $E_{t}$ remains precompact. Then $P_{\mathbb{H}}(E_{t}) = e^{t}P_{\mathbb{H}}(E_{0})$ for any $t> 0$. Moreover, if $E_{0}$ is a horizontal minimizing hull, then $P_{\mathbb{H}}(E_{t}) = e^{t}P_{\mathbb{H}}(E_{0})$ for any $t \geq 0$.
\end{corollary}

Notably, curvature flow serves as a powerful tool for establishing geometric inequalities. This naturally leads to the question of whether the HIMCF can be employed to derive geometric inequalities in the Heisenberg group.

The Heintze-Karcher inequality, originally introduced in \cite{HK1978}, plays a crucial role in exploring the relationship between curvature and geometric properties of manifolds.  It is well known that the classical proof of the Heintze-Karcher inequality relies on Reilly's formula \cite{R1977} (see, e.g., \cite{QX2015,R1980}). In 2014, Montefalcone \cite{M2014} successfully established a generalized integral formula that extends the famous Reilly's formula in \cite{R1977} to the Heisenberg group, and employing it to derive several geometric inequalities. However, as far as we know, the Heintze-Karcher type inequality cannot be derived directly from the generalized Reilly formula in \cite{M2014}. Utilizing the volume-increasing property of the HIMCF, we establish a Heintze-Karcher type inequality in $\mathbb{H}^1$.

\begin{theorem}(open problem)\label{thm1.4}
    Let $E$ be a bounded open set in $\mathbb{H}^{1}$ with smooth boundary $M = \partial E$. If $M \setminus \Sigma (M)$ has positive horizontal mean curvature, then
    \begin{align}\label{1.7}
    \int_{M \setminus \Sigma (M)} H_{0}^{-1} d\sigma_{\mathbb{H}} \geq \frac{4}{3} |E| .
    \end{align}
    Equality holds in (\ref{1.7}) when $M$ is a horizontal constant mean curvature surface.
\end{theorem}

\begin{remark}\label{rem1.6}
    From \cite{CDPT2007}, we know that the bubble set $\mathcal{S}(o,R)$ is a horizontal constant mean curvature surface in $\mathbb{H}^{1}$ and
    \begin{align*}
        H_{0} = \frac{2}{R},\qquad P_{\mathbb{H}}(\mathcal{S}(o,R)) = \frac{1}{2}\pi^{2}R^{3}, \qquad |\mathcal{S}(o,R)| = \frac{3}{16} \pi^{2} R^{4} .
    \end{align*}
    A direct calculation shows that the equality in (\ref{1.7}) holds for bubble sets. Recently, it was shown in \cite{FPV2026} that the constant $\frac{4}{3}$ in (\ref{1.7}) is not sharp in some sense, and the sharp one remains unknown. In future work, we plan to investigate suitable modifications of the HIMCF in order to determine the sharp constant in the Heintze-Karcher type inequality.
\end{remark}

Next, we consider the normalized flow obtained by rescaling the HIMCF through the Heisenberg dilation $\delta_{\lambda(t)}:\mathbb{H}^{1}\to\mathbb{H}^{1}$, where $\lambda(t)=e^{-\frac{1}{3}t}$. Let $E_{t}$ be the weak solution of the HIMCF, and define
\begin{align*}
    \widehat{E}_{t} = \delta_{\lambda(t)}(E_{t}) = \{y \in \mathbb{H}^{1}, y=\delta_{\lambda(t)}(x), x \in E_{t} \}=\{y \in \mathbb{H}^{1}, h(y,t)  < 0\}, \qquad \widehat{M}_{t}= \partial \widehat{E}_{t} .
\end{align*}
Then, the velocity function $\widehat{F}$ of $\widehat{M}_{t}  \setminus  \Sigma(\widehat{M}_{t})$ has the following form
\begin{align*}
    \widehat{F} = H_{0}^{-1}(y)\left\langle\widehat{\nu}_{0}, \widehat{\nu} \right\rangle  - \frac{1}{3}\left\langle y^{'} , \widehat{\nu} \right\rangle ,
\end{align*}
where $y^{'} = (y_{1}, y_{2}, 2y_{3})$, $\widehat{\nu}_{0}$ and $\widehat{\nu} $ are the horizontal unit outer normal vector and  the unit outer normal vector of $\widehat{M}_{t} \setminus \Sigma(\widehat{M}_{t})$ respectively. This normalization ensures that the horizontal perimeter of $\widehat{E}_{t}$ remains invariant (see Lemma \ref{lem8.4}). Let $y : M \times [0,T_{0}) \to \mathbb{H}^{1}$ be a family of smooth embeddings and $\widehat{M}_{t}= y(M,t)$. The $\mathbb{H}$-perimeter preserving flow is
\begin{align}\label{1.8}
    \left\langle \frac{\partial y}{\partial t} , \widehat{\nu} \right\rangle  = H_{0}^{-1}(y)\left\langle\widehat{\nu}_{0}, \widehat{\nu} \right\rangle  - \frac{1}{3}\left\langle y^{'} , \widehat{\nu} \right\rangle, \qquad  y \in \widehat{M}_{t} \setminus \Sigma(\widehat{M}_{t}).
\end{align}
Since the $\mathbb{H}$-perimeter preserving flow (\ref{1.8}) is obtained from the HIMCF (\ref{1.2}) via the Heisenberg dilation, the existence and uniqueness of $\widehat{E}_{t}$ follow directly from Theorem \ref{thm1.1}.

The classical Minkowski formula, as presented in \cite{CY1971,K1962}, is a celebrated integral formula that has proven instrumental in the study of minimal surfaces, curvature flows, and various geometric properties of hypersurfaces. By utilizing the $\mathbb{H}$-perimeter preserving flow (\ref{1.8}), a sub-Riemannian version of the Minkowski type formula is established for closed smooth surfaces in $\mathbb{H}^{1}$.

\begin{theorem}\label{thm1.7}
    Let $M = \partial E$ be a closed smooth surface in $\mathbb{H}^{1}$. Then
    \begin{align}\label{1.9}
        3P_{\mathbb{H}}(E) = \int_{M \setminus \Sigma (M)} H_{0} \left\langle y^{'} , \widehat{\nu} \right\rangle d\mathcal{H}^{2} .
    \end{align}
\end{theorem}
\begin{corollary}\label{cor1.8}
    Let $M$ be a closed, smooth, horizontal positive constant mean curvature surface, i.e., $H_{0}$ is a positive constant on the noncharacteristic locus. Then
    \begin{align}\label{1.10}
        3P_{\mathbb{H}}(E) = 4H_{0}|E| .
    \end{align}
\end{corollary}
\begin{remark}\label{rem1.9}
    In 2008, Manuel Ritor$\acute{\rm e}$ and C$\acute{\rm e}$sar Rosales \cite{RR2008} proved that the Minkowski formula (\ref{1.10}) holds for compact, connected, volume-preserving area-stationary $C^{2}$ surfaces. In particular, from Corollary 4.5 and Theorem 4.17 of \cite{RR2008}, it follows that if $M$ is a volume-preserving area-stationary $C^{2}$ surfaces, then $H_{0}|_{M\setminus \Sigma(M)}$ is constant. Moreover, the converse is true when  $\Sigma(M)$ consists only of isolated points.
\end{remark}

The present paper is organized as follows. In Section 2, we recall some preliminaries about Heisenberg group structure, level set representation and the horizontal perimeter. In Section 3, we introduce two notions of the weak solution of (\ref{0.1}) based on $u$ and $E_{t} =\{u<t\}$, respectively. Moreover, we prove that the two formulations of weak solutions are equivalent and construct an explicit solution to (\ref{0.1}). In Section 4, we introduce the generalized horizontal mean curvature, horizontal jumps, and horizontal minimizing hulls to analyze the behavior of HIMCF at characteristic points and singularities. Several properties of horizontal minimizing hulls are also established. In Sections 5 and 6, we provide the proofs for the existence and uniqueness of weak solutions to (\ref{0.1}), respectively. In Section 7, we study the geometric properties of the HIMCF and prove Theorem \ref{thm1.4}. In Section 8, we construct the $\mathbb{H}$-perimeter preserving flow (\ref{1.8}) and prove Theorem \ref{thm1.7} and Corollary \ref{cor1.8}.


\section{ Preliminaries}

\noindent This section presents the first Heisenberg group, introduces the level set representation of HIMCF, and reviews some well-known properties of the horizontal perimeter of a set that will be used in the present paper.

\subsection{ Heisenberg group structure}
\
\vglue-10pt
 \indent

The first Heisenberg group $\mathbb{H}^{1}$ is an analytic Lie group $(\mathbb{R}^{3},\cdot)$, where the noncommutative group multiplication is defined, for any points $x=\left(x_{1}, x_{2}, x_{3}\right)$, $y=\left(y_{1}, y_{2}, y_{3}\right) \in \mathbb{R}^{3}$, as
\begin{align*}
    \left(x_{1}, x_{2}, x_{3}\right)\cdot \left(y_{1}, y_{2}, y_{3}\right)=\left(x_{1}+y_{1}, x_{2}+y_{2}, x_{3}+y_{3}+\frac{1}{2}\left(x_{1} y_{2}-y_{1} x_{2}\right)\right).
\end{align*}
For $x \in \mathbb{H}^{1}$, the left translation by $x$ is the diffeomorphism $L_{x}(y) = x \cdot y$.  For any $\lambda > 0$, the non-isotropic group dilations $\delta_{\lambda}:\mathbb{H}^{1} \to \mathbb{H}^{1}$ are defined by
\begin{align*}
    \delta_{\lambda}(x) = \left(\lambda x_{1}, \lambda x_{2}, \lambda^{2} x_{3}\right).
\end{align*}
A basis of left-invariant vector fields is given by
\begin{align}\label{2.1}
    X_{1} = \partial_{x_{1}} - \frac{x_{2}}{2}\partial_{x_{3}}, \quad
    X_{2} = \partial_{x_{2}} + \frac{x_{1}}{2}\partial_{x_{3}}, \quad
    T = \partial_{x_{3}}.
\end{align}
We note that the only non-trivial commutator is $[X_{1},X_{2}] = X_{1}X_{2}-X_{2}X_{1} = T$. The horizontal distribution $\mathsf{H} \mathbb{H}^{1}$ in $\mathbb{H}^{1}$ is defined as the distribution spanned by the left-invariant vector fields $X_{1}$ and $X_{2}$, where $\mathsf{H} \mathbb{H}^{1} =  \{\mathsf{H}_{x} \mathbb{H}^{1}\}_{x \in\mathbb{H}^{1}}$ and $ \mathsf{H}_{x} \mathbb{H}^{1} = {\rm span} \{X_{1}(x), X_{2}(x)\}$. An absolutely continuous curve $\gamma : [0,T_{0}]\to \mathbb{H}^{1}$ is a horizontal curve whose tangent vector lies in the horizontal distribution.

Let $ g =\left\langle \cdot,\cdot\right\rangle $ be the Riemannian metric on $\mathbb{H}^{1}$ such that $\left\{X_{1}, X_{2}, T \right\}$ is an orthonormal basis at every point. From \cite[p.~304]{DGN2006}, $det(g_{ij}) =1$. Hence, the Riemannian volume form $dvol$ is given by the standard (Lebesgue) volume form $dx_{1} \wedge dx_{2}\wedge  dx_{3}$ in $\mathbb{R}^{3}$. The Haar measure in $\mathbb{H}^{1}$ is the Lebesgue measure in $\mathbb{R}^{3}$ up to a constant. The Haar measure of a Borel set $E \subset \mathbb{H}^{1}$ is denoted by $|E|$ and the average of an integrable function $u$ over set $E$ is
\begin{align*}
    \fint_{E} u dx = \frac{1}{|E|} \int_{E}  udx
\end{align*}

The Carnot-Carath$\acute{\rm e} $odory ($CC$) distance on $\mathbb{H}^{1}$ between two points is defined by
\begin{align*}
    d(x,y) = \inf \text{Length}(\gamma)
\end{align*}
where the infimum is taken over all horizontal curve connecting $x$ to $y$. Clearly, the $CC$-metric $d$ is invariant under left-translation and homogeneous of degree one, i.e.,
\begin{align*}
    d(L_{x_{0}}(x), L_{x_{0}}(y) ) = d(x,y),\qquad d (\delta_{\lambda}(x), \delta_{\lambda}(y)) = \lambda d(x,y).
\end{align*}
The homogeneous dimension of the Heisenberg group $\mathbb{H}^{n}$ is denoted by $Q$, where $Q=2n+2$. From the left invariance and scaling properties of the $CC $-metric, the Hausdorff measures in $(\mathbb{H}^{n}, d)$ satisfies
\begin{align*}
    \mathcal{H}^{Q}(L_{x_{0}}(E)) = \mathcal{H}^{Q}(E),\qquad \mathcal{H}^{Q}(\delta_{\lambda}(E)) = \lambda^{Q}\mathcal{H}^{Q}(E)
\end{align*}
for all $\lambda >0$, $x_{0} \in \mathbb{H}^{n}$, and $E \subset \mathbb{H}^{n}$. In particular, for all $x \in\mathbb{H}^{n}$ and $r>0$, there is
\begin{align*}
    \mathcal{H}^{Q}(B(x, r)) = r^{Q} \mathcal{H}^{Q}(B(o,1))
\end{align*}
where $B(x, r)$ denotes the $CC$-ball with center $x$ and radius $r$ in  $(\mathbb{H}^{n}, d)$. Indeed, $ (\mathbb{H}^{n}, d, \mathcal{H}^{Q})$ is an Ahlfors $Q$-regular space and $\mathcal{H}^{Q}$ agrees (up to a constant multiplicative factor) with the Haar measure on $\mathbb{H}^{n}$. Thus, $(\mathbb{H}^{n}, d)$ is a doubling metric space. See \cite{H2001}, Chapters 1-3, for more details.
\begin{lemma}\label{lemD}\cite{NSW1985} (Doubling property)
    Given a bounded set $U \subset \mathbb{H}^n$, there exist $R_0>0$ and $c >0$ such that for every $x \in U$, $R \leq R_0$ and $0<t<1$,
    \begin{align}\label{2.d}
        |B(x, t R)| \geq  ct^Q|B(x, R)|,
    \end{align}
    where $c$ is a constant independent of $x \in U$.
\end{lemma}

It is well know that $CC$-distance $d(x,y)$ is equivalent to the Kor$\acute{\text{a}}$nyi distance $\rho(x,y)$, i.e., there exists a constant $C =C(\mathbb{H}^{1}) >0$ such that
\begin{align*}
    C \rho(x,y) \leq d(x,y) \leq C^{-1}\rho(x,y),
\end{align*}
where $\rho(x,y)= \|  y^{-1} \cdot x \|_{\rho}$, and the Kor$\acute{\text{a}}$nyi norm $\|  \cdot \|_{\rho}$ is given by
\begin{align*}
    \|  x \|_{\rho}=\left[ \left(x_{1}^{2} + x_{2}^{2}\right)^{2} + 16x_{3}^{2}\right]^\frac{1}{4}.
\end{align*}
Also, the Kor$\acute{\text{a}}$nyi distance $\rho(x,y)$ is invariant under left translations and homogeneous of degree one under dilations, i.e.,
\begin{align}\label{2.2}
    \rho(L_{x_{0}}(x), L_{x_{0}}(y) ) = \rho(x,y),\qquad \rho (\delta_{\lambda}(x), \delta_{\lambda}(y)) = \lambda \rho(x,y).
\end{align}
The Kor$\acute{\text{a}}$nyi ball of radius $r>0$ centered at $x$ is defined by $\mathsf{B} _{r}(x):=\{y \in \mathbb{H}^{1}: \|y^{-1} \cdot x\|_{\rho}<r \} $, and the corresponding Kor$\acute{\text{a}}$nyi sphere is the boundary $\partial \mathsf{B} _{r}(x)$.

For any $C^{1}$ function $u:\mathbb{H}^{1} \to \mathbb{R}$, the gradient $\nabla u$ and the horizontal gradient $\nabla_{0} u$ are give by
\begin{align*}
    \nabla u =X_{1}u X_{1} + X_{2}u X_{2} + Tu T, \qquad \nabla_{0}u =X_{1}u X_{1} + X_{2}u X_{2},
\end{align*}
and
\begin{align*}
    &|\nabla u|^{2} =  \left\langle \nabla u, \nabla u\right\rangle =  (X_{1}u)^{2} + (X_{2}u)^{2} + (Tu)^{2}, \\
    &|\nabla_{0} u |^{2} =  \left\langle \nabla_{0} u, \nabla_{0} u\right\rangle =  (X_{1}u)^{2} + (X_{2}u)^{2}.
\end{align*}
The symmetrized second horizontal Hessian $u$ is
\begin{align*}
    \left(\nabla^{2}_{0}u\right)^{*} = \left({\begin{array}{cc}
        X^{2}_{1}u & \left(X_{1}X_{2}u + X_{2}X_{1}u\right)/2 \\
        \left(X_{1}X_{2}u + X_{2}X_{1}u\right)/2 & X^{2}_{2}u \\
    \end{array}} \right) : = [u_{,ij}] .
\end{align*}
Also, the horizontal Laplacian $\Delta_{0} u$ and the horizontal $\infty$-Laplacian $\Delta_{0,\infty} u$ are given by
\begin{align*}
    \Delta_{0} u = Tr \left(\nabla^{2}_{0}u\right)^{*} = X^{2}_{1}u + X^{2}_{2}u, \qquad  \Delta_{0,\infty} u = \sum_{i,j = 1}^{m=2}  u_{,ij}X_{i}u X_{j}u  = \frac{1}{2} \left\langle \nabla_{0} | \nabla_{0} u|^{2}, \nabla_{0} u\right\rangle .
\end{align*}
For a horizontal vector field $X = a X_{1} + bX_{2}$, the horizontal divergence of $X$ is defined by
\begin{align*}
    \operatorname{div}_{0}X = X_{1}a + X_{2}b.
\end{align*}

Let $V$ be the normal vector to a $C^{1}$ surface $M \subset \mathbb{H}^{1}$. The horizontal normal vector  $V_0$ is obtained by projecting $V$ onto the horizontal bundle at each point,
\begin{align*}
    V_{0}=V-\left\langle V, T\right\rangle T.
\end{align*}
The horizontal normal vector of smooth surfaces vanishes at some points, which are called characteristic points.
\begin{definition}\label{def2.5}
    Let $M$ be a $C^{1}$ surface in $\mathbb{H}^{1}$, the set of the characteristic points is
    \begin{align*}
        \Sigma(M) = \{x \in \mathbb{H}^{1} : V_{0}(x) = 0 \}.
    \end{align*}
\end{definition}
\begin{remark}\label{rem2.6}
    If $x$ is a characteristic point of $M$, then the tangent hyperplane $T_{x}M$ coincides with the horizontal distribution $\mathsf{H}_{x} \mathbb{H}^{1}$. $\Sigma(M)$ is nowhere dense in $M$ and the 2-dimensional Hausdorff measure of $\Sigma(M)$ is equal to zero.
\end{remark}
In the regular part $M \setminus \Sigma(M)$, the horizontal unit normal vector is defined as (see e.g. \cite{DGN2004,RR2006})
\begin{align}\label{2.3}
    \nu_{0} = \frac{V_{0}}{|V_{0}|}
\end{align}
and the horizontal mean curvature is given by
\begin{align}\label{2.4}
    H_{0}: =\sum_{i = 1}^{2} X_{i} \left\langle  \nu_{0}, X_{i}\right\rangle   =  \operatorname{div}_{0} (\nu_{0}).
\end{align}

\begin{definition}\label{def2.1}
    Let $\Omega \subset \mathbb{H}^{1}$ be an open set and $p \geq 1$. The horizontal Sobolev space $HW^{1,p}(\Omega)$ is the space of functions $u \in L_{p}(\Omega)$ with horizontal derivatives $\{X_{i}u\}_{i=1,2} \in L_{p}(\Omega)$, i.e., $|\nabla_{0} u | \in L_{p}(\Omega)$. The horizontal Sobolev space is a Banach space when it is equipped with the Sobolev norm
    \begin{align*}
        \| u \|_{HW^{1,p}(\Omega)} = \| u\|_{L_{p}(\Omega)} + \| \nabla_{0}u\|_{L_{p}(\Omega)} .
    \end{align*}
\end{definition}

\begin{remark}\label{rem2.2}
    The horizontal Sobolev space $HW^{1,p}(\Omega)$ is reflexive for $1<p<\infty$. By Theorem 2.7 of \cite{R2015}, $C^{\infty}(\Omega)\cap HW^{1,p}(\Omega)$ is dense in $HW^{1,p}(\Omega)$ for $1\leq p<\infty$. The space $HW_0^{1,p}(\Omega)$ is defined as the closure of $C_0^{\infty}(\Omega)$ in $HW^{1,p}(\Omega)$ with respect to the norm $\|\cdot\|_{HW^{1,p}(\Omega)}$. In addition, \cite{F2003} shows that every $u\in HW^{1,p}(\Omega)$ can be approximated by smooth functions. $HW^{1,p}_{\mathrm{loc}}(\Omega)$ denotes the local horizontal Sobolev space.
\end{remark}

The Sobolev inequality in the Heisenberg group is stated as follows.

\begin{theorem}\cite{CDG1993,CDG1994} \label{thmS1}
    Let $U \subset \mathbb{H}^{n}$ be a bounded open set. Let $1<p<Q$ and $1 \leq \kappa \leq \frac{Q}{Q-p}$. Then, there exists $R_0>0$ such that for any $x_{0} \in U, B_{R}=B(x_{0}, R)$, with $R \leq R_0$, we have
    \begin{align}\label{2.s}
        \left(\fint_{B_R}|u|^{\kappa p} d x\right)^{\frac{1}{\kappa p}} \leq \mathcal{C}_{s}(p,c,Q) R\left(\fint_{B_R}\left|\nabla_{0} u\right|^p d x\right)^{\frac{1}{p}}
    \end{align}
    for any $u \in HW^{1, p}\left(B_R\right)$, where $\mathcal{C}_{s}(p,c,Q)$ is a positive constant depending on $p$, $c$ and $Q$.
\end{theorem}

\subsection{ Parametrization by level sets }
\
\vglue-10pt
 \indent

Assume that there exists $u \in C^{1}(\mathbb{H}^{1})$ such that the evolving surfaces $M_{t}$ are represented as
\begin{align*}
    M_{t} = \partial E_{t} ,\qquad E_{t}=\{x\in \mathbb{H}^{1}: u(x)< t \}.
\end{align*}
Then, the normal vector and the horizontal normal vector of $M_{t}$ are given by
\begin{align*}
    V(x) = \nabla u(x) , \qquad V_{0}(x) = \nabla_{0} u(x).
\end{align*}
Also, the unit normal vector at $x \in M_{t}$ and the horizontal unit normal vector $\nu_{0}$ at $x \in M_{t} \setminus \Sigma (M_{t})$ are given by
\begin{align}\label{2.5}
    \nu(x)= \frac{V}{|V|} = \frac{\nabla u(x)}{|\nabla u(x)|}, \qquad \nu_{0}(x)= \frac{V_{0}}{|V_{0}|}= \frac{\nabla_{0} u(x)}{|\nabla_{0} u(x)|} .
\end{align}
Combining (\ref{2.4}) and (\ref{2.5}), for any $x \in M_{t} \setminus \Sigma (M_{t})$, the horizontal mean curvature $H_{0}(x)$ can be expressed as
\begin{align}\label{2.6}
    H_{0}(x) =\operatorname{div}_{0}\left(\frac{\nabla_{0} u(x)}{|\nabla_{0} u(x)|}\right) = \frac{1}{|\nabla_{0} u|} Tr\left[ \left(I - \frac{\nabla_{0}u \otimes \nabla_{0} u}{|\nabla_{0} u|^{2}}\right)\left(\nabla_{0}^{2} u\right)^{*}\right].
\end{align}

Notice that the foregoing geometric quantities associated with the HIMCF have been expressed in terms of the level set function $u$, so we can transform the geometric equation (\ref{1.2}) into a scalar equation for $u$. The process is carried out as follows:
\begin{align*}
    1 &=\frac{d u}{d t} =\left\langle\frac{\partial u}{\partial x}, \frac{\partial x}{\partial t}\right\rangle _{\mathbb{R}^{3}}  = \left\langle \nabla u ,\frac{\partial x}{\partial t}\right\rangle  = |\nabla u| \left\langle \nu ,\frac{\partial x}{\partial t}\right\rangle = |\nabla u| H_{0}^{-1} \left\langle \nu , \nu_{0}\right\rangle \\
    & =\left(\operatorname{div}_{0}\left(\frac{\nabla_{0} u}{|\nabla_{0} u|}\right)\right)^{-1}|\nabla_{0} u | ,
\end{align*}
where the last equality follows from (\ref{2.5}) and (\ref{2.6}). Hence
\begin{align}\label{2.7}
    \operatorname{div}_{0}\left(\frac{\nabla_{0} u}{|\nabla_{0} u|}\right) =|\nabla_{0} u|, \qquad x\in M_{t} \setminus \Sigma (M_{t}).
\end{align}

\begin{lemma}\label{lem2.7}
    Let $v(x,t)= u(x) - t$. If $u(x)$ satisfies (\ref{2.7}) for any $x \in M_{t} \setminus \Sigma (M_{t})$, then $v(x,t)$ satisfies
    \begin{align*}
        \frac{\partial v}{\partial t} = - |\nabla_{0} v| \left(\operatorname{div}_{0}\left(\frac{\nabla_{0} v}{|\nabla_{0} v|}\right)\right)^{-1}, \qquad x \in M_{t} \setminus \Sigma (M_{t}).
    \end{align*}
\end{lemma}
\noindent {\bf Proof}: It is obvious that  $v(x,0) =0$ if and only if $x \in M_{0}$, and $|\nabla_{0} v| = |\nabla_{0} u| \neq  0$ at $x \in M_{t} \setminus \Sigma (M_{t})$. By direct calculation, we have
\begin{align*}
    0 = \frac{d v}{d t} = \left\langle \frac{\partial v}{\partial x} ,\frac{\partial x}{ \partial t} \right\rangle _{\mathbb{R}^{3}} +  \frac{\partial v}{\partial t} = \left\langle \frac{\partial u}{\partial x}, \frac{\partial x}{ \partial t} \right\rangle _{\mathbb{R}^{3}} +  \frac{\partial v}{\partial t}.
\end{align*}
Thus,
\begin{align*}
    \frac{\partial v}{\partial t} = - \left\langle \frac{\partial u}{\partial x}, \frac{\partial x}{ \partial t} \right\rangle _{\mathbb{R}^{3}}  = -|\nabla u| \left\langle \nu ,\frac{\partial x}{\partial t}\right\rangle = - |\nabla_{0} v| \left(\operatorname{div}_{0}\left(\frac{\nabla_{0} v}{|\nabla_{0} v|}\right)\right)^{-1},
\end{align*}
where we use (\ref{1.2}), (\ref{2.5}) and (\ref{2.6}).
\hfill${\square}$

\subsection{ Horizontal perimeter in Heisenberg groups }
\
\vglue-10pt
 \indent

This section reviews the definition of the horizontal perimeter of sets in the Heisenberg group, which is inspired by De Giorgi's theorem (\cite{D1954,D1955}) in Euclidean space, as well as the structure of sets with finite horizontal perimeters (see, e.g., \cite{FSC2002,FSC2003,P1982}). The underlying measure $dx$ used here is the Haar measure on $\mathbb{H}^{1}$, which agrees with both the exponential of the Lebesgue measure on the Lie algebra and the Hausdorff 4-measure associated with the $CC$-metric $d$.

\begin{definition} \cite{CDPT2007}\label{def2.8}
    Let $E \subset \mathbb{H}^{1}$ be a measurable set and $\Omega \subset \mathbb{H}^{1}$ be an open set. The horizontal perimeter ($\mathbb{H}$-perimeter) of $E$ in $\Omega$ is given by
    \begin{align*}
        P_{\mathbb{H}} (E,\Omega)= \sup_{\varphi \in \mathcal{A}(\Omega) } \int_{\Omega} \chi_{E}(x) (X_{1} \varphi_{1} + X_{2}\varphi_{2})  dx,
    \end{align*}
    where $\mathcal{A}(\Omega)= \{\varphi  = (\varphi_{1},\varphi_{2})\in C^{1}_{0}(\Omega, \mathbb{R}^{2} )\text{ and } |\varphi| \leq 1\}$ and $\chi_{E}$ denotes the characteristic function of $E$. Sets with finite $\mathbb{H}$-perimeter are called $\mathbb{H}$-Caccioppoli sets. For $\Omega = \mathbb{H}^{1}$, we let $P_{\mathbb{H}} (E,\mathbb{H}^{1}) = P_{\mathbb{H}} (E)$.
\end{definition}

The above definition of the $\mathbb{H}$-perimeter is quite natural in the geometric problems of Calculus of Variations.
\begin{lemma}\cite{CDPT2007,FSC1996}\label{lem2.9}
    Let $\partial E$ be a $C^{1}$ surface in $\mathbb{H}^{1}$ which bounds an open set $E$, then
    \begin{align*}
        P_{\mathbb{H}} (E,\Omega) = \int_{\partial E \cap \Omega} \left(\sum_{i = 1}^{2} \left\langle X_{i},\nu \right\rangle^{2} \right)^{1/2} d\mathcal{H}^{2} =  \int_{\partial E \cap \Omega} \frac{|V_{0}|}{|V |} d\mathcal{H}^{2} ,
    \end{align*}
    where $\mathcal{H}^{2}$ is the Euclidean 2-dimension Hausdorff measure and $\nu \in C^{0}(\partial E )$ is the outer unit normal vector.
\end{lemma}
\begin{definition}\cite{DGN2006}\label{def2.10}
    Let $M$ be an oriented $C^{2}$ surface in $\mathbb{H}^{1}$. The $\mathbb{H}$-perimeter measure supported on $M$ is denoted by
    \begin{align}\label{2.8}
        d\sigma_{\mathbb{H}} = | V_{0} | d\sigma = \frac{|V_{0}|}{|V |} d\mathcal{H}^{2} = \left\langle \nu_{0}, \nu \right\rangle  d\mathcal{H}^{2},
    \end{align}
    where $d\sigma$ denotes the surface measure on $M$.
\end{definition}

Since the $C^{1}$ smoothness of the boundary of $E_{t}$ is not preserved under the HIMCF, the horizontal unit outer normal vector $\nu_{0}$ is not well defined at characteristic points. Therefore, following the argument of De Giorgi (see e.g. \cite{D1955, EG1992}), we can consider the generalized horizontal outer normal vector and the $\mathbb{H}$-reduced boundaries.
\begin{definition}\label{def2.11}
    Let $E \subset \mathbb{H}^{1}$ be a $\mathbb{H}$-Cacciopoli set, we call $\nu_{0}^{*}$ a generalized horizontal outer unit normal vector to $E$ provided
    \begin{align*}
        \int_{E} \nabla_{0} \varphi dx = \int_{\partial E} \left\langle  \nu_{0}^{*}, \varphi \right\rangle d\sigma_{\mathbb{H}}, \qquad \varphi \in C^{1}_{0}(\mathbb{H}^{1}).
    \end{align*}
\end{definition}
The existence of $\nu_{0}^{*}$ follows from Riesz's representation theorem. Next, we define the $\mathbb{H}$-reduced boundary.
\begin{definition}\label{def2.12}
    The $\mathbb{H}$-reduced boundary of a $\mathbb{H}$-Caccioppoli set $E \subset \mathbb{H}^{1}$ is given by
    \begin{align*}
        \partial_{\mathbb{H}}^{*}E = \left\{x \in \text{spt} \mu : \lim_{r \to 0^{+} } \frac{\int_{B_{r}(x)} \nu_{0}^{*} d\mu}{\mu (B_{r}(x))}  \text{ exists and has length } 1\right\} ,
    \end{align*}
    where $\text{spt} \mu = \left\{ x \in \mathbb{H}: \mu (B_{r}(x)) \neq 0 \text{ for each } r>0 \right\}$.
\end{definition}
Given that Definition \ref{def2.12} involves only the integral of $\nu_{0}^{*}$, the horizontal reduced boundary $\partial_{\mathbb{H}}^{*}E$ depends only on $E$. In addition, modifying a zero measure set of $E$, $\partial_{\mathbb{H}}^{*}E$ is unchanged and has the structure of a 2-dimensional surface.
\begin{lemma}\cite{CDPT2007}\label{lem2.13}
    If $\partial E$ is not a $C^{1}$ surface in $\mathbb{H}^{1}$, then
    \begin{align}\label{2.9}
        P_{\mathbb{H}} (E,\Omega)=\int_{\partial^{*}_{\mathbb{H}} E \cap \Omega} \frac{|V_{0}|}{|V|} d\mathcal{H}^{2}= \int_{\partial^{*}_{\mathbb{H}} E \cap \Omega} |V_{0} | d \sigma  = \int_{\partial^{*}_{\mathbb{H}} E \cap \Omega}  d \sigma_{\mathbb{H}} .
    \end{align}
\end{lemma}
By virtue of De Giorgi's structure theory, $\mathbb{H}$-Caccioppoli sets also have the following properties.
\begin{lemma}\cite{FSC2001}\label{lem2.14}
    If $E \subset \mathbb{H}^{1}$ is a $\mathbb{H}$-Caccioppoli set, then the reduced boundary $\partial_{\mathbb{H}}^{*}E$ is rectifiable.
\end{lemma}
\begin{lemma}\cite{M2012,MC2001}\label{lem2.15}
    Assume that $E, F \subset \mathbb{H}^{1}$ are sets with locally finite $\mathbb{H}$-perimeter, and $\sigma_{\mathbb{H}} \llcorner  \partial^{*}_{\mathbb{H}} E$ is a Radon measure on $\mathbb{H}^{1}$. Then
    \begin{align*}
        \sigma_{\mathbb{H}} \llcorner  \partial^{*}_{\mathbb{H}} (E \cup F) + \sigma_{\mathbb{H}} \llcorner  \partial^{*}_{\mathbb{H}} (E \cap F) \leq \sigma_{\mathbb{H}} \llcorner  \partial^{*}_{\mathbb{H}} E + \sigma_{\mathbb{H}} \llcorner  \partial^{*}_{\mathbb{H}} F .
    \end{align*}
    In particular, if $E \cap F = \emptyset$, then
    \begin{align*}
        \sigma_{\mathbb{H}} \llcorner  \partial^{*}_{\mathbb{H}} (E \cup F)    =   \sigma_{\mathbb{H}} \llcorner  \partial^{*}_{\mathbb{H}} (E) + \sigma_{\mathbb{H}} \llcorner  \partial^{*}_{\mathbb{H}} ( F) .
    \end{align*}
\end{lemma}
\begin{lemma}\cite{DCP1972, G1984}\label{lem2.16}
    The $\mathbb{H}$-perimeter is lower semi-continuous in the sense of $L^{1}_{loc}$ convergence.
\end{lemma}

\begin{lemma}\cite{MC2001}\label{lem2.17}
    Let $E \subset \mathbb{H}^{1}$ be a $\mathbb{H}$-Cacciopoli set, and for any $\lambda >0$,      $\delta_{\lambda}(E) = \{ \delta_{\lambda}(x), x\in E\}$. Then,
    \begin{align*}
        P_{\mathbb{H}}(\delta_{\lambda}(E)) = \lambda^{3} P_{\mathbb{H}}(E),\qquad  |\delta_{\lambda}(E)| = \lambda^{4}|E|.
    \end{align*}
\end{lemma}

\begin{lemma}(Coarea formula)\cite{M2011}\label{lem2.18}
    Let $u : \mathbb{H}^{1} \to \mathbb{R}$ be a Lipschitz mapping. Then
    \begin{align}\label{2.10}
        \int_{\mathbb{H}} \varphi (x) |\nabla_{0} u| dx = \int_{\mathbb{R}} \left(\int_{\{x \in \mathbb{H}^{1}, u(x) = t\}} \varphi (x) d\sigma_{\mathbb{H}} \right) dt ,
    \end{align}
    where $\varphi(x) : \mathbb{H}^{1} \to \mathbb{R}$  is any nonnegative measurable function.
\end{lemma}
Let $E$ be a bounded region enclosed by a $C^{2}$ oriented, immersed surface $M = \partial E$. Consider a $C^{1}$ vector field $\mathcal{X}$ with compact support on $M$, and let $\{\Psi_{t}\}_{t \in \mathbb{R}}$ denote the associated group of diffeomorphisms. For small values of $t$, let
\begin{align}\label{2.11}
    M_{t} := \Psi_{t}(M) = M + t\mathcal{X}
\end{align}
be the variation of $M$ induced by $\mathcal{X}$, and let $E_{t}$ denote the region enclosed by $M_{t}$. Next, we consider the first variation of the horizontal perimeter under the deformation (\ref{2.11}).
\begin{lemma}\cite{CDPT2007}\label{lem2.19}
    Let $M \subset \mathbb{H}^{1}$ be an oriented $C^{2}$ immersed surface, and $\mathcal{X}$ be a $C^{1}$ vector field with compact support on $M$. Then
    \begin{align}\label{2.12}
        \frac{d}{dt}P_{\mathbb{H}}(E_{t})\Big|_{t = 0} =  \int_{M \setminus \Sigma(M) } H_{0}  \left\langle \mathcal{X}, \nu \right\rangle d\mathcal{H}^{2} -  \int_{M \setminus \Sigma(M)} \operatorname{div}_{M}\left(\left\langle \mathcal{X}, \nu \right\rangle (\nu_{0})^{\top}\right) d\mathcal{H}^{2},
    \end{align}
    where $\nu_{0}^{\top}$ denotes the tangential component of $\nu_{0}$. Moreover, if $H_{0} \in \mathcal{L}^{1}(M, d\mathcal{H}^{2})$, then
    \begin{align}\label{2.13}
        \frac{d}{dt}P_{\mathbb{H}}(E_{t})\Big|_{t = 0} =  \int_{M } H_{0} \left\langle \mathcal{X}, \nu \right\rangle  d\mathcal{H}^{2} -  \int_{M } \operatorname{div}_{M}\left(\left\langle \mathcal{X}, \nu \right\rangle (\nu_{0})^{\top}\right) d\mathcal{H}^{2}.
    \end{align}
\end{lemma}
\begin{remark}\label{rem2.20}
    From the proof of the first variation formula in Lemma \ref{lem2.19} (see \cite[p.~134]{CDPT2007}), it follows that
    \begin{align}\label{2.14}
        \int_{M} H_{0}\left\langle \mathcal{X}, \nu \right\rangle   d\mathcal{H}^{2} -& \int_{M} \operatorname{div}_{M} (\left\langle \mathcal{X}, \nu \right\rangle (\nu_{0})^{\top}) d\mathcal{H}^{2} \\
        & =\int_{M} \left[\mathcal{X}^{\bot} (|V_{\mathcal{H}}|) + |V_{\mathcal{H}}| \operatorname{div}_{M} \mathcal{X}^{\bot}  \right] d\mathcal{H}^{2} \notag ,
    \end{align}
    where $V_{\mathcal{H}} = \nu - \left\langle  \nu, T\right\rangle T$.
\end{remark}

The well-known result on the integrability of $H_{0}$ with respect to the $\mathbb{H}$-perimeter measure $d\sigma_{\mathbb{H}}$ is recalled here.
\begin{lemma}\label{lem2.21} \cite{DGN2012}
    Let $M \subset \mathbb{H}^{1}$ be an oriented $C^{2}$ immersed surface, and $\mathcal{X}$ be a $C^{1}$ vector field with compact support on $M$. Then, $H_{0} \in L^{1}(M, d\sigma_{\mathbb{H}})$.
\end{lemma}
This result is crucial, as it allows one to consider horizontal variations satisfying $\left\langle \mathcal{X}, \nu \right\rangle = F \left\langle \nu_{0}, \nu \right\rangle$, where $F \in C^{1}_{0}(M)$. Clearly, the horizontal variations are compactly supported in $M \setminus \Sigma(M)$. Using (\ref{2.8}), we obtain
\begin{lemma}\cite{CDPT2007}\label{lem2.22}
    Let $M \subset \mathbb{H}^{1}$ be a $C^{2}$ immersed surface. Assume that $\mathcal{X}$ is a horizontal variation, i.e.,  $\left\langle \mathcal{X}, \nu \right\rangle  = F \left\langle \nu_{0}, \nu \right\rangle$, where $F\in C^{1}_{0}(M)$. Then, the equation (\ref{2.12}) can be reduced to
    \begin{align}\label{2.15}
        \frac{d}{dt}P_{\mathbb{H}}(E_{t})\Big|_{t = 0} =  \int_{M } \left\langle \mathcal{X}, \nu \right\rangle H_{0} d\mathcal{H}^{2} =  \int_{M } F H_{0} d\sigma_{\mathbb{H}} .
    \end{align}
\end{lemma}

\section{ Notions of solutions}

\noindent Note that even though $u$ is smooth, the associated level set $E_{t}$ may contain singularities and characteristic points. Therefore, it is necessary to introduce a notion of weak solution, which is a generalization of the weak solution introduced in \cite{HI2001} to the Heisenberg group. The weak solution of the HIMCF can be defined as a minimizer of either functional (\ref{1.5}) or functional (\ref{1.6}). We also give an explicit example of then HIMCF when $E_{0}$ is a unit Kor$\acute{\rm a}$nyi ball.

\subsection{ Weak solutions}
\
\vglue-10pt
 \indent

Let $\Omega \subset \mathbb{H}^{1} $ be an open set, $K$ be a compact set with $K \subset \Omega$, and $v$ be a locally Lipschitz function. We consider the functional
\begin{align*}
    J_{u}^{K}(v):=\int_{K} |\nabla_{0} v|+ v|\nabla_{0} u| dx.
\end{align*}
In general, $K$ can be omitted, since it does not matter which such set is chosen.

\begin{lemma}\label{lem3.1}
    $u$ is a minimizer of the functional $J_{u}^{K}(v)$ if and only if $u$ satisfies the evolution equation (\ref{2.7}).
\end{lemma}
\noindent {\bf Proof}: If $x \in \Sigma (M_{t})$, we have  $| \nabla_{0} u(x)| = 0$ and $J_{u}(u) = 0$, thus $J_{u}(u) \leq J_{u}(v)$ for any $v \in HW^{1,1}(K)$. Suppose that $u$ satisfies (\ref{2.7}) for any $x \in M_{t} \setminus \Sigma (M_{t})$. Then
\begin{align*}
    \int_{K} | \nabla_{0} u |  \varphi = -  \int_{K} \left\langle  \frac{\nabla_{0} u}{|\nabla_{0} u|}, \nabla_{0} \varphi  \right\rangle , \qquad \text{for any } \varphi \in HW_{0}^{1,1}(K).
\end{align*}
Let $\varphi  = u-v$, we obtain
\begin{align*}
    \int_{K} ( u-v) | \nabla_{0} u | &= \int_{K} \left\langle  \frac{\nabla_{0} u}{| \nabla_{0} u |}, \nabla_{0} (v-u) \right\rangle  = \int_{K}  | \nabla_{0} u |^{-1} \left\langle \nabla_{0}u , \nabla_{0}v \right\rangle  -  \int_{K}  | \nabla_{0} u | \\
    & \leq \int_{K}  | \nabla_{0} v | - \int_{K}  | \nabla_{0} u | .
\end{align*}
Thus
\begin{align*}
    \int_{K}  | \nabla_{0} u | + u |\nabla_{0} u| \leq \int_{K}  | \nabla_{0} v| + v  |\nabla_{0} u|.
\end{align*}

Conversely, suppose that $u$ is a minimizer of the functional $J_{u}^{K}(v)$ and $|\nabla_{0} u| \neq 0$, we then show that $u$ satisfies (\ref{2.7}). Let
\begin{align*}
    \mathcal{F}(t) : =J_{u}^{K}(u+ t \varphi), \qquad f(x,t) :=  | \nabla_{0} (u + t \varphi) | + (u + t \varphi) | \nabla_{0} u |.
\end{align*}
where $\varphi \in C_{0}^{\infty}(\Omega)$, $u+t \varphi \in HW^{1.1}_{loc}(\Omega)$ and $t \in [0,\infty)$. Then, $\mathcal{F}(0) = J_{u}^{K}(u)$ and $\mathcal{F}(0) \leq \mathcal{F}(t)$ holds for any $t$.  Since $\mathcal{F}(t)$ attains a minimum at $t=0$, we have
\begin{align}\label{3.1}
    0 = \mathcal{F}^{'}(0) = \lim_{t \to 0} \frac{\mathcal{F}(t)-\mathcal{F}(0)}{t} =  \lim_{t \to 0} \int_{K} \frac{f(x,t) - f(x,0)}{t} dx.
\end{align}
From the mean value theorem, there exists a constant $\varepsilon \in (0,1)$ such that
\begin{align*}
    \frac{f(x,t) - f(x,0)}{t}  &= \frac{\partial}{\partial t}f(x, \varepsilon t) = \varphi \frac{\partial f (x, \varepsilon t)}{\partial (u+t\varphi)}  + \left\langle \frac{\partial f (x, \varepsilon t)}{\partial (\nabla_{0} u+t\nabla_{0} \varphi)},  \nabla_{0} \varphi \right\rangle\\
    &= \varphi | \nabla_{0} u|  + \left\langle \frac{\nabla_{0} u  +\varepsilon t\nabla_{0} \varphi }{|\nabla_{0} u  +\varepsilon t\nabla_{0}  \varphi |}, \nabla_{0} \varphi \right\rangle .
\end{align*}
Hence, since $\varphi \in C^{\infty}_{0}(\Omega)$ and $K$ is compact, we obtain
\begin{align*}
    \frac{1}{t}  \left| f(x,t)-f(x,0) \right| \leq  \left\| \varphi \right\|_{L^{\infty}(K)} |\nabla_{0} u| + \left\|\nabla_{0} \varphi \right\|_{L^{\infty}(K)}  \in L^{1}(K) .
\end{align*}
By Lebesgue's dominated convergence theorem, we obtain
\begin{align}\label{3.2}
    \lim_{t \to 0} \int_{K} \frac{f(x,t) - f(x,0)}{t} &= \int_{K} \lim_{t \to 0} \frac{f(x,t) - f(x,0)}{t} \notag \\
    &=\int_{K} |\nabla_{0} u|\varphi + \frac{1}{|\nabla_{0} u|} \left\langle \nabla_{0} u, \nabla_{0} \varphi \right\rangle .
\end{align}
Combining (\ref{3.1}) with (\ref{3.2}), we get
\begin{align*}
    0= \int_{K} |\nabla_{0} u|\varphi + |\nabla_{0} u|^{-1}\left\langle \nabla_{0} u, \nabla_{0} \varphi \right\rangle , \qquad \text{for any } \varphi \in C_{0}^{\infty}(\Omega),
\end{align*}
which yields that $u(x)$ satisfies (\ref{2.7}). This completes the proof of Lemma \ref{lem3.1}.
\hfill${\square}$

In the following, we state the definition of the weak solution to (\ref{1.2}).

\begin{definition}\label{def3.2}
   A locally Lipschitz function $u : \Omega \to \mathbb{R}$ is called a weak subsolution (supersolution) of (\ref{1.2}) in $\Omega$ if, for any locally Lipschitz function $v$, where $v\leq u$ (or $v\geq u$) and $\{ v \neq u \} \subset \subset \Omega$, the following inequality holds
    \begin{align}\label{3.3}
        J_{u}^{K}(u) \leq J_{u}^{K}(v),
    \end{align}
    where $K$ is a compact set and $\{ v \neq u \} \subset K$. A function $u$ is called a weak solution of (\ref{1.2}) if it is both a subsolution and a supersolution.
\end{definition}
Given a bounded open set $E_{0} =  \{x \in \mathbb{H}^{1} : u(x) < 0\}$ with a smooth boundary $M_{0} =\partial E_{0}$. We say that $u$ is a weak solution of (\ref{1.2})  with initial condition $M_{0}$ if
\begin{equation}\label{3.4}
    \begin{cases}
        J_{u}^{K}(u) \leq J_{u}^{K}(v),  \qquad K \subset  \mathbb{H}^{1} \setminus E_{0}, \\
        u|_{M_{0}} = 0.
    \end{cases}
\end{equation}

\subsection{ An equivalent definition }
\
\vglue-10pt
 \indent

Next, the minimization principle (\ref{3.3}) for $u$ is extended to the corresponding level set $E_{t} = \{u < t\}$. The symmetric difference of two sets $E$ and $F$ is defined as $E \Delta F := (E \setminus F) \cup (F \setminus E)$. Given a measurable set $F \subset \Omega$ with locally finite $\mathbb{H}$-perimeter, we define the functional of the level set $F$ by
\begin{align*}
    J_{u}^{K}(F)=P_{\mathbb{H}}(F,K) - \int_{F \cap K} |\nabla_{0} u|dx.
\end{align*}
From Lemma \ref{lem2.16}, $J_{u}^{K}(F)$ is lower semi-continuous.

We say that $E$ minimizes $J_{u}^{K}(F)$ in $\Omega$ on the outside (inside) if
\begin{align}\label{3.5}
     J_{u}^{K}(E) \leq  J_{u}^{K}(F)
\end{align}
for any set $F$ satisfying $F\supseteq E$ ($F\subseteq E$) and  $E \Delta F \subseteq K \subset \subset \Omega$. Using Lemma \ref{lem2.15}, we obtain
\begin{align}\label{3.6}
     J_{u}(E \cup F) +J_{u}(E \cap F) \leq J_{u}(E) + J_{u}(F).
\end{align}
Thus, $E$ minimizes $J_{u}^{K}(F)$ in $\Omega$ if and only if $E$ minimizes $J_{u}^{K}(F)$ in $\Omega$ on the outside and on the inside.
\begin{lemma}\label{lem3.3}
    For each $t>0$, $E_{t} = \{ u(x) < t\}$ minimizes the functional (\ref{1.6}) in $\Omega$ (on the outside, inside respectively) if and only if the locally Lipschitz function $u: \Omega \to \mathbb{R}$ is a weak solution (subsolution, supersolution respectively) of (\ref{1.2}) in $\Omega$.
\end{lemma}
\noindent {\bf Proof}: Let $F_{t}=\{v < t\}$, and $v$ be a locally Lipschitz function. Let $\{v\neq u\} \subset \subset \Omega$ and $\{v \neq u\}$ contained in a compact set $K$. It follows that $\{F_{t} \Delta E_{t}\}_{t>0} \subset \subset \Omega$ and $\{F_{t} \Delta E_{t}\}_{t>0} \subseteq K$. The proof is divided into two parts.

First we show that if $E_{t}$ minimizes $J_{u}(F_{t})$ in $\Omega$, then $u$ is a weak solution of (\ref{1.2}) on $\Omega$. Denote $a := \min_{x\in K} u(x)$ and $b := \max_{x\in K} v(x)$, by Lemma \ref{lem2.18}, we have
\begin{align}\label{3.7}
    \int_{a}^{b} J_{u}(F_{t}) dt &=  \int_{a}^{b} P_{\mathbb{H}}(F_{t}, K) dt - \int_{a}^{b} \int_{F_{t} \cap K} |\nabla_{0} u| dx dt  \notag \\
    &=\int_{K} |\nabla_{0} v| dx - \int_{K} \left(\int_{v(x)}^{b} |\nabla_{0} u| dt\right) dx \notag \\
    &= J_{u}^{K}(v) - b\int_{K} |\nabla_{0} u|  dx.
\end{align}
Similarly, we have
\begin{align}\label{3.8}
    \int_{a}^{b} J_{u}(E_{t}) dt = J_{u}^{K}(u) - b\int_{K} |\nabla_{0} u|  dx.
\end{align}
Since $J_{u}(E_{t}) \leq J_{u}(F_{t})$, we deduce $J_{u}^{K}(u) \leq J_{u}^{K}(v)$ from (\ref{3.7}) and (\ref{3.8}).

Next we prove that if $u$ is a weak solution of (\ref{1.2}) in $\Omega$, then $E_{t}$ minimizes $J_{u}(F_{t})$ in $\Omega$. The details are provided only for weak subsolutions, as the case of weak supersolutions follows similarly. Fix $t_{0}$ and let $F$ be a set such that
\begin{align*}
    E_{t_{0}} \subseteq F , \qquad F \setminus E_{t_{0}} \subset \subset \Omega.
\end{align*}
By the lower semi-continuity of the functional (\ref{1.6}), we may assume that for all $G$ with $ G \Delta E_{t_{0}} \subseteq  F \Delta E_{t_{0}}$, there is $J_{u}(F) \leq  J_{u}(G)$. Let $G= E_{t} \cap F$ and
\begin{equation*}
    F_{t} :=
    \begin{cases}
        F \cup E_{t}, \qquad &t\leq t_{0}, \\
        E_{t},        \qquad &t > t_{0}.
    \end{cases}
\end{equation*}
From (\ref{3.6}), we obtain
\begin{align}\label{3.9}
    J_{u}(F \cap E_{t}) + J_{u}(F \cup E_{t}) \leq J_{u}(E_{t}) +J_{u}(F).
\end{align}
Thus
\begin{align}\label{3.10}
    J_{u}(F_{t}) \leq J_{u}(E_{t}), \qquad \text{for all } t.
\end{align}
By virtue of $F_{t}=\{v<t\}$, we can define $v$
\begin{equation*}
    v :=
    \begin{cases}
        t_{0}, \qquad &\text{ on }  F \setminus E_{t_{0}}, \\
        u,    \qquad &\text{ elsewhere}.
    \end{cases}
\end{equation*}
Note that $v \in BV_{loc}\cap L^{\infty}_{loc}$ and $\{v\neq u\} \subset \subset \Omega$, so $J_{u}(v)$ make sense. From Remark \ref{rem2.2}, there exists a sequence functions $\{v_{k}\}_{k\in \mathbb{N}} $ such that $\left\|v_{k} - v\right\|_{HW^{1,1}(\Omega)} \to 0$. Thus, $\lim_{k \to \infty} \int_{K} |\nabla_{0}v_{k} - \nabla_{0}v| dx =0$. Since $J_{u}^{K}(u) \leq J_{u}^{K}(v_{k}) $ holds for $k \in\mathbb{N}$, we have $J_{u}^{K}(u) \leq J_{u}^{K}(v) $. Combining (\ref{3.7}) and (\ref{3.8}), we have
\begin{align*}
    \int_{a}^{b} J_{u}(E_{t}) dt \leq \int_{a}^{b} J_{u}(F_{t}) dt.
\end{align*}
From (\ref{3.10}), we obtain $J_{u}(E_{t}) =J_{u}(F_{t})$ for a.e. $t$. In view of the definition of $F_{t}$, we have
\begin{align*}
    J_{u}(E_{t}) = J_{u}(F \cup E_{t}),\qquad \text{ for a.e. } t\leq t_{0}.
\end{align*}
Also, by (\ref{3.9}), we get
\begin{align*}
    J_{u}(F)  \geq J_{u}(F \cap E_{t}) , \qquad \text{ for a.e. } t\leq t_{0}.
\end{align*}
From the lower semi-continuity of the functional (\ref{1.6}), we have
\begin{align*}
    J_{u}(F) \geq  \liminf_{t \to t_{0}} J_{u}(F \cap E_{t}) \geq J_{u}(E_{t_{0}})  .
\end{align*}
Hence, $E_{t_{0}}$ minimizes $J_{u}(F)$ on the outside.
\hfill${\square}$

Let $\{E_{t}\}_{t>0}$ be a family of open sets that is closed under ascending unions. The family $\{E_{t}\}_{t\geq 0}$ is called a weak solution of (\ref{1.2}) with initial condition $E_{0}$ if
\begin{equation}\label{3.11}
    \begin{cases}
        J_{u}(E_{t}) \leq J_{u}(F) , \qquad &F\Delta E_{t} \subset \subset \mathbb{H}^{1} \setminus E_{0}  \text{ for each } t>0 ,\\
        E_{0},  \qquad &t=0.
    \end{cases}
\end{equation}
Furthermore, by approximating the boundary and applying Lemma \ref{lem3.3}, we can obtain the following conclusion.

\begin{lemma}\label{lem3.4}
    The initial value problem (\ref{3.4}) is equivalent to (\ref{3.11}).
\end{lemma}
Next we give the definition of generalized horizontal inverse mean curvature flow.
\begin{definition}\label{def3.5}
    Let $M_{0} =\partial \{x \in \mathbb{H}^{1} : u(x) < 0 \} = \partial E_{0}$ be a smooth surface in $\mathbb{H}^{1}$. We say that $M_{t} = \partial \{x \in \mathbb{H}^{1} : u(x) < t \} = \partial E_{t}$ is a generalized horizontal inverse mean curvature flow with initial surface $M_{0}$ if either $u$ is a locally Lipschitz function satisfying (\ref{3.4}), or $E_{t}$ is precompact and satisfies (\ref{3.11}).
\end{definition}

\subsection{ An explicit solution }
\
\vglue-10pt
 \indent

\begin{proposition}\label{prop3.6}
    If $E_{0}$ is a unit Kor$\acute{a}$nyi ball with center $o = (0,0,0)$, then under the HIMCF, the evolving surface $M_{t} = \partial E_{t}$ is given by $M_{t} = \left\{x \in \mathbb{H}^{1}, \|x\|_{\rho} = e^{\frac{1}{3}t}\right\}$, where $\|x\|_{\rho}=\left[\left(x_{1}^{2}+x_{2}^{2}\right)^{2}+16 x_{3}^{2}\right]^{1/4}$.
\end{proposition}
\noindent {\bf Proof}: According to Section 2.2, if $u$ satisfies (\ref{2.7}), then $M_{t} = \{ u = t \}$ evolves under the HIMCF. To simplify the calculations, let $f(x):=\|x\|_{\rho}^{4}$. By direct calculation, we obtain
\begin{align*}
    X_{1} u &= \frac{3}{4} f^{-1} X_{1} f , &
    X_{2} u &= \frac{3}{4} f ^{-1} X_{2} f ,\\
    X^{2}_{1} u &= -  \frac{3}{4} f ^{-2} X_{1} f X_{1} f + \frac{3}{4} f ^{-1} X^{2}_{1} f , &
    X^{2}_{1} u &= -  \frac{3}{4} f ^{-2} X_{2} f X_{2} f + \frac{3}{4} f ^{-1} X^{2}_{2} f, \\
    X_{1} X_{2}u &= -  \frac{3}{4} f ^{-2} X_{1} f X_{2} f + \frac{3}{4} f ^{-1} X_{1} X_{2} f , &
    X_{2} X_{1}u &= -  \frac{3}{4} f ^{-2} X_{2} f X_{1} f + \frac{3}{4} f ^{-1} X_{2} X_{1} f,
\end{align*}
where
\begin{align*}
    X_{1} f &= 4x_{1}\left(x_{1}^{2} + x_{2}^{2}\right) - 16x_{2}x_{3} , &
    X_{2} f &= 4x_{2}\left(x_{1}^{2} + x_{2}^{2}\right) + 16x_{1}x_{3} ,\\
    X^{2}_{1} f &= 12 x_{1}^{2} + 12 x_{2}^{2} , & X^{2}_{1} f &= 12 x_{1}^{2} + 12 x_{2}^{2}, \\
    X_{1} X_{2} f &= 16x_{3}, & X_{2} X_{1} f &= - 16x_{3}.
\end{align*}
Furthermore
\begin{align*}
    \left(\nabla^{2}_{0}u\right)^{*} = \left({\begin{array}{cc}
        -  \frac{3}{4} f ^{-2} X_{1} f X_{1} f &  -  \frac{3}{4} f ^{-2} X_{1} f X_{2} f\\
        -  \frac{3}{4} f ^{-2} X_{2} f X_{1} f &  -  \frac{3}{4} f ^{-2} X_{2} f X_{2} f \\
    \end{array}} \right)
    +
    \left({\begin{array}{cc}
        \frac{3}{4} f ^{-1} X^{2}_{1} f &  0 \\
        0 & \frac{3}{4} f ^{-1} X^{2}_{2} f \\
    \end{array}} \right)
    := A_{1} + A_{2} .
\end{align*}
The above equalities imply that
\begin{align*}
    Tr &\left[ \left(I - \frac{\nabla_{0} u \otimes \nabla_{0} u }{|\nabla_{0} u|^{2}}\right) \left(\nabla_{0}^{2} u\right)^{*}\right] = Tr\left[ \left(I - \frac{\nabla_{0} u \otimes \nabla_{0} u }{|\nabla_{0} u|^{2}}\right) (A_{1}+ A_{2})\right] \\
    &= Tr\left[ \left(I - \frac{\nabla_{0} u \otimes \nabla_{0} u }{|\nabla_{0} u|^{2}}\right) A_{1}\right] + Tr\left[ \left(I - \frac{\nabla_{0} u \otimes \nabla_{0} u }{|\nabla_{0} u|^{2}}\right) A_{2}\right] .
\end{align*}
Note that
\begin{align*}
    \left(I - \frac{\nabla_{0} u \otimes \nabla_{0} u }{|\nabla_{0} u|^{2}}\right) A_{1} =  \left({\begin{array}{cc}
        0 &  0 \\
        0 & 0 \\
    \end{array}} \right)
\end{align*}
and
\begin{align*}
    A_{2} =    \frac{3}{4} f^{-1} \left({\begin{array}{cc}
        12 x_{1}^{2} + 12 x_{2}^{2} &  0 \\
        0 & 12 x_{1}^{2} + 12 x_{2}^{2} \\
    \end{array}} \right).
\end{align*}
Hence,
\begin{align}\label{3.12}
    Tr\left[ \left(I - \frac{\nabla_{0} u \otimes \nabla_{0} u }{|\nabla_{0} u|^{2}}\right) \left(\nabla_{0}^{2} u\right)^{*}\right] =  9 f^{-1} \left(x_{1}^{2} + x_{2}^{2}\right)
\end{align}
and
\begin{align}\label{3.13}
    |\nabla_{0} u|^{2} &= (X_{1} u)^{2} + (X_{2} u)^{2}  = \frac{9}{16} f^{-2} \left[(X_{1} f )^{2} + (X_{2} f)^{2}\right] = 9 f^{-1} \left(x_{1}^{2} +  x_{2}^{2} \right) .
\end{align}
Combining (\ref{3.12}) and (\ref{3.13}), we have
\begin{align*}
    |\nabla_{0} u| =  |\nabla_{0} u|^{-1}Tr\left[ \left(I - \frac{\nabla_{0} u \otimes \nabla_{0} u }{|\nabla_{0} u|^{2}}\right) \left(\nabla_{0}^{2} u\right)^{*}\right]  = \operatorname{div}_{0}\left(\frac{\nabla_{0} u}{|\nabla_{0} u|}\right) .
\end{align*}
Moreover, $\Sigma(M_{t}) = \left\{x \in \mathbb{H}^{1} :|\nabla_{0} u(x)| = 0\right\} =\left\{\left(0,0,\frac{1}{4}e^{\frac{2}{3} t} \right), \left(0,0,-\frac{1}{4}e^{\frac{2}{3} t} \right)\right\}$.
\hfill${\square}$

\begin{remark}\label{rem3.7}
    Notice that $M_{t}$ cannot be generated by performing a homothetic deformation on $M_{0}= \left\{x \in \mathbb{H}^{1}, \left(x_{1}^{2}+x_{2}^{2}\right)^{2}+16 x_{3}^{2} = 1\right\} $.
\end{remark}

\section{ Behavior of singularity }

\noindent In this section, we introduce the generalized horizontal mean curvature and the horizontal minimizing hull. Additionally, we describe the mechanism of horizontal jumps at singularities under the HIMCF and establish the properties of the horizontal minimizing hull.

The generalized horizontal mean curvature is introduced to interpret weak solutions from a PDE perspective. Let $E$ be a $\mathbb{H}$-Caccioppoli set whose boundary $M = \partial E$ is a $C^2$ surface such that $\Sigma(M)$ is nowhere dense in $M$. The generalized horizontal mean curvature is defined by the first variation formula (\ref{2.14}).
\begin{definition}\label{def4.1}
    A function $H_{0} \in L^{1}_{loc}(M, d\mathcal{H}^{2})$ is called the generalized horizontal mean curvature of $M$ if for any $\mathcal{X} \in C^{1}_{0}(M)$
    \begin{align*}
        \int_{M} H_{0}\left\langle \mathcal{X}, \nu \right\rangle  d\mathcal{H}^{2}
        =\int_{M} \left[\mathcal{X}^{\bot} (|V_{\mathcal{H}}|) + |V_{\mathcal{H}}| \operatorname{div}_{M} \mathcal{X}^{\bot} + \operatorname{div}_{M} (\left\langle \mathcal{X}, \nu \right\rangle (\nu_{0})^{\top}) \right] d\mathcal{H}^{2} ,
    \end{align*}
    where $\mathcal{X}^{\bot}$ and $\nu_0^{\top}$ denote the normal and tangential components of $\mathcal{X}$ and $\nu_0$, respectively.
\end{definition}
\begin{lemma}\label{lem4.2}
    If $E_{t} = \{u(x) < t\}$ is the weak solution of (\ref{1.2}) and  $H_{0}$ is the generalized horizontal mean curvature of $M_{t}$, then for a.e. $x \in M_{t} $,
    \begin{align}\label{4.1}
        H_{0} = | \nabla_{0} u|, \qquad \text{ for a.e.} t .
    \end{align}
\end{lemma}
\noindent {\bf Proof}:
From Lemma \ref{lem3.3}, $u$ is a minimizer of the functional (\ref{1.5}). Thus,
\begin{align*}
    0&= \frac{d}{dt}\Big|_{t=0} J_{u}\left(u(\Psi_{t}(x))\right)\\
    & = \frac{d}{dt}\Big|_{t=0}  \left(\int_{-\infty}^{+\infty} P_{\mathbb{H}}(E_{t}, K) dt + \int_{E_{t}\cap K} u(\Psi_{t}(x)) |\nabla_{0} u| dx \right),
\end{align*}
where $\Psi_{t}(x) = x +t \mathcal{X}$ and $\mathcal{X}$ is a horizontal variation. The coarea formula (\ref{2.10}) is applied in the second equality. By Lebesgue's dominated convergence theorem and Lemma \ref{lem2.22}, we have
\begin{align}\label{4.2}
    0&= \int_{-\infty}^{+\infty} \frac{d}{dt} P_{\mathbb{H}}(E_{t}, K)\Big|_{t=0} dt + \int_{E_{t}\cap K} \frac{d}{dt}\Big|_{t=0} u(\Psi_{s}(x)) |\nabla_{0} u| dx \notag \\
    &= \int_{-\infty}^{+\infty}\left(\int_{M_{t}} H_{0} \frac{\left\langle  \mathcal{X} , \nu \right\rangle}{\left\langle  \nu_{0} , \nu \right\rangle} d\sigma_{\mathbb{H}}  +
    \int_{M_{t}}  \left\langle \mathcal{X}, \nabla u \right\rangle  d\sigma_{\mathbb{H}} \right)dt \notag \\
    & =  \int_{-\infty}^{+\infty} \left(\int_{M_{t}} \left\langle \mathcal{X}, \nabla u + \frac{H_{0} \nu }{\left\langle  \nu_{0} , \nu \right\rangle}\right\rangle  d\sigma_{\mathbb{H}} \right) dt ,
\end{align}
where $\nu$ and $\nu_{0}$ are the unit normal vector and the horizontal unit normal vector of $M_{t}$ respectively. Applying the Lebesgue differentiation theorem to (\ref{4.2}), we obtain (\ref{4.1}).
\hfill${\square}$

For the generalized horizontal inverse mean curvature flow, $H_0$ is allowed to approach zero, in which case the velocity function blows up to infinity. Following the approach in \cite{HI1997, HI2001}, we introduce a framework that guarantees the flow exists for all time and that the $\mathbb{H}$-perimeter of $E_t$ evolves continuously, whereby the surface performs horizontal jumps. Next, we introduce the horizontal minimizing hull which is important for understanding horizontal jumps.
\begin{definition}\label{def4.3}
    Let $\Omega \subset \mathbb{H}^{1}$ be an open set. A set $E$ is a horizontal minimizing hull if
    \begin{align}\label{4.3}
        P_{\mathbb{H}} (E, K) \leq P_{\mathbb{H}}(F, K)
    \end{align}
    for any set $F$ such that $E \subseteq F$ and $F \setminus E \subset \subset \Omega$, where $K$ is a compact set and $F\setminus E \subseteq K$.  $E$ is a strictly horizontal minimizing hull if $P_{\mathbb{H}} (E, K) = P_{\mathbb{H}}(F, K)$ implies $E \cap \Omega = F \cap \Omega$ almost everywhere. Moreover,
    the boundary of the horizontal minimizing hull is called horizontal outward minimizing.
\end{definition}
\begin{lemma}\label{lem4.4}
    Suppose $E$ is a horizontal minimizing hull and $M=\partial E$ is $C^{2}$, then the generalized horizontal mean curvature $H_{0}\big|_{M} \geq 0$ in $\mathcal{H}^{2}$- almost everywhere.
\end{lemma}
\noindent {\bf Proof}: Let $M_{t}$ be a variation induced by a vector field $\mathcal{X} \in C_{0}^{1}(M)$, and let $E_t$ denote the region enclosed by $M_t$, with $E \subset E_{t}$ for $t>0$. From the inequality (\ref{4.3}), we have
\begin{align*}
    \frac{d}{dt}P_{\mathbb{H}}(E_{t})\Big|_{t = 0} \geq 0.
\end{align*}
If $\mathcal{X}$ is a horizontal variation, it follows from (\ref{2.15}) that
\begin{align*}
    \frac{d}{dt}P_{\mathbb{H}}(E_{t})\Big|_{t=0} =  \int_{M} \left\langle \mathcal{X},\nu  \right\rangle H_{0}  d\mathcal{H}^{2} \geq 0 .
\end{align*}
This completes the proof.
\hfill${\square}$
\begin{remark}\label{rem4.5}
    Let $\Omega \subset \mathbb{H}^{1}$ be an open set, the horizontal minimizing hull $E \subset \mathbb{H}^{1}$ (in $\Omega$) is a locally subsolution to the plateau problem in $\mathbb{H}^{1}$ which stated in \cite{P2004}.
\end{remark}
In the following we state some properties of the horizontal minimizing hull.
\begin{proposition}\label{prop4.6}
    Let $\Omega \subset \mathbb{H}^{1}$ be an open set and $E \subset \mathbb{H}^{1}$ be a locally $\mathbb{H}$-Cacciopoli set in $\Omega$. Suppose there are two countable families of locally $\mathbb{H}$-Cacciopoli sets $\{E_{i}\}$ and open sets $\{\Omega_{i}\}$, such that for every $i$, $E_{i}$ is a horizontal minimizing hull (in $\Omega_{i}$), then $E=\bigcap_{i} E_{i}$ is a horizontal minimizing hull in every open sets $\Omega \subset \bigcap_{i} \Omega_{i}$.
\end{proposition}
\noindent {\bf Proof}: Let
\begin{align*}
    \mathcal{E} (E,\Omega) = P_{\mathbb{H}}(E,\Omega) - \inf \{P_{\mathbb{H}}(F,\Omega): E\Delta F \subset\subset \Omega, E \subseteq  F \} .
\end{align*}
Clearly, $\mathcal{E} (E,\Omega) \geq 0$ and equality holds if and only if $E$ is a horizontal minimizing hull in $\Omega$. From Lemma \ref{lem2.16}, we know that $\mathcal{E} (\cdot,\Omega)$ is lower semi-continuous with respect to  $L^{1}_{loc}$ convergence.

Suppose that $E_{1}$ and $E_{2}$ are horizontal minimizing hulls in $\Omega_{1}$ and $\Omega_{2} $, respectively, and there is a compact set $K$ such that $ K \subset \subset \Omega \subset \subset \Omega_{1} \cap \Omega_{2}$. Then, using Lemma \ref{lem2.15} to  $E= \left(E_{1} \cap E_{2}\right)\cup K$ and $F=E_{1}$, and combining $P_{\mathbb{H}}(E_{1}) \leq P_{\mathbb{H}}(E_{1} \cup K)$, we have
\begin{align}\label{4.4}
    P_{\mathbb{H}}\left((E_{1} \cap E_{2})\cup K\right) \geq P_{\mathbb{H}}\left((E_{1} \cap E_{2})\cup (K\cap E_{1})\right).
\end{align}
Similarly, using Lemma \ref{lem2.15} to  $E^{'}= \left(E_{1} \cap E_{2}\right)\cup (K\cap E_{1})$ and $F^{'}=E_{2}$, and combining $P_{\mathbb{H}}(E_{2}) \leq P_{\mathbb{H}}(E_{2} \cup (K\cap E_{1}))$, we obtain
\begin{align}\label{4.5}
    P_{\mathbb{H}}\left((E_{1} \cap E_{2})\cup (K\cap E_{1}) \right) \geq P_{\mathbb{H}}\left(E_{1} \cap E_{2}\right) .
\end{align}
Combining (\ref{4.4}) and (\ref{4.5}), we have
\begin{align*}
    P_{\mathbb{H}}\left((E_{1} \cap E_{2})\cup K\right) \geq  P_{\mathbb{H}}\left(E_{1} \cap E_{2}\right).
\end{align*}
Thus $E_{1} \cap E_{2}$ is a horizontal minimizing hull in $\Omega$.

Now suppose that $E_{i}$ are horizontal minimizing hulls in $\Omega_{i}$, and let $A_{j} = \bigcap_{i=1}^{j} E_{i}$. It follows from the above arguments that $A_{j}$ is a horizontal minimizing hull in both $\bigcap_{i=1}^{j} \Omega_{i}$ and $\Omega$, where $\Omega \subset \bigcap_{i} \Omega_{i}$. Moreover, let $E= \bigcap_{j} A_{j}$, we have
\begin{align*}
    \int_{ K} |\chi_{A_{j}} -  \chi_{E} | d \mathcal{H}^{2} \to 0 , \qquad \text{ as } j \to \infty .
\end{align*}
Thus, $A_{j} \to E$ in $L^{1}_{loc}(\Omega)$. Since $\mathcal{E}(\cdot,\Omega)$ is lower semi-continuous, we have
\begin{align*}
    \mathcal{E}(E,\Omega) \leq \lim_{j \to \infty} \inf \mathcal{E}(A_{j},\Omega) = 0.
\end{align*}
Hence, $\mathcal{E}(E,\Omega) = 0$ and $E= \bigcap_{j} A_{j}$ is a horizontal minimizing hull in $\Omega$.
\hfill${\square}$
\begin{remark}\label{rem4.7}
    Analogous results of Proposition \ref{prop4.6} in Euclidean space and complete Riemannian manifolds have been established in \cite{RI1984} and \cite{FM2022}, respectively.
\end{remark}

Let $E_{t} = \{ u < t \}$, $E_{t}^{+} = int\{u \leq t\}$, and $E_{t}^{'}$ be the intersection of all the strictly horizontal minimizing hulls in $\Omega$ containing $E_{t}$. The horizontal minimizing hull has the following important properties, analogous to the Riemannian case in \cite{HI2001}.
\begin{proposition}\label{prop4.8}
    Assume that $E_{t}$ is a minimizer of the functional (\ref{1.6}) in $\Omega$ with initial condition $E_{0}$, and $\mathbb{H}^{1}$ has no compact component. Then \\
    (1) $\{E_{t}\}_{t>0}$ are horizontal minimizing hulls in $\mathbb{H}^{1}$.
    \\
    (2) $\{E_{t}^{+}\}_{t\geq 0}$ are strictly horizontal minimizing hulls in $\mathbb{H}^{1}$.
    \\
    (3) $\{E_{t}^{+}\}_{t\geq 0}= \{E^{'}_{t}\}_{t\geq 0}$ holds as long as $E^{+}_{t}$ is precompact.
    \\
    (4) $P_{\mathbb{H}}(E_{t}) = P_{\mathbb{H}}(E^{+}_{t})$ holds for $t > 0$ as long as $E^{+}_{t}$ is precompact. In addition, if $E_{0}$ is a horizontal minimizing hull, the equality also holds for $t=0$.
\end{proposition}
\noindent {\bf Proof}:
(1) Since $J_{u}(E_{t}) \leq J_{u}(F)$ holds for any $E_{t} \subseteq F$ and $t>0$, it follows that
\begin{align*}
    P_{\mathbb{H}}(E_{t}) + \int_{F\setminus E_{t}} |\nabla_{0} u| dx \leq  P_{\mathbb{H}}(F).
\end{align*}
Thus, $P_{\mathbb{H}}(E_{t}) \leq  P_{\mathbb{H}}(F)$ for any $t>0$.
\\
(2) Since $|\nabla_0 u| = 0$ almost everywhere on the level set $\{u(x) = t \}$ and $E_{t} =\{u < t\}$ minimizes the functional (\ref{1.6}) on $\Omega$, it follows that $E^{+}_{t}$ is also a minimizer of the functional (\ref{1.6}). Specifically, for any set $F$ such that $F\Delta E^{+}_{t}\subset K \subset \subset \Omega$ for some compact set $K$, we have
\begin{align*}
    J_{u}(E^{+}_{t}) \leq J_{u}(F).
\end{align*}
If $E^{+}_{t} \subseteq  F$, we get
\begin{align}\label{4.6}
    P_{\mathbb{H}}(E^{+}_{t}) + \int_{F\setminus E^{+}_{t}} |\nabla_{0} u| dx \leq  P_{\mathbb{H}}(F).
\end{align}
Hence, $E_{t}^{+}$ is a horizontal minimizing hull for any $t \geq 0$.

Suppose $P_{\mathbb{H}}(E^{+}_{t}) =  P_{\mathbb{H}}(F)$. Then, $F$ is a horizontal minimizing hull. It follows from (\ref{4.6}) that $|\nabla_{0} u |= 0$ almost everywhere on $F\setminus E^{+}_{t}$. This ensures that $u$ is a constant on each connected component of $F\setminus (E^{+}_{t} \cup \partial E^{+}_{t} )$. By modifying a zero measure set, we may assume that $F$ is an open set. Since $\mathbb{H}^{1}$ has no compact component and $F$ is a horizontal minimizing hull with no wasted $\mathbb{H}$-perimeter, it follows that  no such component of $F\setminus (E^{+}_{t} \cup \partial E^{+}_{t} )$ is disjoint from  $E^{+}_{t} \cup \partial E^{+}_{t}$. Hence, $u = t$ almost everywhere on $F \setminus E_{t}^{+}$, and $F \setminus E_{t}^{+} \subseteq \partial E_{t}^{+}$. Consequently, $F \cap \Omega = E_{t}^{+} \cap \Omega$ holds almost everywhere, which shows that $E_{t}^{+}$ is a strictly horizontal minimizing hull.
\\
(3) By the definition of $E_{t}^{'}$ and the precompactness of $E^{+}_{t}$, we have $E_{t}^{'} \subseteq E^{+}_{t}$ and $P_{\mathbb{H}}(E^{'}_{t}) \leq P_{\mathbb{H}}(E^{+}_{t})$. If $P_{\mathbb{H}}(E^{'}_{t}) < P_{\mathbb{H}}(E^{+}_{t})$, then
\begin{align*}
    J_{u}(E_{t}^{'}) &= P_{\mathbb{H}}(E_{t}^{'}) - \int_{E_{t}^{'}} | \nabla_{0} u| dx  \\
    &<  P_{\mathbb{H}}(E^{+}_{t}) - \int_{E_{t}^{'}} | \nabla_{0} u| dx  \leq  P_{\mathbb{H}}(E^{+}_{t}) - \int_{E_{t}} | \nabla_{0} u| dx \\
    &= P_{\mathbb{H}}(E^{+}_{t}) - \int_{E_{t}^{+}} | \nabla_{0} u | dx = J_{u}(E_{t}^{+}),
\end{align*}
which contradicts the fact that $E_{t}^{+}$ is a minimizer of the functional (\ref{1.6}). Thus $P_{\mathbb{H}}(E^{'}_{t})=P_{\mathbb{H}}(E^{+}_{t})$. From (2), we obtain $E_{t}^{'} = E^{+}_{t}$.
\\
(4) According to (1) and (3), we have $E_{t} \subseteq E_{t}^{'} = E _{t}^{+}$ and $P_{\mathbb{H}}(E_{t}) \leq P_{\mathbb{H}}(E^{+}_{t})$. Since $E_{t}^{+}$ is a minimizer of the functional (\ref{1.6}), we obtain
\begin{align*}
    P_{\mathbb{H}}(E_{t}^{+}) \leq P_{\mathbb{H}}(E_{t}) + \int_{E_{t}^{+}\setminus E_{t}} | \nabla_{0} u | dx
    \leq P_{\mathbb{H}}(E_{t}) +  \int_{\{u=t\}} | \nabla_{0} u | dx = P_{\mathbb{H}}(E_{t})
\end{align*}
Thus, $P_{\mathbb{H}}(E^{+}_{t}) = P_{\mathbb{H}}(E_{t})$ for all $t>0$. Obviously, if $E_{0}$ is a horizontal minimizing hull, the above equality also holds at $t=0$.
\hfill${\square}$

Motivated by Huisken and Ilmanen's description of the jump under IMCF in \cite{HI2001}, the horizontal jump under the HIMCF can be represented as\\
(1) Whenever $E_{t}$ remains a horizontal minimizing hull, it evolves according to the generalized horizontal inverse mean curvature flow. \\
(2) When $E_{t}$ is not a horizontal minimizing hull, $E_{t}$ horizontal jumps to $E_{t}^{'}$ and continues.

It is obvious that $E_{t}^{+} = int \{ u \leq t\} = \{\min\{u,t\}< t\}$. If $E^{+}_{t}$ is precompact, from Proposition \ref{prop4.8}, we have $E^{+}_{t} = E^{'}_{t}$, the horizontal jumping $E_{t} \to E_{t}^{+}$ can be expressed as $u \to \min\{u,t\}$.
\begin{lemma}\label{lem4.9}
    For any $t \in [0,\infty)$, $ \min \{u, t\}$ is a weak solution of (\ref{1.2}), provided that $u$ is a weak solution of (\ref{1.2}).
\end{lemma}
\noindent {\bf Proof }: Given that $J_{u}(u) \leq J_{u}(v)$, it follows from Lemma \ref{lem3.3} that $E_{t}$ minimizes $J_{u}(F)$ in $\mathbb{H}^{1} \setminus E_{0}$. Thus,
\begin{align*}
    P_{\mathbb{H}}(E_{t})- \int_{E_{t} \cap K} |\nabla_{0} u| \leq P_{\mathbb{H}}(F)- \int_{F \cap K} |\nabla_{0} u| .
\end{align*}
For $s \leq t$, we have $E_{s} \subseteq E_{t}$ and
\begin{align*}
    P_{\mathbb{H}}(E_{s}) -\int _{\{u < s\} \cap K} |\nabla_{0} u| \leq P_{\mathbb{H}}(F)- \int_{F \cap K} |\nabla_{0} u| .
\end{align*}
Replacing $u$ by $\min\{u,t\}$ only makes the left-hand side enlarge in the above inequality, the right-hand side remains unchanged. Therefore, $\min\{u,t\}$ is a weak solution of (\ref{1.2}).
\hfill${\square}$

\section{Existence of weak solutions }

\noindent  This section proves the existence of weak solutions to the initial value problem (\ref{0.1}), obtained as the limit of solutions to the Dirichlet problem (\ref{5.2a})-(\ref{5.2b}).

\subsection{ Approximation of the horizontal $p$-Laplace equation }
\
\vglue-10pt
 \indent

Let $E_{0}$ be a bounded open set with smooth boundary $\partial E_{0}$, and let $\Omega= \mathbb{H}^{1} \setminus E_{0}$. Then $\partial \Omega = \partial E_{0}$. The horizontal inverse mean curvature flow with initial condition $M_{0}=\partial E_{0}$ can be expressed as the Dirichlet equation (\ref{0.1}):
\begin{eqnarray*}
    \begin{cases}
        \operatorname{div}_{0}\left(\frac{\nabla_{0} u}{|\nabla_{0} u|}\right) = |\nabla_{0} u|,& \qquad \text{in } \Omega, \\
        u=0,& \qquad \text{on } \partial \Omega.
    \end{cases}
\end{eqnarray*}
Introducing the parameter $p \in (1,Q)$, where $Q$ is the homogeneous dimension of $\mathbb{H}^{1}$,
\begin{subequations}
    \begin{empheq}[left=\empheqlbrace]{align}
        &\operatorname{div}_{0}\left(|\nabla_{0} u_{p}|^{p-2}\nabla_{0}u_{p}\right) =|\nabla_{0}u_{p}|^{p},
       &&\text{in }\Omega, \label{5.2a}
        \\
        &u_p=0, &&\text{on }\partial\Omega.\label{5.2b}
    \end{empheq}
\end{subequations}

\begin{lemma}\label{lem5.1}
    If $u_{p}$ satisfies (\ref{5.2a}), then $u_{p}$ is a minimizer of the functional
    \begin{align*}
        \mathcal{J}_{u_{p}}^{K}(v)=\int_{K}\left(\frac{1}{p}|\nabla_{0} v|^{p} + v |\nabla_{0} u_{p}|^{p}\right) dx ,
    \end{align*}
    where $v \in HW^{1,p}_{loc}(\Omega)$, $\{v \neq u_{p}\} \subset \subset \Omega$, and $\{v \neq u_{p}\}$ contained in a compact set $K \subset \Omega$.
\end{lemma}
\noindent {\bf Proof }: Since $u_{p}$ satisfies (\ref{5.2a}), the following identity holds for any test function $\varphi \in HW^{1,p}_{\text{loc}}(\Omega)$:
\begin{align*}
    \int_{K} |\nabla_{0} u_{p}|^{p-2}\left\langle \nabla_{0}u_{p}, \nabla_{0} \varphi \right\rangle dx = - \int_{K} |\nabla_{0} u_{p}|^{p} \varphi dx.
\end{align*}
Let $\varphi = (v-u_{p})$, we have
\begin{align}\label{5.3}
     \int_{K} |\nabla_{0} u_{p}|^{p}  (u_{p}-v)dx + \int_{K} |\nabla_{0} u_{p}|^{p} dx  \leq \int_{K} |\nabla_{0} u_{p}|^{p-1} |\nabla_{0} v| dx.
\end{align}
Applying Young's inequality with exponent $\frac{p}{p-1}$ for the right-hand side of (\ref{5.3}), we get
\begin{align*}
    \int_{K} |\nabla_{0} u_{p}|^{p}  (u_{p}-v)dx \leq \int_{K} \frac{1}{p} \left(|\nabla_{0} v|^{p} - |\nabla_{0} u_{p}|^{p}\right) dx .
\end{align*}
Thus
\begin{align*}
    \int_{K}\left(\frac{1}{p}|\nabla_{0} u_{p}|^{p} + u_{p} |\nabla_{0} u_{p}|^{p}\right) dx \leq \int_{K}\left(\frac{1}{p}|\nabla_{0} v|^{p} + v |\nabla_{0} u_{p}|^{p}\right) dx.
\end{align*}
This completes the proof of Lemma \ref{lem5.1}.
\hfill${\square}$

If $w_{p} = \exp \left(\frac{u_{p}}{1-p}\right)$ and $u_p$ satisfies (\ref{5.2a})–(\ref{5.2b}), then $w_{p}$ satisfies the horizontal $p$-Laplace equation with the boundary condition
\begin{subequations}
    \begin{empheq}[left=\empheqlbrace]{align}
        &\operatorname{div}_{0}\left(|\nabla_{0} w_{p}|^{p-2} \nabla_{0} w_{p}\right) = 0,
        \qquad &&\text{in }\Omega, \label{5.4a}
        \\
        &w_{p}=1,\qquad &&\text{on }\partial\Omega.\label{5.4b}
    \end{empheq}
\end{subequations}
which is the same as the Dirichlet equation (\ref{0.3}). Consider the elliptic regularization of the horizontal $p$-Laplace equation:
\begin{subequations}
    \begin{empheq}[left=\empheqlbrace]{align}
        &\operatorname{div}_{0}\left(\left(\delta^{2} + |\nabla_{0} w_{p}|^{2}\right)^{\frac{p-2}{2}} \nabla_{0} w_{p}\right) = 0,
        \qquad \text{in }\Omega, \label{5.5a}
        \\
        &w_{p}-\psi \in HW^{1,p}_{0}(\Omega).\label{5.5b}
    \end{empheq}
\end{subequations}
where $\psi \in HW^{1,p}(\Omega)$. Let $\mathcal{A}_{1}= \left\{v \in HW^{1,p}(\Omega) , v =1 \text{ a.e. on } \partial \Omega \right\}$. Obviously, the case $\delta=0$ in (\ref{5.5a})-(\ref{5.5b}) with $\psi \in \mathcal{A}_{1}$ corresponds to the Dirichlet problem (\ref{5.4a})-(\ref{5.4b}). We say that a function $w_{p} \in HW^{1,p}(\Omega)$ is a weak solution of (\ref{5.5a}) if the following holds
\begin{align}\label{5.6}
    \int_{\Omega} \left(\delta^{2} + |\nabla_{0} w_{p}|^{2}\right)^{\frac{p-2}{2}}\left\langle \nabla_{0}w_{p}, \nabla_{0}\varphi \right\rangle  dx =0, \qquad \forall \varphi \in HW^{1,p}_{0}(\Omega).
\end{align}
From Theorem 3.2 in \cite{R2015}, for $\delta \geq 0$, $w_{p}$ is a weak solution to (\ref{5.5a}) if and only if $w_{p}$ is a minimum of the functional
\begin{align}\label{5.7}
    \mathcal{L}_{w_{p},\delta}(v) = \frac{1}{p} \int_{\Omega}  \left(\delta^{2} + | \nabla_{0} v|^{2}\right)^{\frac{p}{2}} dx, \qquad v \in HW^{1,p}(\Omega).
\end{align}
Moreover, for $1 < p < \infty$ and $\delta \geq 0$, the Dirichlet problem (\ref{5.5a})-(\ref{5.5b}) has a unique weak solution $w_{p} \in HW^{1,p}(\Omega)$, as detailed in \cite{D1995, HKM2006}. Hence, by the transformation $u_{p} = (1-p) \ln w_{p}$, $u_{p}$ is the unique weak solution of (\ref{5.2a})-(\ref{5.2b}).

Next, we present some well-known results that can be employed to establish a priori estimates for nonnegative solutions $w_{p}$ of (\ref{5.4a}). First, we recall the Harnack inequality.
\begin{theorem}\cite{CDG1993}
    Let $1 < p \leq Q$ and $x_{0} \in \Omega$. Assume that $w_{p}\in HW_{loc}^{1,p}(\Omega)$ is a nonnegative weak solution to (\ref{5.4a}). Then, there exist constants $\mathcal{C}^{h}(p,c,Q)>0$ and $R_{0} >0$ such that  for any $B_{R}=B(x_{0},R)$, with $B_{4R} \subset \Omega$, and $R \leq R_{0}$,
    \begin{align}\label{5.8}
        \esssup_{B_{R}} w_{p} \leq \mathcal{C}^{h}(p,c,Q)\essinf_{B_{R}} w_{p}
    \end{align}
    where $\mathcal{C}^{h}(p,c,Q)$ depending on $p$, the doubling constant $c$ in (\ref{2.d}) and $Q$.
\end{theorem}

To state the subsequent regularity results, it is necessary to introduce the following class of domains. Let $D \subset \mathbb{H}^{1}$ be a bounded open set satisfying the condition: there exists a constant $\iota > 0$ such that,  for each $y \in \partial D$, one can find a unit vector $\vec{b} \in \mathbb{R}^{3}$ satisfying
\begin{align}\label{5.9}
    \langle x-y, \vec{b} \rangle_{E} \geq \iota |x-y|_{E}^2, \qquad \forall x \in D.
\end{align}
Let $w_{p} \in HW^{1,p}(D)$ be a weak solution of (\ref{5.5a}) with $\delta >0$. Then, \cite{C1997,R2015,Z2017} shows that
\begin{align}\label{5.10}
    \left\|\nabla_{0} w_{p}\right\|_{L^{\infty}(D)} \leq C_{\iota},\quad \nabla_{0} w_{p} \in HW^{1,2}_{loc}(D),\quad Tw_{p} \in HW^{1,2}_{loc}(D) \cap L_{loc}^{\infty}(D)
\end{align}
where $C_{\iota}$ is a positive constant that depend on $\iota$, $\left\|\nabla \psi\right\|_{L^{\infty}(D)} $ and $\left\|\nabla^{2} \psi\right\|_{L^{\infty}(D)} $. Moreover, they also established Lipschitz estimates for $w_{p}$ in $D$, which can be extended to the case $\delta=0$ and general domains by a standard approximation argument. More precisely, the following Lipschitz estimate holds.
\begin{theorem}\cite{R2015,Z2017}
    Let $1 <p < \infty$ and $w_{p}\in HW_{loc}^{1,p}(\Omega)$ be a nonnegative weak solution to (\ref{5.4a}). Then,
    \begin{align}\label{5.11}
        \sup_{B_{R}} |\nabla_{0} w_{p}| \leq \mathcal{C}^{l}(p,c,Q)\left( \fint_{B_{2R}} |\nabla_{0} w_{p}|^{p}\right)^{\frac{1}{p}}
    \end{align}
    for every ball $B_{R}$ such that the concentric ball $B_{2R} \subset \Omega$, where $\mathcal{C}^{l}(p,c,Q) $ is a positive constant depending on $p$, $c$ and $Q$.
\end{theorem}

\subsection{ Refined Harnack inequality for solutions to (\ref{5.4a})}
\
\vglue-10pt
 \indent

Let $\mathsf{B}_{r}\left(y\right)$ denote the Kor$\acute{\rm a}$nyi ball in $\mathbb{H}^{1}$ centered at $y$ with radius $r$. Since $\Omega^{c}$ is a bounded open set, there exists $r>0$ such that $\mathsf{B} _{r}\left(y\right) \subset \Omega^{c}$. From \cite{CDG1996}, there exists an explicit solution $\Gamma_{p}(x,y)$ to (\ref{5.4a}) with $1<p<Q$, given by
\begin{equation*}
    \begin{cases}
        \Gamma_{p}(x,y)=\left(\frac{r}{\rho(x,y)}\right)^{(Q-p)/(p-1)} \qquad &\text{in } \Omega, \\
        \Gamma_{p}(x,y)\leq 1 \qquad &\text{on } \partial\Omega,
    \end{cases}
\end{equation*}
where $\rho(x,y)$ is the Kor$\acute{\text{a}}$nyi distance between $x$ and $y$. By the comparison principle for the solutions of the Dirichlet problem (\ref{5.4a})-(\ref{5.4b}) (see, e.g., Theorem 2.5 in \cite{D1995} or Proposition 2.5 in \cite{FP2021}), we have
\begin{align}\label{5.12}
    0 < \left(\frac{r}{\rho(x,y)}\right)^{\frac{Q-p}{p-1}} \leq w_{p}(x) < 1, \qquad \forall x \in \Omega.
\end{align}
Combining the above inequality with the Harnack inequality (\ref{5.8}), one obtains a uniform $C^{0}$-estimate for $w_{p}$. However, since the constant $\mathcal{C}^{h}(p,c,Q)$ in (\ref{5.8}) is not given explicitly in terms of $p$, the $C^{0}$-estimate derived from (\ref{5.8}) cannot be used to analyze the limit as $p\to 1$.

In the following, we will derive an explicit expression for $\mathcal{C}^{h}(p,c,Q)$ in terms of $p$. In particular, unless otherwise specified, $C$ denotes a generic absolute constant, and $C(c,Q)$ denotes a generic positive constant depending only on $c$ and $Q$, but independent of $p$. Both constants may vary from line to line.

The Sobolev inequality (\ref{2.s}) plays a fundamental role in the proof of the Harnack inequality. Therefore, it is necessary to obtain an explicit expression for the constant $\mathcal{C}_{s}(p,c,Q)$ in (\ref{2.s}) in terms of the parameter $p$. Although the argument is standard, the details are included in the appendix (see Section 9.1) for the convenience of the reader.

\begin{theorem}\label{thmS}
    Let $U \subset \mathbb{H}^{1}$ be a bounded open set. Let $1<p<Q$ and $1 \leq \kappa \leq \frac{Q}{Q-p}$. Then, there exists $R_0>0$ such that for any $x_{0} \in U, B_{R}=B(x_{0}, R)$, with $R \leq R_0$, we have
    \begin{align}\label{5.s}
        \left(\fint_{B_R}|u|^{\kappa p} d x\right)^{\frac{1}{\kappa p}} \leq \mathcal{C}_{s} R\left(\fint_{B_R}\left|\nabla_{0} u\right|^p d x\right)^{\frac{1}{p}}
    \end{align}
    for any $u \in HW^{1, p}\left(B_R\right)$, where the constant $\mathcal{C}_{s}$ is given by
    \begin{equation*}
        \mathcal{C}_{s} =
        \begin{cases}
            \frac{C(c,Q)}{\left(2^{\frac{Q-p}{\kappa(p-1)+1}\left(\frac{Q}{Q-p} -\kappa\right)}-1\right)^{\frac{\kappa(p-1)+1}{\kappa p}}}, \qquad  1 \leq \kappa < \frac{Q}{Q-p}\\
            \frac{C(c, Q)}{Q-p}\left(1 + \frac{1}{2^{\frac{Q-p}{p-1}} -1 }\right)^{Q-1/Q} \left(\frac{p2^{p-1}}{p-1}\right)^{\frac{Q-p}{pQ}}, \qquad \kappa =  \frac{Q}{Q-p}
        \end{cases}
    \end{equation*}
\end{theorem}

\begin{remark}\label{remS}
    In this paper, the inequality (\ref{5.s}) will be used in the following three cases:
    \begin{align*}
        &p=2 \text{ with } \kappa = \frac{Q}{Q-2},  \quad \mathcal{C}_{s} = \frac{C(c,Q)}{Q-2} \left(1 + \frac{1}{2^{Q-2} -1}\right)^{Q-1/Q}4^{\frac{Q-2}{2Q}};\\
        &p \in(1,Q) \text{ with } \kappa =1, \quad \mathcal{C}_{s} =C(c,Q); \\
        &p \in(1,2) \text{ with } \kappa =\frac{2Q}{2Q-1}, \quad \mathcal{C}_{s} =\frac{C(c,Q)}{\left(2^{\frac{Q(2p-1)}{2pQ-1}} -1\right)^{\frac{2pQ-1}{2pQ}}} <36C(c,Q)
    \end{align*}
    In the above cases, $\mathcal{C}_{s}$ can be bounded by a constant depending only on $c$ and $Q$, but independent of $p$.
\end{remark}

Let $u_{B} =\fint_{B_{R}} u(x) dx $ and $\bar{u}(x) =  u(x)-u_{B}$. Then, applying (\ref{5.s}) with $\kappa =1$ for $\bar{u}(x)$, together with H$\ddot{\rm o} $lder inequality, we obtain the following Poincar$\acute{\rm e}$ inequality.
\begin{corollary}\label{cor5.4}
    There exists a constant $C(c,Q) >0$, such that for all $ u \in C^{1}(B_{R})$, we have
    \begin{align}\label{5.14}
        \fint_{B_{R}}|u(x)-u_{B}| dx \leq C(c,Q)R\left(\fint_{B_{R}}|\nabla_{0} u|^{p}dx\right)^{\frac{1}{p}}
    \end{align}
    where $u_{B} =\fint_{B_{R}} u(x) dx $.
\end{corollary}

Next, we establish a refined Harnack inequality for weak solutions to (\ref{5.4a}) by deriving an explicit expression for the constant $\mathcal{C}^{h}(p,c,Q)$ in (\ref{5.8}) in terms of $p$. To analyze the limiting case $p \to 1$, it suffices to consider the case $p \in (1,2)$.

\begin{theorem}\label{thm5.5}(Refined Harnack inequality)
    Let $\tau_{1} >0$ be a constant satisfying $\tau_{1} < \min\left\{1,  \left(7QC^{*}(c,Q)\right)^{-1/2}\right\}$. Assume that $1 < p \leq 1+\tau_{1}$ and $w_{p}\in HW_{loc}^{1,p}(\Omega)$ is a nonnegative weak solution to (\ref{5.4a}). Then, there exists a constant $R_{0}>0 $ such that for any $x_{0} \in \Omega$, $B_{R} = B(x_{0}, R)$ with $B_{4R} \subset \Omega$ and $R \leq R_{0}$,
    \begin{align*}
        \esssup_{B_{R}} w_{p} \leq \left[\mathcal{C}^{h}(c,Q)\left(\frac{p-1}{2-p}\right)^{2}\right]^{1/\alpha} \essinf_{B_{R}} w_{p}
    \end{align*}
    where $\alpha = C^{*}(c,Q)\frac{(p-1)^{2}}{p} $, $C^{*}(c,Q)$ and $\mathcal{C}^{h}(c,Q)$ are positive constants depending only on $c$ and $Q$.
\end{theorem}
\noindent {\bf Proof }: Firstly, we derive a Caccioppoli type inequality for arbitrary powers of $w_{p}$. Let $\beta \in \mathbb{R}\backslash \{-1/p\}$, $\eta \in C_{0}^{\infty}(B_{4R})$ and $\varphi = \eta^{p}w_{p}^{1+\beta p}$. Using $\varphi$ as a test function in weak formula of (\ref{5.4a}) and Young's inequality with exponent $\frac{p}{p-1}$, we have
\begin{align}\label{5.15}
    \int_{\Omega}w_{p}^{\beta p}|\nabla_{0} w_{p}|^{p}\eta^{p}  \leq \left(\frac{p}{|1+\beta p|}\right)^{p} \int_{\Omega}|\nabla_{0}\eta|^{p} w_{p}^{p+\beta p}
\end{align}
For $\beta =-1$, (\ref{5.15}) yields
\begin{align}\label{5.16}
    \int_{\Omega} \eta^{p}  |\nabla_{0}(\ln w_{p})|^{p}
    \leq \left(\frac{p}{p-1}\right)^{p} \int_{\Omega}|\nabla_{0}\eta|^{p}.
\end{align}
For $\beta \in \mathbb{R}\backslash \{-\frac{1}{p}, -1\}$, (\ref{5.15}) yields
\begin{align*}
    \int_{\Omega} |\nabla_{0}(\eta w_{p}^{1+\beta})|^{p} \leq 2^{p-1}  \left( 1+ \left(\frac{p|1+\beta|}{|1+\beta p|}\right)^{p}\right)\int_{\Omega}|\nabla_{0}\eta|^{p}w_{p}^{(1+\beta)p}.
\end{align*}
When $\beta \in [0, +\infty)$, $\left(\frac{p|1+\beta|}{|1+\beta p|}\right)^{p} < p^{p} <4$. When $\beta \in (-\infty, -1)$, $\left(\frac{p|1+\beta|}{|1+\beta p|}\right)^{p}<1$. Hence, there exists an absolute constant $C$ such that
\begin{align}\label{5.17}
    \int_{\Omega} |\nabla_{0}(\eta w_{p}^{1+\beta})|^{p} \leq C \int_{\Omega}|\nabla_{0}\eta|^{p}w_{p}^{(1+\beta)p},\qquad \beta \in (-\infty,-1) \cup[0, +\infty).
\end{align}
Next, by specializing inequality (\ref{5.16}), we show that there exist constants $C>0$ and $\alpha >0$ such that $\left( \fint_{B_{2R}} w_{p}^{\alpha} dx \right)^{\frac{1}{\alpha}} \leq C \left(\fint_{B_{2R}} w_{p}^{-\alpha} dx\right)^{-\frac{1}{\alpha}}$ holds. From Lemma 3.6 of \cite{CGL1993}, we know that there exists a cut-off function $\eta \in C^{\infty}_{0}(B_{4R})$ such that $\eta \equiv 1$ in $B_{2R}$ and  $|\nabla_{0} \eta| \leq \frac{C}{R}$, where $C$ independent of $R$. Let $v = \ln w_{p}$. Then, by applying (\ref{5.14}) to $\eta v$ and using Lemma \ref{lemD} together with (\ref{5.16}), we get
\begin{align*}
    \fint_{B_{2R}}|v-v_{B}| &\leq C(c,Q)R \left(\fint_{B_{2R}} \eta^{p}|\nabla_{0} v|^{p}\right)^{\frac{1}{p}} \leq C(c,Q)R \left(\frac{|B_{4R}|}{|B_{2R}|}\right)^{\frac{1}{p}}\left(\fint_{B_{4R}} \eta^{p}|\nabla_{0} v|^{p}\right)^{\frac{1}{p}} \\
    &\leq C(c,Q)\frac{pR}{p-1} \left(\fint_{B_{4R}} |\nabla_{0} \eta|^{p}\right)^{\frac{1}{p}} \leq C(c,Q)\frac{p}{p-1}
\end{align*}
where the last equality follows from the properties of $\eta$. Thus, $ v \in BMO(B_{2R},d)$ and
\begin{align}\label{5.18}
    \| v\|_{BMO} = \sup \fint_{B_{2R}}|v-v_{B}|  = C(c,Q)\frac{p}{p-1}.
\end{align}
From Theorem 0.4 and Theorem 2.2 in \cite{Bu1999}, the John-Nirenberg inequality is given by
\begin{align}\label{5.19}
    |\{x\in B_{2R}, |v-v_{B}| > t\}| \leq C_{1}e^{-\frac{C_{2}t}{\| v\|_{BMO}}} |B_{2R}|,
\end{align}
where $C_{1}$ and $C_{2}$ depending  only on $c$. Let $ t = \frac{\ln s}{\alpha}$, we have $\chi_{\left\{e^{\alpha |v-v_{B}|} > s\right\}} = \chi_{\left\{|v-v_{B}| > t\right\}}$ and $ds = \alpha e^{\alpha t}dt$. Hence, by (\ref{5.18}) and (\ref{5.19}), we obtain
\begin{align}\label{5.20}
    \int_{B_{2R}} e^{\alpha |v-v_{B}|} dx &= \int_{B_{2R}}  \int_{0}^{\infty}\chi_{\{e^{\alpha |v-v_{B}|} > s\}} ds dx  = \int_{0}^{\infty}\alpha e^{\alpha t} \left|\{x\in B_{2R},|v-v_{B}| > t\}\right| dt \notag \\
    &\leq C_{1}|B_{2R}| \alpha \int_{0}^{\infty}e^{\left(\alpha -\frac{C_{2}}{\| v\|_{BMO}}\right)t}dt \notag \\
    &= C_{1}|B_{2R}| \alpha\left(\alpha -\frac{C_{2}}{C(c,Q)}\frac{p-1}{p}\right)^{-1} e^{\left(\alpha -\frac{C_{2}}{C(c,Q)}\frac{p-1}{p}\right)t}\Big|^{\infty}_{0}.
\end{align}
Let $C^{*}(c,Q) := \frac{C_{2}}{C(c,Q)}$ and choose $\alpha = C^{*}(c,Q)\frac{(p-1)^{2}}{p}$, we have
\begin{align*}
    \int_{B_{2R}} e^{\alpha |v-v_{B}|} dx \leq C_{1}|B_{2R}|\frac{p-1}{2-p} .
\end{align*}
Substituting $v = \ln w_{p}$ into the above inequality, we obtain
\begin{align*}
    e^{-\alpha v_{B}} \int_{B_{2R}} w_{p}^{\alpha} dx \leq C_{1}|B_{2R}|\frac{p-1}{2-p}, \quad e^{\alpha v_{B}} \int_{B_{2R}} w_{p}^{-\alpha} dx \leq C_{1}|B_{2R}|\frac{p-1}{2-p}.
\end{align*}
Thus
\begin{align}\label{5.21}
    \left( \fint_{B_{2R}} w_{p}^{\alpha} dx \right)^{\frac{1}{\alpha}} &\leq \left(C_{1}\frac{p-1}{2-p}\right)^{\frac{2}{\alpha}} \left(\fint_{B_{2R}} w_{p}^{-\alpha} dx\right)^{-\frac{1}{\alpha}} ,
\end{align}
where $\alpha = C^{*}(c,Q)\frac{(p-1)^{2}}{p}$.

In the following, we establish the $L^{\infty}$-estimate for $w_{p}$ and $w_{p}^{-1}$ on $B_{R}$ respectively. We begin with the $L^{\infty}$-estimate for $w_{p}$. Let $R \leq R_{1} < R_{2} \leq 2R$, and suppose the cut-off function $\eta \in C_{0}^{\infty} (B_{2R})$ satisfies  $\eta  \equiv 1$ in $B_{R_{1}}$, $spt(\eta) \subset B_{R_{2}}$ , $ 0< \eta \leq 1 $ in $B_{R_{2}}$ and $|\nabla_{0} \eta| \leq\frac{C}{R_{2}-R_{1}}$, where $C$ independent of $R$. Applying (\ref{5.s}) with $\kappa = \frac{2Q}{2Q -1}$ to $\eta w_{p}^{1+\beta}$ on $B_{R_{2}}$, we have
\begin{align}\label{5.22}
    \left( \fint_{B_{R_{2}}}|\eta w_{p}^{1+\beta}|^{\kappa p} d x\right)^{\frac{1}{\kappa p}} & \leq \mathcal{C}_{s} R_{2} \left( \fint_{B_{R_{2}}}\left|\nabla_{0} (\eta w_{p}^{1+\beta})\right|^p d x\right)^{\frac{1}{p}}  \notag\\
    &\leq C(c,Q) R_{2} \left( \fint_{B_{R_{2}}}\left|\nabla_{0} (\eta w_{p}^{1+\beta})\right|^p d x\right)^{\frac{1}{p}}
\end{align}
where the last inequality follows from Remark \ref{remS}. By Lemma \ref{lemD}, (\ref{5.17}), (\ref{5.22}) and the properties of  $\eta$, we obtain that, for any $\beta \in (-\infty, -1) \cup [0, +\infty)$
\begin{align}\label{5.23}
    &\left( \fint_{B_{R_{1}}}  w_{p}^{(1+\beta)\kappa p} d x\right)^{\frac{1}{\kappa p}} \notag \\
    &\qquad \leq C(c,Q)R_{2}\left( \fint_{B_{R_{2}}}\left|\nabla_{0} (\eta w_{p}^{1+\beta})\right|^p d x\right)^{\frac{1}{p}}
    \leq C(c,Q)R_{2}\left(\fint_{B_{R_{2}}}|\nabla_{0}\eta|^{p}w_{p}^{(1+\beta)p}\right)^{\frac{1}{p}} \notag \\
    &\qquad \leq C(c,Q)\frac{R_{2}}{R_{2} -R_{1}} \left( \fint_{B_{R_{2}}}w_{p}^{(1+\beta)p}\right)^{\frac{1}{p}}
\end{align}
To perform the Moser iteration, define
\begin{align*}
    &\beta_{0} = 0,\quad \beta_{i+1} = \frac{1}{2Q-1}+\beta_{i}\frac{2Q}{2Q-1}; \\
    &s_{0}= p,\quad  s_{i} = (1+\beta_{i})p,  \quad s_{i+1} = \kappa s_{i},\quad s_{i}= \kappa ^{i}s_{0};\\
    &R_{0} = 2R,\quad  R_{i} = R\left(1 + 2^{-i}\right).
\end{align*}
For any $i \geq 0$, $\beta_{i} \geq 0$. Using (\ref{5.23}) with $(R_{1}, R_{2},\beta) = (R_{i+1}, R_{i}, \beta_{i})$, we have
\begin{align*}
    \left(\fint_{B_{R_{i+1}}}w_{p}^{s_{i+1}} d x\right)^{\frac{1}{s_{i+1}}} &\leq C(c,Q)^{\frac{p}{s_{i}}}\left(\frac{R_{i}}{R_{i} - R_{i+1}}\right)^{\frac{p}{s_{i}}} \left(\fint_{B_{R_{i}}}w_{p}^{s_{i}}\right)^{\frac{1}{s_{i}}} \\
    &\leq C(c,Q)^{\frac{p}{s_{i}}}2^{\frac{p(i+2)}{s_{i}}} \left(\fint_{B_{R_{i}}}w_{p}^{s_{i}}\right)^{\frac{1}{s_{i}}}.
\end{align*}
Iterating the above inequality with exponent $\{ s_{i}\}_{i \geq 0}$ yields
\begin{align}\label{5.24}
    \|w_{p}\|_{L^{\infty}(B_{R })} \leq (C(c,Q))^{\sum_{i = 0}^{\infty} \frac{p}{s_{i}} }  2^{\sum_{i = 0}^{\infty} \frac{p(i+2)}{s_{i}} } \left(\fint_{B_{2R}}w_{p}^{p}\right)^{\frac{1}{p}}.
\end{align}
A direct computation shows that
\begin{align}\label{5.25}
   \sum_{i = 0}^{\infty} \frac{p}{s_{i}} = 2Q,\quad \sum_{i = 0}^{\infty} \frac{p(i+1)}{s_{i}} = 4Q^{2},\quad
    \sum_{i = 0}^{\infty} \frac{p(i+2)}{s_{i}} = 2Q(2Q+1).
\end{align}
Plugging (\ref{5.25}) into (\ref{5.24}), we get
\begin{align}\label{5.26}
    \|w_{p}\|_{L^{\infty}(B_{R})} \leq (C(c,Q))^{2Q} 2^{2Q(2Q+1)}\left(\fint_{B_{2R}}w_{p}^{p}\right)^{\frac{1}{p}} \leq  C(c,Q)\left( \fint_{B_{2R}}w_{p}^{p}\right)^{\frac{1}{p}}.
\end{align}
If $ \alpha \geq p $, using H$\ddot{\rm o} $lder inequality with exponent $\alpha /p$, we obtain
\begin{align}\label{5.27}
    \|w_{p}\|_{L^{\infty}(B_{R})} \leq  C(c,Q) \left( \fint_{B_{2R}}w_{p}^{p}\right)^{\frac{1}{p}} \leq  C(c,Q) \left( \fint_{B_{2R}}w_{p}^{\alpha}\right)^{\frac{1}{\alpha}}.
\end{align}
If $\alpha < p$, by adapting an argument of Dahlberg and Kenig, we can decrease the exponent on the right-hand side of (\ref{5.26}) from $p$ to $\alpha $. For any $\frac{1}{2} \leq a < b \leq 1$, let $R_{i} = 2R(a+(b-a)2^{-i})$, we have
\begin{align*}
    R_{0} = 2bR,\quad R_{\infty} = 2aR ,\quad \frac{R_{i}}{R_{i}-R_{i+1}} =  \frac{a+(b-a)2^{-i}}{(b-a)2^{-i-1}}  \leq (b-a)^{-1} 2^{i+2},
\end{align*}
Repeating the above iteration procedur, we obtain
\begin{align*}
    \|w_{p}\|_{L^{\infty}(B_{2aR})} &\leq (C(c,Q)(b-a)^{-1})^{\sum_{i = 0}^{\infty} \frac{p}{s_{i}} } 2^{\sum_{i = 0}^{\infty} \frac{p(i+2)}{s_{i}} } \left( \fint_{B_{2bR}}w_{p}^{p}\right)^{\frac{1}{p}}  \\
    &\leq C(c,Q) (b-a)^{-2Q}\left(\fint_{B_{2bR}}w_{p}^{p}\right)^{\frac{1}{p}}.
\end{align*}
From Lemma \ref{lemD}, $|B_{2bR}|\geq cb^{Q}|B_{2R}| \geq c2^{-Q} |B_{2R}|$. Thus,
\begin{align}\label{5.28}
    \|w_{p}\|_{L^{\infty}(B_{2aR})} &\leq C(c,Q)(b-a)^{-2Q}\left(\frac{|B_{2R}|}{|B_{2bR}|}\right)^{\frac{1}{p}}\left(\frac{1}{|B_{2R}|}\int_{B_{2bR}}w_{p}^{p}\right)^{\frac{1}{p}} \notag \\
    &\leq C(c,Q)(b-a)^{-2Q}\left(\frac{1}{|B_{2R}|}\int_{B_{2bR}}w_{p}^{p}\right)^{\frac{1}{p}}.
\end{align}
Define
\begin{align}\label{5.29}
    J(s) = |B_{2R}|^{\frac{p-\alpha}{\alpha}} \int_{B_{2sR}} w_{p}^{p} dx \left(\int_{B_{2R}} w_{p}^{\alpha} dx\right)^{-\frac{p}{\alpha}}.
\end{align}
Obviously, $J(s)$ is monotonically increasing with respect to $s \in [1/2,1]$ and
\begin{align*}
    J\left(\frac{2}{3}\right)^{\frac{\alpha}{p}}= |B_{2R}|^{\frac{p-\alpha}{p}} \left(\int_{B_{4R/3}} w_{p}^{p} dx \right)^{\frac{\alpha}{p}}\left(\int_{B_{2R}} w_{p}^{\alpha} dx\right)^{-1}.
\end{align*}
Choose $a =1/2$ and $b=2/3$ in (\ref{5.28}), we have
\begin{align}\label{5.30}
    \|w_{p}\|^{\alpha}_{L^{\infty}(B_{R})} &\leq  C(c,Q)^{\alpha}6^{2Q\alpha}\left(\frac{1}{|B_{2R}|} \int_{B_{4R/3}}w_{p}^{p}\right)^{\frac{\alpha}{p}}  \notag \\
    &\leq  C(c,Q)^{\alpha}J\left(\frac{2}{3}\right)^{\frac{\alpha}{p}} \left(\fint_{B_{2R}} w_{p}^{\alpha} dx\right).
\end{align}
Assume $J\left(\frac{2}{3}\right) > 1$, otherwise (\ref{5.27}) follows. It remains to prove that $J\left(\frac{2}{3}\right)$ is bounded. Combining (\ref{5.28}) and (\ref{5.29}), we have
\begin{align}\label{5.31}
    J(a) &\leq \|w_{p}\|^{p-\alpha}_{L^{\infty}(B_{2aR})} |B_{2R}|^{\frac{p-\alpha}{\alpha}} \int_{B_{2aR}} w_{p}^{\alpha} dx \left(\int_{B_{2R}} w_{p}^{\alpha} dx\right)^{-\frac{p}{\alpha}} \notag \\
    & \leq C(c,Q)^{p-\alpha}(b-a)^{-2Q(p-\alpha)} \left(|B_{2R}|^{\frac{p-\alpha}{\alpha}} \int_{B_{2bR}} w_{p}^{p} dx \left(\int_{B_{2R}} w_{p}^{\alpha} dx\right)^{-\frac{p}{\alpha}}\right)^{\frac{p-\alpha}{p}}  \notag \\
    &\leq C(c,Q)^{p-\alpha}(b-a)^{-2Q(p-\alpha)}J(b)^{\frac{p-\alpha}{p}}.
\end{align}
Let $\beta \in (0,1/3]$, $1<\theta \leq \frac{3}{2}$, $a = b^{\theta}$ and $a \in \left[1-\beta,1\right]$. Then, $J(b) \geq J(a) \geq J(\frac{2}{3}) >1$ and (\ref{5.31}) implies that
\begin{align}\label{5.32}
    \ln J(b^{\theta})  \leq (p-\alpha) \ln C(c,Q) - 2Q(p-\alpha) \ln (b-b^{\theta}) +\frac{p-\alpha}{p} \ln J(b).
\end{align}
Integrating both sides of (\ref{5.32}) with respect to $db$ over $\left[\left(1-\beta \right)^{1/\theta},1\right]$ , we obtain
\begin{align}\label{5.33}
    \int_{\left(1-\beta \right)^{\frac{1}{\theta}}}^{1} \ln J(b^{\theta}) db
   \leq&  (p-\alpha) \ln C(c,Q) \int_{\left(1-\beta \right)^{\frac{1}{\theta}}}^{1}db -2Q(p-\alpha)\int_{\left(1-\beta \right)^{\frac{1}{\theta}}}^{1} \ln (b-b^{\theta})db \\
   &+  \frac{p-\alpha}{p} \int_{\left(1-\beta\right)^{\frac{1}{\theta}}}^{1} \ln J(b)db \notag.
\end{align}
Since $db^{\theta} = \theta b^{\theta -1} db $ and $b \in [(1-\beta)^{\frac{1}{\theta}},1]\subset [2/3,1]$, we have
\begin{align}\label{5.34}
    db \geq \frac{1}{\theta} db^{\theta}, \qquad db \leq  \frac{1}{\theta} \left(\frac{3}{2}\right)^{\theta -1}db^{\theta}.
\end{align}
Let $f(b) = (1-b^{\theta-1}) - (\theta -1)^{2} (1-b^{\theta})$. It is straightforward to verify that $f(1)=0$ and $f(b)$ is monotonically decreasing with respect to $b \in [2/3, 1]$. Thus,
\begin{align}\label{5.35}
    (1-b^{\theta-1}) \geq (\theta -1)^{2} (1-b^{\theta}), \qquad \forall  b \in [2/3, 1].
\end{align}
Substituting (\ref{5.34}) and (\ref{5.35}) into (\ref{5.33}), we have
\begin{align}\label{5.36}
    \frac{1}{\theta}\int_{\left(1-\beta\right)}^{1} \ln J(b) db
    \leq&  (p-\alpha) \ln C(c,Q) \beta - 2Q(p-\alpha)\int_{\left(1-\beta \right)^{\frac{1}{\theta}}}^{1} (\ln b) db \notag \\
    &-  2Q(p-\alpha)\int_{\left(1-\beta\right)^{\frac{1}{\theta}}}^{1} \ln (1-b^{\theta-1}) db + \frac{p-\alpha}{p} \int_{\left(1-\beta\right)^{\frac{1}{\theta}}}^{1} \ln J(b)db \notag \\
    \leq & (p-\alpha) \ln C(c,Q)\beta + 2Q(p-\alpha)\ln (3/2)\beta - 4Q(p-\alpha)\beta\ln(\theta -1) \\
    &- 2Q(p-\alpha) \int_{\left(1-\beta\right)^{\frac{1}{\theta}}}^{1} \ln (1-b^{\theta}) db +  \frac{p-\alpha}{p} \int_{\left(1-\beta\right)^{\frac{1}{\theta}}}^{1} \ln J(b)db  \notag
\end{align}
where we used $\left(1-\beta\right)^{1/\theta} \geq \left(1-\beta\right)$ and $b^{-1} \leq (1-\beta)^{-1} \leq 3/2$. Also, by (\ref{5.34}), we have
\begin{align}\label{5.361}
    - \int_{\left(1-\beta\right)^{\frac{1}{\theta}}}^{1} \ln (1-b^{\theta}) db &\leq  - \frac{1}{\theta}\left(\frac{3}{2}\right)^{\theta -1}\int_{\left(1-\beta\right)^{\frac{1}{\theta}}}^{1} \ln (1-b^{\theta}) db^{\theta} = - \frac{1}{\theta}\left(\frac{3}{2}\right)^{\theta -1}\int_{\left(1-\beta\right)}^{1} \ln (1-b) db  \notag \\
    &= \frac{1}{\theta}\left(\frac{3}{2}\right)^{\theta -1} \beta (1-\ln \beta).
\end{align}
Since $1<p\leq1+\tau_{1}$ and $\tau_{1} < \left(7QC^{*}(c,Q)\right)^{-1/2}$, it follows that $\alpha = C^{*}(c,Q) \frac{(p-1)^{2}}{p} \leq C^{*}(c,Q) \frac{\tau_{1}^{2}}{p} = \frac{1}{7pQ} < \frac{1}{3}$ and $\frac{\alpha}{2p-\alpha} <\frac{1}{5}$. Choose $\theta= \frac{2p}{2p-\alpha}$, (\ref{5.36}) can be reduced to
\begin{align}\label{5.37}
    \frac{\alpha}{2p}\int_{\left(1-\beta\right)}^{1}& \ln J(b) db  \notag  \\
    \leq & (p-\alpha) \ln C(c,Q) \beta + 2Q(p-\alpha)\ln (3/2)\beta + 4Q(p-\alpha)\beta \ln \frac{2p-\alpha}{\alpha} \\
    &+ \frac{2Q}{\theta}\left(\frac{3}{2}\right)^{\theta -1}(p-\alpha) \beta\left(1 -\ln \beta \right)-  \frac{p-\alpha}{p} \int_{1-\beta}^{\left(1-\beta\right)^{\frac{1}{\theta}}} \ln J(b)db  \notag
\end{align}
where we used (\ref{5.361}).
Let $g(b) = -\frac{\ln(1-b)}{b^{\frac{1}{\theta}}- b}$. It can be readily verified that $g(b)$ is monotonically increasing with respect to $b \in [2/3, 1]$. To prove that $J(2/3)$ is bounded, we distinguish between two cases.

Case 1: $\ln J(b) \geq g(b)$ for any $b \in [2/3, (1-\beta)^{1/\theta}]$. Then,
\begin{align}\label{5.38}
    \int_{1-\beta}^{\left(1-\beta\right)^{\frac{1}{\theta}}} \ln J(b)db  \geq  \int_{1-\beta}^{\left(1-\beta\right)^{\frac{1}{\theta}}} g(b) db  \geq g(1-\beta) \int_{1-\beta}^{\left(1-\beta\right)^{\frac{1}{\theta}}} db  = \ln \beta^{-1}
\end{align}
Combining (\ref{5.37}) and (\ref{5.38}), we have
\begin{align*}
    \frac{\alpha}{2p} \ln J(1-\beta) \beta
    \leq&  (p-\alpha) \ln C(c,Q) \beta + 2Q(p-\alpha)\ln (3/2)\beta + 4Q(p-\alpha)\beta \ln (2p-\alpha)  \\
    &+ 4Q(p-\alpha)\beta \ln (\alpha^{-1})+ \frac{2Q}{\theta}\left(\frac{3}{2}\right)^{\theta -1}(p-\alpha)\beta \left(1 +\ln \beta^{-1}\right)-  \frac{p-\alpha}{p}  \ln \beta ^{-1}
\end{align*}
Since $\beta \in (0,\frac{1}{3}]$, we can choose  $\beta = \alpha$. Hence, the above inequality yields that
\begin{align*}
    \frac{\alpha}{p} \ln J(1-\alpha)
     \leq & 2p \ln C(c,Q) +  4pQ\ln (3/2) +  8pQ\ln4 + 6pQ \\
    &+ \left[14Q\alpha  - \frac{2}{p}\right](p-\alpha)\ln \alpha^{-1} \\
     \leq & 4 \ln C(c,Q) +  8Q\ln (3/2) +  16Q\ln 4 + 12Q
\end{align*}
where we used $0<\alpha < \frac{1}{7pQ}$ and $1<p <2$. Therefore, from $J(s)$ is monotonically increasing with respect to $s \in [1/2,1]$, we have
\begin{align*}
     J(2/3)^{\frac{1}{p}} \leq  J(1-\alpha)^{\frac{1}{p} } \leq \left(C(c,Q)^{4} \left(3/2\right)^{8Q}4^{16Q}e^{12Q}\right)^{\frac{1}{\alpha}} \leq  C(c,Q)^{\frac{1}{\alpha}}.
\end{align*}

Case 2: there exists $b \in [2/3, (1-\beta)^{\frac{1}{\theta}}]$ such that  $\ln J(b) < g(b)$. Since both $J(b)$ and $g(b)$ are monotone increasing with respect to $b$, we obtain
\begin{align*}
    \ln J(2/3) \leq \ln J(b) < g(b) \leq g((1-\beta)^{\frac{1}{\theta}})
\end{align*}
Choose $\beta = \frac{1}{3}$, we get
\begin{align*}
    \ln J(2/3) &\leq g((2/3)^{\frac{1}{\theta}}) = \frac{\ln(1-(2/3)^{\frac{1}{\theta}})^{-1}}{\left(\frac{2}{3}\right)^{\frac{1}{\theta^{2}}} -\left(\frac{2}{3}\right)^{\frac{1}{\theta}}} \leq \frac{\ln(1-(2/3)^{\frac{5}{6}})^{-1}}{\ln(\frac{3}{2})(\frac{3}{2})^{\frac{1}{\theta}}\frac{\theta -1}{\theta^{2}} } \\
    &\leq  \left(\frac{6}{5}\right)^{2}\frac{\ln(1-(2/3)^{\frac{5}{6}})^{-1}}{\ln(3/2)} \frac{2p-\alpha}{\alpha} \leq \frac{C}{\alpha}
\end{align*}

Combining the above two cases, it follows that there exists a positive constant $C(c,Q)$ such that $J(2/3)^{\frac{1}{p}} \leq C(c,Q)^{\frac{1}{\alpha}}$. Hence, from (\ref{5.30}), we get
\begin{align}\label{5.39}
    \|w_{p}\|_{L^{\infty}(B_{R})} &\leq  C(c,Q)J\left(\frac{2}{3}\right)^{\frac{1}{p}} \left(\fint_{B_{2R}} w_{p}^{\alpha} dx\right)^{\frac{1}{\alpha}} \leq C(c,Q)^{\frac{1}{\alpha }}\left(\fint_{B_{2R}} w_{p}^{\alpha} dx\right)^{\frac{1}{\alpha}}.
\end{align}

Finally, we estimate $\| w_{p}^{-1} \|_{L^{\infty}(B_{R})} $. Define the sequences $\{\beta_i\}$, $\{s_i\}$ and $\{R_i\}$ by
\begin{align*}
    &\beta_{0} = -\frac{\alpha}{p}-1 ,\quad \beta_{i+1} = \frac{1}{2Q-1}+\beta_{i}\frac{2Q}{2Q-1};  \\
    &s_{0}= -\alpha, \quad s_{i} = (1+\beta_{i})p,  \quad s_{i+1} = ks_{i},\quad s_{i}= \kappa^{i}s_{0};\\
    &R_{0} = 2R,\quad R_{i} = R\left(1+ 2^{-i}\right).
\end{align*}
For any $i \geq 0$, $\beta_{i} < -1$. Using the inequality (\ref{5.23}) with $(R_{1}, R_{2}, \beta) = (R_{i+1}, R_{i},\beta_{i})$, we get
\begin{align*}
    \left(\fint_{B_{R_{i+1}}}w_{p}^{s_{i+1}} d x\right)^{\frac{1}{|s_{i+1}|}} \leq (C(c,Q))^{\frac{p}{|s_{i}|}}2^{\frac{p(i+2)}{|s_{i}|}} \left(\fint_{B_{R_{i}}}w_{p}^{s_{i}}\right)^{\frac{1}{|s_{i}|}}.
\end{align*}
Applying the Moser iteration with the exponents $\{s_{i}\}_{i \geq 0}$ yields
\begin{align*}
    \|w_{p}^{-1}\|_{L^{\infty}(B_{R})} &\leq (C(c,Q))^{\sum_{i = 0}^{\infty} \frac{p}{|s_{i}|}} 2^{\sum_{i = 0}^{\infty}\frac{p(i+2)}{s_{i}}} \left(\fint_{B_{2R}}w_{p}^{-\alpha}\right)^{\frac{1}{\alpha}} \\
    &\leq (C(c,Q))^{\frac{4Q}{\alpha}} 2^{\frac{4Q(2Q+1)}{\alpha}}\left(\fint_{B_{2R}}w_{p}^{-\alpha}\right)^{\frac{1}{\alpha}} \leq C(c,Q)^{\frac{1}{\alpha}}\left(\fint_{B_{2R}}w_{p}^{-\alpha}\right)^{\frac{1}{\alpha}}
\end{align*}
where we used $\sum_{i = 0}^{\infty} \frac{p}{|s_{i}|}  = \frac{2pQ}{\alpha} < \frac{4Q}{\alpha}$ and $2^{\sum_{i = 0}^{\infty}\frac{p(i+2)}{s_{i}}}  = 2^{\frac{2pQ(2Q+1)}{\alpha}} < 2^{\frac{4Q(2Q+1)}{\alpha}}$. Thus,
\begin{align}\label{5.40}
    \left( \fint_{B_{2R}}w_{p}^{-\alpha}\right)^{-\frac{1}{\alpha}}   \leq C(c,Q)^{\frac{1}{\alpha}}  \left(\esssup_{B_{R}} w_{p}^{-1}\right)^{-1}
    \leq \left(C(c,Q)\right)^{\frac{1}{\alpha}}\essinf_{B_{R}} w_{p}.
\end{align}
Combining (\ref{5.21}), (\ref{5.39}) and (\ref{5.40}), it follows that there exists a constant $\mathcal{C}^{h}(c,Q)$ such that
\begin{align*}
    \esssup_{B_{R}} w_{p} \leq  \left[\mathcal{C}^{h}(c,Q)\left(\frac{p-1}{2-p}\right)^{2}\right]^{1/\alpha}\essinf_{B_{R}} w_{p}
\end{align*}
where $\mathcal{C}^{h}(c,Q)$ depending only on $c$ and $Q$, and $\alpha = C^{*}(c,Q) \frac{(p-1)^{2}}{p}$.
\hfill${\square}$

Combining (\ref{5.12}) with Theorem \ref{thm5.5}, we obtain the following $C^{0}$-estimate for $w_{p}$.
\begin{proposition}\label{prop5.6}
    Let $1<p \leq 1+\tau_{1}$ and $w_{p}$ be a weak solution of the Dirichlet problem (\ref{5.4a})-(\ref{5.4b}). Then, there exists a constant $R_{0} >0$ such that for any $x_{0} \in \Omega$, $B_{R} = B(x_{0}, R)$ with $B_{4R} \subset \Omega$ and $R \leq R_{0}$,
    \begin{align*}
        \left(\mathcal{C}^{h}(c,Q)\frac{(p-1)^{2}}{(2-p)^{2}}\right)^{-1/\alpha} \leq w_{p}(x) \leq 1
    \end{align*}
    holds for any $x \in B_{R}$.
\end{proposition}

\subsection{ Refined Lipschitz estimate for solutions to (\ref{5.4a})}
\
\vglue-10pt
 \indent

This section is devoted to establishing a refined Lipschitz estimate for weak solutions to (\ref{5.4a}) by deriving an explicit expression for the constant $\mathcal{C}^{l}(p,c,Q)$ in (\ref{5.11}) in terms of $p$. The first step is to establish suitable Caccioppoli type inequalities, with the associated constants expressed explicitly in terms of $p$.

Recall the structure condition of equation (\ref{5.5a}). Let $z=\left(z_1, z_2\right) \in \mathbb{R}^2$ and $ a_i(z)=\left(\delta^2+|z|^2\right)^{\frac{p-2}{2}} z_i$. For $1 <p <2$,  there holds
\begin{align}
    |a_{i}(\nabla_{0} w_{p})| &\leq \left(\delta^{2} + |\nabla_{0} w_{p}|^{2}\right)^{\frac{p-1}{2}} \label{5.41} \\
    \sum_{i,j = 1}^{2}  \partial_{z_{j}}a_{i}(\nabla_{0} w_{p})\eta_{i}\eta_{j} &\geq (p-1)\left(\delta^{2} + |\nabla_{0} w_{p}|^{2}\right)^{\frac{p-2}{2}}|\eta|^{2} \label{5.42}\\
    |\partial_{z_{j}}a_{i}(\nabla_{0} w_{p})| &\leq \left(\delta^{2} + |\nabla_{0} w_{p}|^{2}\right)^{\frac{p-2}{2}} \label{5.43}
\end{align}

\begin{lemma}\cite{R2015} \label{lem5.7}
    The functions $w_{p}^1=X_1 w_{p}, w_{p}^2=X_2 w_{p}$ and $w_{p}^3=T w_{p}$ are weak solutions respectively of the following equations (in $D$):
    \begin{align}
        & \sum_{i=1}^2 X_i\left(\sum_{j=1}^2 \partial_{z_j} a_i\left(\nabla_{0} w_{p}\right) X_j w_{p}^1\right)+\sum_{i=1}^2 X_i\left(\partial_{z_2} a_i\left(\nabla_{0} w_{p}\right) T w_{p}\right)+T\left(a_2\left(\nabla_{0} w_{p}\right)\right)=0 \label{5.44}, \\
        & \sum_{i=1}^2 X_i\left(\sum_{j=1}^2 \partial_{z_j} a_i\left(\nabla_{0} w_{p}\right) X_j w_{p}^2\right)-\sum_{i=1}^2 X_i\left(\partial_{z_1} a_i\left(\nabla_{0} w_{p}\right) T w_{p}\right)-T\left(a_1\left(\nabla_{0} w_{p}\right)\right)=0 \label{5.45}, \\
        & \sum_{i=1}^2 X_i\left(\sum_{j=1}^2 \partial_{z_j} a_i\left(\nabla_{0} w_{p}\right) X_j w_{p}^3\right)=0 . \label{5.46}
    \end{align}
\end{lemma}

To simplify the notation, we will write $\Phi = \delta^{2} + |\nabla_{0} w_{p}|^{2}$ and $a_{i} = a_{i}(\nabla_{0} w_{p})$. For completeness, the proofs of Lemmas \ref{lem5.8}-\ref{lem5.11} are included in the appendix (see Section 9.2).
\begin{lemma}\label{lem5.8}
    Let $\beta \geq 0$ and $\eta \in C^{\infty}_{0}(D)$. Then
    \begin{align*}
        \int_{D}  \eta^{\beta+2}\Phi^{\frac{p-2}{2}} |Tw_{p}|^{\beta} |\nabla_{0}(Tw_{p})|^{2} \leq \frac{(\beta +2)^{2}}{(p-1)^{2}(\beta +1)^{2}} \int_{D}\eta^{\beta}\Phi^{\frac{p-2}{2}} |\nabla_{0} \eta|^{2}|Tw_{p}|^{\beta +2}.
    \end{align*}
\end{lemma}

\begin{lemma}\label{lem5.9}
    Let $\beta \geq 0$ and $\eta \in C^{\infty}_{0}(D)$. Then, there exist absolute constants $C_{3}$ and $C_{4}$ such that
    \begin{align*}
        \int_{\Omega}  \eta^{2}\Phi^{\frac{p-2+\beta}{2}} \left|\nabla^{2}_{0}w_{p}\right|^{2}
        \leq  C_{3}\frac{(\beta +1)}{(p-1)^{2}}K_{\eta}\int_{\Omega} \Phi^{\frac{p+\beta}{2}}
        + C_{4}\frac{(\beta +1)^{2}}{(p-1)^{2}} \int_{\Omega} \eta^{2}  \Phi^{\frac{p-2+\beta}{2}} |Tw_{p}|^{2}
    \end{align*}
    where $K_{\eta} = \left\|\eta T \eta\right\|_{L^{\infty}(D)} + \left\|\nabla_{0} \eta\right\|^{2}_{L^{\infty}(D)} $.
\end{lemma}

\begin{lemma}\label{lem5.10}
    Let $\beta \geq 2$ and $\eta \in C_0^{\infty}(D)$. Then, there exists an absolute constant $C$ such that
    \begin{align*}
        \int_{D}  \eta^{\beta+2}\Phi^{\frac{p-2}{2}} |Tw_{p}|^{\beta} |\nabla^{2}_{0}w_{p}|^{2}
        \leq \left(\frac{C(\beta +1)^{2}}{(p-1)^{4}}\right)^{\frac{\beta}{2}}\left\|\nabla_{0} \eta\right\|^{\beta}_{L^{\infty}(D)} \int_{D} \eta^{2}\Phi^{\frac{p-2+\beta}{2}}  |\nabla^{2}_{0}w_{p}|^{2}.
    \end{align*}
\end{lemma}

\begin{lemma}\label{lem5.11}
    Let $\beta \geq 2$ and $\eta \in C^{\infty}_{0}(D)$. Then, there exists an absolute constant $C$ such that
    \begin{align*}
        \int_{D} \eta^{2}  \Phi^{\frac{p-2+\beta}{2}} |\nabla^{2}_{0}w_{p}|^{2} dx \leq  C K_{\eta}\frac{(\beta +1)^{5}}{(p-1)^{8}} \int_{D}  \Phi^{\frac{p+\beta}{2}} dx
    \end{align*}
    where $K_{\eta} = \left\|\eta T \eta\right\|_{L^{\infty}(D)} + \left\|\nabla_{0} \eta\right\|^{2}_{L^{\infty}(D)} $.
\end{lemma}

\begin{theorem}\label{thm5.12}
    Let $1 < p < 2$ and $w_{p} \in HW^{1,p}(D)$ be a weak solution of (\ref{5.5a}) with $\delta >0$. Assume that $D$ satisfies (\ref{5.9}). Then, there exists a positive constant $\mathcal{C}^{l}(c,Q)$ such that
    \begin{align*}
        \left\|\nabla_{0} w_{p}\right\|_{L^{\infty}(B_{R})} \leq \mathcal{C}^{l}(c,Q)(p-1)^{-\frac{4Q}{p}}\left(\fint_{B_{2R}}  \left(\delta^{2} + |\nabla_{0} w_{p}|^{2}\right)^{\frac{p}{2}} dx \right)^{\frac{1}{p}}
    \end{align*}
    for every $CC$-ball $B_{R}$ with the concentric ball $B_{2R} \subset D$, and $\mathcal{C}^{l}(c,Q)$ independent of $p$.
\end{theorem}
\noindent {\bf Proof }:  Let $\Phi_{1} = \left(\delta^{2} + |\nabla_{0} w_{p}|^{2}\right)^{\frac{1}{2}}$ and $\eta \in C^{\infty}_{0} (D)$. A direct calculation yield
\begin{align*}
    \nabla_{0} (\eta^{2}\Phi_{1}^{\frac{p+\beta}{2}}) = 2\eta \nabla_{0} \eta\Phi_{1}^{\frac{p+\beta}{2}} + \frac{p+\beta}{2}\eta^{2}\Phi_{1}^{\frac{p+\beta -2}{2}}\nabla_{0} \Phi_{1}
\end{align*}
and
\begin{align*}
    |\nabla_{0} (\eta^{2}\Phi_{1}^{\frac{p+\beta}{2}})|^{2} &\leq 2\left(4\eta^{2} |\nabla_{0} \eta|^{2}\Phi_{1}^{p+\beta} + \frac{(p+\beta)^{2}}{4}\eta^{2} \Phi_{1}^{p+\beta -2} |\nabla_{0} \Phi_{1}|^{2}\right)\\
    &\leq 8\eta^{2} |\nabla_{0} \eta|^{2}\Phi_{1}^{p+\beta} + \frac{(p+\beta)^{2}}{2}\eta^{2} \Phi_{1}^{p+\beta -2} |\nabla^{2}_{0} w_{p}|^{2}.
\end{align*}
Then, from $1<p<2 \leq \beta$ and Lemma \ref{lem5.11}, we have
\begin{align}\label{5.47}
    \int_{D} |\nabla_{0} (\eta^{2}\Phi_{1}^{\frac{p+\beta}{2}})|^{2} &\leq   8 \int_{D} \eta^{2} |\nabla_{0} \eta|^{2}\Phi_{1}^{p+\beta} + \frac{(p+\beta)^{2}}{2} \int_{\Omega} \eta^{2} \Phi_{1}^{p+\beta -2} |\nabla^{2}_{0} w_{p}|^{2} \notag \\
    & \leq 8K_{\eta}\int_{spt(\eta)}\Phi_{1}^{p+\beta} + \frac{(p+\beta)^{2}}{2}\left(C K_{\eta}\frac{(\beta +1)^{5}}{(p-1)^{8}} \int_{spt(\eta)}  \Phi_{1}^{p+\beta}\right) \notag \\
    &\leq C K_{\eta}\frac{(p+\beta)^{7}}{(p-1)^{8}}\int_{spt(\eta)}  \Phi_{1}^{p+\beta}.
\end{align}
Let $R <R_{1} < R_{2} \leq 2R$. From \cite{FSS1997,GN1998}, there exists a cut-off function $\eta \in C^{\infty}_{0}(B_{R_{2}})$ satisfies $\eta \equiv 1$ on $B_{R_{1}}$, $|\nabla_{0} \eta| \leq \frac{C}{R_{2}-R_{1}}$ and $|T \eta| \leq \frac{C}{(R_{2}-R_{1})^{2}} $. Then, using (\ref{5.s}) for $\eta^{2} \Phi_{1}^{\frac{p+\beta}{2}}$ with $p=2$ and $\kappa=\frac{Q}{Q-2}$, we have
\begin{align*}
    \left( \fint_{B_{R_{2}}}|\eta^{2}\Phi_{1}^{\frac{p+\beta}{2}}|^{2\kappa} d x\right)^{\frac{1}{\kappa}} \leq \mathcal{C}^{2}_{s} R^{2}_{2} \fint_{B_{R_{2}}}\left|\nabla_{0} (\eta^{2}\Phi_{1}^{\frac{p+\beta}{2}})\right|^2 d x \leq  C(c,Q) R^{2}_{2} \fint_{B_{R_{2}}}\left|\nabla_{0} (\eta^{2}\Phi_{1}^{\frac{p+\beta}{2}})\right|^2 d x
\end{align*}
where the last inequality follows from Remark \ref{remS}. Applying Lemma \ref{lemD}, we get
\begin{align}\label{5.48}
    \fint_{B_{R_{1}}}\Phi_{1}^{(p+\beta)\kappa} d x \leq  \frac{|B_{R_{2}}|}{|B_{R_{1}}|}\left(\fint_{B_{R_{2}}}|\eta^{2}\Phi_{1}^{\frac{p+\beta}{2}}|^{2\kappa} d x\right) \leq C(c,Q) \left(\fint_{B_{R_{2}}}|\eta^{2}\Phi_{1}^{\frac{p+\beta}{2}}|^{2\kappa} d x\right).
\end{align}
Combining Lemma \ref{5.11}, (\ref{5.47}), (\ref{5.48}) and the properties of the cut-off function $\eta$, we obtain
\begin{align}\label{5.49}
    \left( \fint_{B_{R_{1}}}\Phi_{1}^{(p+\beta)\kappa} d x\right)^{\frac{1}{\kappa}} &\leq C(c,Q) R^{2}_{2} K_{\eta}\frac{(p+\beta )^{7}}{(p-1)^{8}}\left(\fint_{B_{R_{2}}}  \Phi_{1}^{p+\beta} dx \right) \notag \\
    & \leq C(c,Q) \left(\frac{R_{2}}{R_{2}-R_{1}}\right)^{2} \frac{(p+\beta)^{7}}{(p-1)^{8}}\left(\fint_{B_{R_{2}}}  \Phi_{1}^{p+\beta} dx \right).
\end{align}
Consider the following sequences of exponents and radii
\begin{align*}
    &\beta_{0} = 2,\quad \beta_{i+1} = (p+2)\kappa^{i} -p,\quad \beta_{i}+p = (p+2)\kappa^{i}; \\
    &s_{0}= p+2, \quad s_{i} = \beta_{i}+p, \quad s_{i+1} = \kappa s_{i}= \kappa^{i+1}s_{0};\\
    &R_{0} = 2R ,\quad R_{i} = R\left(1 + 2^{-i}\right).
\end{align*}
Then, $\left(\frac{R_{i}}{R_{i}-R_{i+1}} \right)^{2}\leq (2^{i+1}+2)^{2} \leq 4^{i+2}$. Using (\ref{5.49}) with $(R_{1}, R_{2}, \beta) = (R_{i+1}, R_{i}, \beta_{i})$, we have
\begin{align}\label{5.50}
    \left(\fint_{B_{R_{i+1}}}\Phi_{1}^{s_{i+1}} d x\right)^{\frac{1}{s_{i+1}}}&\leq \left[C(c,Q) \left(\frac{R_{i}}{R_{i}-R_{i+1}} \right)^{2} \frac{s_{i}^{7}}{(p-1)^{8}}\right]^{\frac{1}{s_{i}}}\left(\fint_{B_{R_{i}}}  \Phi_{1}^{s_{i}}\right)^{\frac{1}{s_{i}}} \notag \\
    &\leq \left(\frac{C(c,Q)4^{i+2} }{(p-1)^{8}}\right)^{\frac{1}{s_{i}}}  s_{i}^{\frac{7}{s_{i}}}\left(\fint_{B_{R_{i}}}  \Phi_{1}^{s_{i}}\right)^{\frac{1}{s_{i}}}.
\end{align}
Iterating with the exponents $\{s_i\}_{i \geq 0}$, we obtain
\begin{align*}
    \left\|\Phi_{1}\right\|_{L^{\infty}(B_{R})} &\leq \left(\frac{C(c,Q) }{(p-1)^{8}}\right)^{\sum_{i = 0}^{\infty}  \frac{1}{s_{i}}} 4^{\sum_{i = 0}^{\infty}\frac{i+2}{s_{i}}} \prod_{i=0}^{\infty}s_{i}^{\frac{7}{s_{i}}}\left(\fint_{B_{2R}}  \Phi_{1}^{p+2}\right)^{\frac{1}{p+2}} .
\end{align*}
A straightforward calculation leads to
\begin{align}\label{5.51}
    \sum_{i = 0}^{\infty}  \frac{1}{s_{i}} = \frac{Q}{2(p+2)},\; \sum_{i = 0}^{\infty}\frac{i+2}{s_{i}}  = \frac{Q(Q+2)}{4(p+2)},\; \prod_{i=0}^{\infty}s_{i}^{\frac{7}{s_{i}}} = (p+2)^{\frac{7Q}{2(p+2)}}\left(\frac{Q}{Q-2}\right)^{\frac{7Q(Q-2)}{4(p+2)}}.
\end{align}
Hence,
\begin{align} \label{5.52}
    \left\|\Phi_{1}\right\|_{L^{\infty}(B_{R})}  &\leq \left(\frac{C(c,Q) }{(p-1)^{8}}\right)^{\frac{Q}{2(p+2)}} 4^{\frac{Q(Q+2)}{4(p+2)}} (p+2)^{\frac{7Q}{2(p+2)}}\left(\frac{Q}{Q-2}\right)^{\frac{7Q(Q-2)}{4(p+2)}}\left(\fint_{B_{2R}}  \Phi_{1}^{p+2}\right)^{\frac{1}{p+2}} \notag \\
    &\leq C(c,Q) (p-1)^{-\frac{4Q}{p+2}}\left(\fint_{B_{2R}}  \Phi_{1}^{p+2}\right)^{\frac{1}{p+2}}.
\end{align}
To decrease the exponent on the right-hand side of (\ref{5.52}) from $p+2$ to $p$, we adopt the interpolation argument in Theorem 5.1 of \cite{R2015}. Let $0< \lambda <1$ and
\begin{align*}
    R_{0} = 2R, \qquad R_{i} = 2R (\lambda + (1-\lambda)2^{-i}).
\end{align*}
Then, we have
\begin{align*}
    \frac{R_{i}}{R_{i} - R_{i+1}} = \frac{\lambda + (1-\lambda)2^{-i}}{(1-\lambda)2^{-i-1}} = \frac{\lambda}{1-\lambda}2^{i+1} + 1 \leq (1-\lambda)^{-1}2^{i+2}.
\end{align*}
Combining (\ref{5.50}) and the above inequality, we obtain
\begin{align}\label{5.53}
    \left(\fint_{B_{2\lambda R}}\Phi_{1}^{s_{j}} d x\right)^{\frac{1}{s_{j}}}  & \leq  \left(\frac{4^{i+2}}{(1-\lambda )^{2}} \frac{C(c,Q)}{(p-1)^{8}}\right)^{\sum_{i= 0}^{j}  \frac{1}{s_{i}}}\prod_{i=0}^{j} (s_{i})^{\frac{7}{s_{i}}}
    \left(\fint_{B_{2R}} \Phi_{1}^{p+2}\right)^{\frac{1}{p+2}}.
\end{align}
For any $j>0$, one can easily check that
\begin{align}\label{5.54}
    \sum_{i= 0}^{j}  \frac{1}{s_{i}} < \sum_{i = 0}^{\infty}  \frac{1}{s_{i}},\quad \sum_{i= 0}^{j}  \frac{i+2}{s_{i}} < \sum_{i= 0}^{\infty}  \frac{i+2}{s_{i}}, \quad \prod_{i=0}^{j}s_{i}^{\frac{7}{s_{i}}} <  \prod_{i=0}^{\infty}s_{i}^{\frac{7}{s_{i}}}.
\end{align}
Combining (\ref{5.51}), (\ref{5.53}) and (\ref{5.54}), we find that for any $s> p+2$, there is
\begin{align}\label{5.55}
    \left(\int_{B_{2\lambda R}}\Phi_{1}^{s} d x\right)^{\frac{1}{s}}  &\leq C(c,Q) (1-\lambda)^{-\frac{Q}{p+2}}(p-1)^{-\frac{4Q}{p+2}}\left(\int_{B_{2R}}  \Phi_{1}^{p+2}\right)^{\frac{1}{p+2}} .
\end{align}
Choose $\epsilon$ such that $\frac{1}{p+2} = \frac{\epsilon }{ p} + \frac{1-\epsilon}{s}$. Let $\zeta = \frac{Q}{p+2} $ and $\lambda^{'} = \frac{1+\lambda}{2} \in (1/2,1)$, we have
\begin{align}\label{5.56}
    \frac{(1-\lambda)^{\frac{\zeta}{\epsilon}}}{\left(1-\lambda/\lambda^{'}\right)^{\zeta}} \leq 2^{\frac{\zeta}{\epsilon}}\left(1-\lambda^{'}\right)^{{\zeta \frac{1-\epsilon}{\epsilon}}}, \quad
    \frac{(1-\lambda)^{\frac{\zeta}{\epsilon}}}{\left(1-\lambda/\lambda^{'}\right)^{\zeta}\left(1-\lambda^{'}\right)^{{\zeta \frac{1-\epsilon}{\epsilon}}}} = (1+\lambda)^{\zeta} 2^{\zeta\frac{1-\epsilon}{\epsilon}} \leq 2^{\frac{\zeta}{\epsilon}}.
\end{align}
Let
\begin{align*}
    \psi  = \sup_{1/2 < \lambda <1} (1-\lambda)^{\zeta \frac{1-\epsilon}{\epsilon}}\left(\fint_{B_{2\lambda R}}\Phi_{1}^{p+2} d x\right)^{\frac{1}{p+2}} .
\end{align*}
From (\ref{5.55}) and (\ref{5.56}), we have
\begin{align}\label{5.57}
    (1-\lambda)^{\frac{\zeta}{\epsilon}}\left(\fint_{B_{2\lambda R}}\Phi_{1}^{s} d x\right)^{\frac{1}{s}}
    &\leq C(c,Q) (p-1)^{-\frac{4Q}{p+2}}\left(1-\lambda/\lambda^{'}\right)^{-\zeta} (1-\lambda)^{\frac{\zeta}{\epsilon}}\left(\fint_{B_{2\lambda^{'}R}}  \Phi_{1}^{p+2}\right)^{\frac{1}{p+2}} \notag \\
    &\leq C(c,Q) (p-1)^{-\frac{4Q}{p+2}}2^{\frac{\zeta}{\epsilon}}\psi.
\end{align}
By the definition of $\psi $, there exist $\lambda^{'} $ and $\varepsilon > 0$ such that
\begin{align*}
    \psi < \left(1-\lambda^{'}\right)^{{\zeta \frac{1-\epsilon}{\epsilon}}}\left(\fint_{B_{2\lambda^{'}R}}  \Phi_{1}^{p+2}\right)^{\frac{1}{p+2}} +\varepsilon
    \leq \left(1-\lambda^{'}\right)^{{\zeta \frac{1-\epsilon}{\epsilon}}}\left(\fint_{B_{2\lambda^{'}R}}  \Phi_{1}^{p} \right)^{\frac{\epsilon}{p}}\left(\fint_{B_{2\lambda^{'}R}}  \Phi_{1}^{s} \right)^{\frac{1-\epsilon}{s}} +\varepsilon
\end{align*}
where we used H$\ddot{\rm o} $lder’s inequality with exponent $\frac{p}{\epsilon(p+2)}$. Letting $\varepsilon \to 0 $ and applying  Young's inequality with a parameter $\theta>0$ to be chosen later we get
\begin{align*}
    \psi &\leq \theta (1-\epsilon) \left(1-\lambda^{'}\right)^{{\frac{\zeta }{\epsilon}}}\left(\fint_{B_{2\lambda^{'}R}}  \Phi_{1}^{s} \right)^{1/s}  + \epsilon \theta^{-\frac{1-\epsilon}{\epsilon}}\left(\fint_{B_{2\lambda^{'}R}}  \Phi_{1}^{p} \right)^{1/p}\\
    & \leq \theta (1-\epsilon) C(c,Q)(p-1)^{-\frac{4Q}{p+2}}2^{\frac{\zeta}{\epsilon}}\psi  + \epsilon \theta^{-\frac{1-\epsilon}{\epsilon}}\left(\fint_{B_{2\lambda^{'}R}}  \Phi_{1}^{p} \right)^{1/p}.
\end{align*}
where the second inequality follows from (\ref{5.57}). Choose $\theta = \left(2(1-\epsilon) C(c,Q) (p-1)^{-\frac{4Q}{p+2}}2^{\frac{\zeta}{\epsilon}}\right)^{-1}$, the above inequality yields that
\begin{align*}
    \psi  \leq  2\epsilon \left(2(1-\epsilon) C(c,Q) (p-1)^{-\frac{4Q}{p+2}}2^{\frac{\zeta}{\epsilon}}\right)^{\frac{1-\epsilon}{\epsilon}}\left(\fint_{B_{2\lambda^{'}R}}  \Phi_{1}^{p} \right)^{\frac{1}{p}}.
\end{align*}
Applying (\ref{5.57}) once again implies that
\begin{align*}
    (1-\lambda)^{\frac{\zeta}{\epsilon}}&\left(\fint_{B_{2\lambda R}}\Phi_{1}^{s} d x\right)^{\frac{1}{s}} \notag \\
    &\leq 2\epsilon \left(2(1-\epsilon)\right)^{\frac{1-\epsilon}{\epsilon}}\left(C(c,Q) (p-1)^{-\frac{4Q}{p+2}}2^{\frac{\zeta}{\epsilon}}\right)^{\frac{1}{\epsilon}}\left(\frac{|B_{2R}|}{|B_{2\lambda^{'}R}|}\right)^{\frac{1}{p}}\left(\fint_{B_{2R}}  \Phi_{1}^{p} \right)^{\frac{1}{p}}.
\end{align*}
Let $\lambda = 1/2$ and $s \to \infty$, then $\epsilon = \frac{p}{p+2}$, $\frac{\zeta}{\epsilon} = \frac{Q}{p}$ and
\begin{align*}
    \left\|\Phi_{1}\right\|_{L^{\infty}(B_{R})} &\leq 2^{\frac{Q}{p}}\frac{2p}{p+2} \left(\frac{4}{p+2}\right)^{\frac{2}{p}} \left(\frac{|B_{2R}|}{|B_{2\lambda^{'}R}|}\right)^{\frac{1}{p}}\left(C(c,Q)(p-1)^{-\frac{4Q}{p+2}} 2^{\frac{Q}{p}}\right)^{\frac{p+2}{p}} \left(\fint_{B_{2R}} \Phi_{1}^{p}\right)^{\frac{1}{p}}.
\end{align*}
From Lemma \ref{lemD}, $|B_{2\lambda^{'}R}| \geq c(\lambda^{'})^{Q} |B_{2R}| \geq c2^{-Q} |B_{2R}|$. Thus, the above inequality implies that there exists a positive constant $\mathcal{C}^{l}(c,Q)$ such that
\begin{align*}
    \left\|\Phi_{1}\right\|_{L^{\infty}(B_{R})} \leq \mathcal{C}^{l}(c,Q)(p-1)^{-\frac{4Q}{p}}\left(\fint_{B_{2R}}  \Phi_{1}^{p} \right)^{\frac{1}{p}}.
\end{align*}
where $\mathcal{C}^{l}(c,Q)$ independent of $p$. Since $\Phi_{1} = \left(\delta^{2} + |\nabla_{0} w_{p}|^{2}\right)^{\frac{1}{2}}$, we have
\begin{align*}
    \left\|\nabla_{0} w_{p}\right\|_{L^{\infty}(B_{R})} \leq \left\|\Phi_{1}\right\|_{L^{\infty}(B_{R})}\leq \mathcal{C}^{l}(c,Q)(p-1)^{-\frac{4Q}{p}}\left(\fint_{B_{2R}}  (\delta^{2} + |\nabla_{0} w_{p}|^{2})^{\frac{p}{2}} \right)^{\frac{1}{p}}.
\end{align*}
This completes the proof.
\hfill${\square}$

Using an approximation argument,  Theorem \ref{thm5.12} can be extended to general domain. The proof is given in the appendix.
\begin{theorem}\label{thm5.13}
    Let $1 <p <2$ and $\Omega$ be an open set of $\mathbb{H}^{1}$. Assume that $w_{p} \in HW^{1,p}(\Omega)$ is a weak solution of (\ref{5.5a}) with $\delta >0$. Then
    \begin{align*}
        \left\|\nabla_{0} w_{p}\right\|_{L^{\infty}(B_{R})} \leq \mathcal{C}^{l}(c,Q)(p-1)^{-\frac{4Q}{p}}\left(\fint_{B_{R}}  \left(\delta^{2} + |\nabla_{0} w_{p}|^{2}\right)^{\frac{p}{2}} \right)^{\frac{1}{p}}
    \end{align*}
    for every $CC$-ball $B_{R}$ such that the concentric ball $B_{2R} \subset \Omega$.
\end{theorem}

By a similar approximation argument, the conclusion of Theorem \ref{thm5.13} can be extended to weak solutions of (\ref{5.4a}). The details are presented in the appendix.
\begin{theorem}\label{thm5.14}(Refined Lipschitz estimate)
    Let $1 < p < 2$ and $w_{p} \in HW^{1,p}(\Omega)$ be a weak solution of (\ref{5.4a}). Then
    \begin{align*}
        \|\nabla_{0} w_{p}\|_{L^{\infty}(B_{R})} \leq \mathcal{C}^{l}(c,Q)(p-1)^{-\frac{4Q}{p}} \left(\fint_{B_{2R}}  |\nabla_{0} w_{p}|^{p} \right)^{\frac{1}{p}}
    \end{align*}
    for every $CC$-ball $B_{R}$ such that the concentric ball $B_{2R} \subset \Omega$.
\end{theorem}

\begin{proposition}\label{prop5.15}
    Let $1 < p < 2$ and $w_{p} \in HW^{1,p}(\Omega)$ be a weak solution of the Dirichlet problem (\ref{5.4a})-(\ref{5.4b}). Then, there exists a constant $\mathcal{\widetilde{C} }^{l}(c,Q)>0$ independent of $p$ such that
    \begin{align*}
        \|\nabla_{0} w_{p}\|_{L^{\infty}(B_{R/2})} \leq \frac{\mathcal{\widetilde{C} }^{l}(c,Q)}{(p-1)^{\frac{4Q}{p}} R}
    \end{align*}
    for every $CC$-ball $B_{R}$ with the concentric ball $B_{2R} \subset \Omega$.
\end{proposition}
\noindent {\bf Proof }: From  Lemma 3.2 of \cite{CDG1993}, there exists a cut-off function $\eta \in C_{0}^{\infty} (B_{2R})$, $\eta \equiv 1$ on $B_{R}$ and $|\nabla_{0} \eta| \leq \frac{C}{R}$, where $C$ independent of $R$. Let $\varphi = \eta^{p}w_{p}$ as test function on weak formula (\ref{5.6}) with $\delta =0$, we have
\begin{align*}
    \int_{B_{2R}} \eta^{p}|\nabla_{0} w_{p}|^{p} dx  &= -\int_{B_{2R}} p \eta^{p-1}|\nabla_{0} w_{p}|^{p-2} \left\langle \nabla_{0}w_{p},\nabla_{0} \eta \right\rangle  w_{p} dx \\
    &\leq \int_{B_{2R}} p\eta^{p-1}|\nabla_{0} \eta||\nabla_{0} w_{p}|^{p-1} w_{p}
\end{align*}
Using Young's inequality with exponent $\frac{p}{p-1}$, we get
\begin{align*}
    \int_{B_{2R}} \eta^{p}|\nabla_{0} w_{p}|^{p} \leq  \frac{p-1}{p}\int_{B_{2R}} \eta^{p}|\nabla_{0} w_{p}|^{p}  + p^{p-1}\int_{B_{2R}} |\nabla_{0} \eta|^{p}w_{p}^{p},
\end{align*}
i.e., $\int_{B_{2R}} \eta^{p}|\nabla_{0} w_{p}|^{p} \leq p^{p}\int_{B_{2R}} |\nabla_{0} \eta|^{p} w_{p}^{p}$. Therefore, using (\ref{5.12}), Lemma \ref{lemD} and the properties of $\eta$, we have
\begin{align}\label{5.58}
    \fint_{B_{R}} |\nabla_{0} w_{p}|^{p} &\leq \frac{1}{|B_{R}|} \int_{B_{2R}} \eta^{p}|\nabla_{0} w_{p}|^{p} \leq p^{p} \frac{1}{|B_{R}|} \int_{B_{2R}} |\nabla_{0} \eta|^{p} w_{p}^{p} \notag \\
    &\leq p^{p} \left(\frac{C}{R} \right)^{p}\frac{|B_{2R}|}{|B_{R}| } \leq C(c,Q) \left(\frac{p}{R}\right)^{p}
\end{align}
Combining (\ref{5.58}) and Theorem \ref{thm5.14}, we obtain
\begin{align*}
    \|\nabla_{0} w_{p}\|_{L^{\infty}(B_{R/2})} &\leq \frac{\mathcal{C}^{l}(c,Q)}{(p-1)^{\frac{4Q}{p}}}\left(\fint_{B_{R}}  |\nabla_{0} w_{p}|^{p} \right)^{\frac{1}{p}} \leq  \frac{\mathcal{C}^{l}(c,Q)}{(p-1)^{\frac{4Q}{p}}} \frac{C(c,Q)p}{R}
\end{align*}
Let $\mathcal{\widetilde{C} }^{l}(c,Q) = \max \{1, 2\mathcal{C}^{l}(c,Q)C(c,Q)\}$. Then the proof of Proposition \ref{prop5.15} is complete.
\hfill${\square}$

\subsection{ Limit of approximate solutions }
\
   \vglue-10pt
    \indent

Combining Propositions \ref{prop5.6} and \ref{prop5.15}, we can establish the uniform bound for $|\nabla_{0} u_{p}|$ that is independent of $p$.
\begin{theorem}\label{thm5.16}
    Let $u_{p} \in HW^{1,p}_{loc}(\Omega)$ be a weak solution of the Dirichlet problem (\ref{5.2a})-(\ref{5.2b}). Assume that $1 < p \leq 1+\tau$, where $\tau > 0$ is a constant defined by (\ref{5.59}). Then, there exist constants $R_{0}>0$ and $\mathcal{\widetilde{C} }(c,Q,R_{0})>0$ such that
    \begin{align*}
        |\nabla_{0} u_{p} (x)| \leq  \mathcal{\widetilde{C} }(c,Q,R_{0}), \qquad \forall x \in \Omega
    \end{align*}
    where $R_{0}$ and $\mathcal{\widetilde{C} }(c,Q,R_{0})$ independent of $p$.
\end{theorem}
\noindent {\bf Proof }:  Let $x_{0} \in \Omega$, $R>0$, $B_{R}= B(x_{0},R)$ and $B_{4R} \subset\subset \Omega$. Denote
\begin{align}\label{5.59}
  \mathcal{\widetilde{C}}(c,Q) = \left[\mathcal{\widetilde{C} }^{l}(c,Q) \mathcal{C}^{h}(c,Q)\right]^{-1},\quad \tau  = \min\left\{\tau_{1}, \frac{2\mathcal{\widetilde{C} }(c,Q) +1 - \sqrt{4\mathcal{\widetilde{C} }(c,Q)+1}}{2\mathcal{\widetilde{C} }(c,Q)}\right\}
\end{align}
where $\tau_{1} < \min\left\{1,  \left(7QC^{*}(c,Q)\right)^{-1/2}\right\}$ is the constant in Theorem \ref{thm5.5}. From Propositions \ref{prop5.6} and \ref{prop5.15}, we have
\begin{align*}
    \frac{|\nabla_{0}w_{p}(x)|}{w_{p}(x)} &\leq \frac{\|\nabla_{0} w_{p}\|_{L^{\infty}(B_{R/2})}}{w_{p}(x)} \leq \frac{\mathcal{\widetilde{C}}^{l}(c,Q)}{(p-1)^{\frac{4Q}{p}}R} \left[\mathcal{C}^{h}(c,Q)\frac{(p-1)^{2}}{(2-p)^{2}}\right]^{\frac{1}{\alpha}} \\
    &\leq \left[\frac{({\widetilde{C}}^{l}(c,Q))^{\alpha}\mathcal{C}^{h}(c,Q) \tau^{2-\frac{4Q\alpha}{p}}}{(1-\tau)^{2}}\right]^{\frac{1}{\alpha}} \frac{1}{R} \leq \frac{1}{R}
\end{align*}
holds for any $x \in B_{R/2}$. Hence,
\begin{align*}
    |\nabla_{0} u_{p}(x)|  = (p-1)\frac{|\nabla_{0}w_{p}(x)|}{w_{p}(x)} \leq \frac{\tau}{R},\qquad \forall x \in B_{R/2}.
\end{align*}
In particular, the above inequality implies that
\begin{align*}
    \lim_{d(x,x_{0})\to \infty} |\nabla_{0}u_{p}(x)|  = 0.
\end{align*}
Then, there exists a constant $R_{0} > 0$, independent of $p$, such that $|\nabla_{0}u_{p}(x)| \leq 1$ whenever $d(x,x_{0}) > R_{0} /2$. In summary,
\begin{equation*}
    \begin{cases}
    |\nabla_{0}u_{p}(x)| \leq \frac{\tau}{R_{0}}, \qquad &x\in B_{R_{0}/2},\\
    |\nabla_{0}u_{p}(x)| \leq 1, \qquad &x\in B_{R_{0}/2}^{c} \cap \Omega.
    \end{cases}
\end{equation*}
Let $\mathcal{\widetilde{C}}(c,Q,R_{0}) = \max\left\{1, \frac{\tau}{R_{0}}\right\}$, we obtain $|\nabla_{0}u_{p}(x)| \leq \mathcal{\widetilde{C}}(c,Q,R_{0})$ for any $x \in \Omega$.
\hfill${\square}$

\begin{proposition}\label{prop5.17}
    Let $\{ u_{p}\}_{1<p \leq 1+\tau}$ be a sequence of weak solutions of the Dirichlet equation (\ref{5.2a})-(\ref{5.2b}). Then, there exists a locally Lipschitz function $u \in HW^{1,\infty}_{loc}(\Omega)$ and a subsequence $\{p_{i}\} \to 1$ such that $u_{p_{i}} \to u$ locally uniformly in $\Omega$.
\end{proposition}
\noindent {\bf Proof }: Let $u_{p}=(1-p) \ln w_{p}$, where $w_{p} \in HW_{loc}^{1, p}\left(\Omega\right)$ is the weak solution to the Dirichlet problem (\ref{5.2a})-(\ref{5.2b}). Let $r >0$ and $\mathsf{B}_{r}(y) \subset \Omega^{c}$, then from (\ref{5.12}), we have
\begin{align*}
    0< u_{p}(x) \leq (Q-p)\ln \left(\frac{\rho (x,y)}{r}\right) < 3 \ln \left(\frac{\rho (x,y)}{r}\right)
\end{align*}
holds for any $x \in \Omega$. Thus, $u_{p}(x)$ is uniformly bounded with respect to $p$ on any bounded set $U \subset \Omega$. From Theorem \ref{thm5.16}, $|\nabla_{0} u_{p}|$ is uniformity bounded witeh respect to $p$ in $\Omega$. Therefore, $\{u_{p}\}_{1<p\leq1+\tau }$ is uniformly bounded in $HW^{1,\infty}_{loc}(\Omega)$. By Arzel$\grave{\rm a} $-Ascoli theorem, there exist a convergent subsequence $\{u_{p_{i}}\}$ and a locally Lipschitz function $u \in HW^{1,\infty}_{loc}(\Omega)$ such that $u_{p_{i}} \to u$ locally uniformly as $p_{i} \to 1$.
\hfill${\square}$

\begin{proposition}\label{prop5.18}
    $u = \lim_{p_{i} \to 1} u_{p_{i}}$ is a weak solution of (\ref{1.2}).
\end{proposition}
\noindent {\bf Proof }: Let $u_{p}$ and $u$ be the functions discussed in Proposition \ref{prop5.17}, and let $v\in HW^{1,p}_{loc}(\Omega)$ satisfy $\{v \neq u\} \subset K$, where $K \subset \Omega$ is a compact set. From Lemma \ref{lem5.1}, we have
\begin{align}\label{5.60}
    \mathcal{J}_{u_{p}}(u_{p}) \leq \mathcal{J}_{u_{p}}(v).
\end{align}
Let $\mathcal{A}_{K} = \{ \varphi \in HW_{loc}^{1,\infty}(\Omega),  K \subset \{\text {spt }\varphi \} \subset \subset \Omega, 0 \leq\varphi\leq 1 \}$ and $\eta \in \mathcal{A}_{K}$ be a cut-off function. Replacing $v$ by $\eta v+(1-\eta) u_{p}$ in (\ref{5.60}), we have
\begin{align}\label{5.61}
    \int_{\text {spt } \eta} & \left(\frac{1}{p}|\nabla_{0} u_{p}|^{p} + \eta \left(u_{p}-v\right) |\nabla_{0} u_{p}|^{p}\right) d x  \notag \\
    & \leq \frac{1}{p} \int_{\text {spt } \eta} | \eta \nabla_{0} v + (1-\eta) \nabla_{0} u_{p}+\left(v-u_{p}\right) \nabla_{0} \eta|^{p} d x  \notag  \\
    & \leq \frac{3^{p-1}}{p} \int_{\text {spt } \eta}\left(\eta^{p}|\nabla_{0} v|^{p}+(1-\eta)^{p}|\nabla_{0} u_{p}|^{p}+(v-u_{p})^{p}|\nabla_{0} \eta|^{p}\right) d x .
\end{align}
Also, replacing $v$ by $u$ in (\ref{5.61}) to obtain
\begin{align*}
    \int_{\text {spt } \eta} & \left(\frac{1}{p}|\nabla_{0} u_{p}|^{p} + \eta \left(u_{p}-u\right) |\nabla_{0} u_{p}|^{p}\right) d x  \\
    \leq
    & \frac{3^{p-1}}{p} \int_{\text {spt } \eta}\left(\eta^{p}|\nabla_{0} u|^{p}+(1-\eta)^{p}|\nabla_{0} u_{p}|^{p}+(u-u_{p})^{p}|\nabla_{0} \eta|^{p}\right) d x
\end{align*}
Since $0 < \eta \leq 1$, the above inequality reduces to
\begin{align}\label{5.62}
    \int_{\text {spt } \eta} \eta |\nabla_{0} u_{p}|^{p} dx \leq& 3^{p-1}\int_{\text {spt } \eta}  \eta |\nabla_{0} u|^{p}dx + p \int_{\text {spt } \eta}  \eta (u-u_{p})|\nabla_{0} u_{p}|^{p} dx  \\
    &+3^{p-1} \int_{\text {spt } \eta}  (u-u_{p})^{p} |\nabla_{0} \eta |^{p}dx + (3^{p-1} -1)\int_{\text {spt } \eta}  (1-\eta)|\nabla_{0} u_{p}|^{p} dx \notag
\end{align}
Now, we estimate the right-hand side of (\ref{5.62}) as $p \to 1$. From Theorem \ref{thm5.16} and Proposition \ref{prop5.17}, we have
\begin{align}\label{5.63}
    \limsup_{p \to 1} \int_{\text {spt } \eta} p\eta \left(u_{p}-u\right) |\nabla_{0} u_{p}|^{p} d x  \leq \limsup_{p \to 1} \mathcal{\widetilde{C}}(c,Q,R_{0}) \int_{\text {spt } \eta}  p|u_{p}-u| d x = 0
\end{align}
Similarly, for the third term on the right-hand side of (\ref{5.62}), we get
\begin{align}\label{5.64}
    \limsup_{p \to 1} 3^{p-1}\int_{\text {spt } \eta}(u-u_{p})^{p} |\nabla_{0} \eta|^{p} d x  =  0
\end{align}
Moreover, since $(3^{p-1} -1) \to 0$ as $p \to 1$, we have
\begin{align}\label{5.65}
    \lim_{p \to 1}(3^{p-1} -1)\int_{\text {spt } \eta}(1-\eta)  |\nabla_{0} u_{p}|^{p} dx = 0.
\end{align}
Plugging (\ref{5.63})-(\ref{5.65}) into (\ref{5.62}), we have
\begin{align}\label{5.66}
    \limsup_{p \to 1} \int_{\text {spt } \eta}  \eta |\nabla_{0} u_{p}|^{p}  d x  \leq  \int_{\text {spt } \eta} \eta|\nabla_{0} u| dx.
\end{align}
On the other hand, by Theorem \ref{thm5.16} and Proposition \ref{prop5.17}, it follows that $\nabla_0 u_p \stackrel{*}{\rightharpoonup} \nabla_0 u$ in $L^{\infty}\left(\Omega , \mathbb{R}^{2}\right) \subset L^{1}_{loc}\left(\Omega,\mathbb{R}^{2}\right)$. Hence, by the lower semi-continuity of the $L^{1}$ norm, we obtain
\begin{align}\label{5.67}
    \int_{\operatorname{spt}\eta} \eta\left|\nabla_0 u\right| \leq \liminf _{ p \to 1} \int_{\operatorname{spt}\eta} \eta\left|\nabla_0 u_p\right|.
\end{align}
Combining (\ref{5.66}) with (\ref{5.67}), we get
\begin{align}\label{5.671}
    \lim _{ p \to 1} \int_{\operatorname{spt}\eta} \eta  \left|\nabla_{0} u_p\right|^p = \int_{\operatorname{spt}\eta} \eta\left|\nabla_0 u\right| , \qquad \forall  \eta \in \mathcal{A}_{K}.
\end{align}

In the following, we fix $\eta \in C^{\infty}_{0}(\Omega)$ such that $\text {spt } \eta \subset \Omega$, $0 \leq \eta(u_{p} -v) \leq 1$ on $\Omega$ and $\eta \equiv 1$ on $K \subset \{\text {spt } \eta\}$.  Let $\eta_{p} = \eta (u_{p} -v)$. Then $\eta_{p} \in \mathcal{A}_{K}$, $\lim_{p \to 1} \eta_{p} = \eta (u -v)$ and $\eta (u -v) \in \mathcal{A}_{K}$. By (\ref{5.671}), we obtain
\begin{align*}
    \lim _{ p \to 1} \int_{\operatorname{spt}\eta} \eta_{p}  \left|\nabla_{0} u_p\right|^p &= \lim _{ p \to 1}  \left(\int_{\operatorname{spt}\eta} \eta(u-v)  \left|\nabla_{0} u_p\right|^p + \int_{\operatorname{spt}\eta} \eta(u_{p}-u)  \left|\nabla_{0} u_p\right|^p\right) \\
    &= \int_{\operatorname{spt}\eta} \eta(u-v)  \left|\nabla_{0} u\right|.
\end{align*}
Using the above equality and letting $p \rightarrow 1$ on both sides of (\ref{5.61}), we get
\begin{align*}
    \int_{K}  \left(|\nabla_{0} u| + u |\nabla_{0} u |\right) d x
    \leq
    \int_{K} \left(|\nabla_{0} v|+ v |\nabla_{0} u |\right) d x .
\end{align*}
where we used $\{u \neq v\} \subset K$ and $\eta \equiv 1$ on $K$.
\hfill${\square}$

The proof of existence of the weak solution in Theorem \ref{thm1.1} is complete.

\section{ Uniqueness of weak solutions }

\noindent In this section, we establish a comparison principle for the weak solution $E_{t}$ to (\ref{3.11}), which yields the uniqueness of weak solutions to the HIMCF. The comparison theorem is stated below.

\begin{theorem}\label{thm6.1}
    Assume that $\mathbb{H}^{1}$ has no compact component. Let $\{E_{t}\}_{t>0}$ and $\{F_{t}\}_{t>0}$ be weak solutions of (\ref{1.2}) in $ \Omega \subset \mathbb{H}^{1}$ and $E_{t}$ is precompact. If $E_{0} \subseteq F_{0}$, then $E_{t} \subseteq F_{t}$ holds for every time $t\geq 0$.
\end{theorem}

\noindent {\bf Proof }: Let $\Omega : = \mathbb{H}^{1} \setminus F_{0}$, $E_{t} = \{ u(x) < t\}$ and $F_{t} = \{ v(x) < t\}$. From Lemma \ref{lem3.3}, we know that $u$ and $v$ are weak solutions of (\ref{1.2}). Also, from Lemma \ref{lem4.9}, $\min \{v,t\}$ is a weak solution of (\ref{1.2}) in $\mathbb{H}^{1} \setminus (F_{0} \cup \partial F_{0})$. Suppose that $u > 0$, from $E_{0} \subseteq F_{0}$, we have $\min \{v,t\} < \frac{ u}{ 1- \varepsilon} $ for each $\varepsilon \in (0,1)$ near $\partial (E_{t}\setminus F_{0})$. Since $E_{t}$ is precompact, $E_{t}\setminus F_{0}$ is a precompact open set and $\left\{ \min \{v,t\} > \frac{ u}{ 1- \varepsilon}\right\} \subset \subset E_{t}\setminus F_{0}$.

Denote $v_{1} = \min \{v,t\}$, $u_{1} = \frac{ u}{ 1- \varepsilon}$ and $w_{1} = \frac{ w}{ 1- \varepsilon}$, where $w$ is a locally Lipschitz function satisfying $ u \leq w$. Since $u$ is a weak supersolution of (\ref{1.2}), it follows that
\begin{align*}
    \int | \nabla_{0} u|  + u | \nabla_{0} u| \leq  \int | \nabla_{0} w|  + w | \nabla_{0} u|.
\end{align*}
Hence,
\begin{align}\label{6.1}
    \int |\nabla_{0} u_{1}| + u_{1}|\nabla_{0} u_{1}| + \varepsilon   \int \left( w_{1}- u_{1}\right)|\nabla_{0} u_{1}| \leq \int |\nabla_{0} w_{1}| + w_{1}|\nabla_{0} u_{1}| .
\end{align}
Replacing $w_{1}$ by $u_{1}+\left(v_{1}-u_{1}\right)_{+}$ in (\ref{6.1}) to get
\begin{align}\label{6.2}
    \int_{\{v_{1}>u_{1}\}} |\nabla_{0} u_{1}| + u_{1}|\nabla_{0} u_{1}| + \varepsilon \int \left( v_{1}- u_{1}\right)|\nabla_{0} u_{1}| \leq \int_{\{v_{1}>u_{1}\}} |\nabla_{0} v_{1}| + v_{1}|\nabla_{0} u_{1}| .
\end{align}
Also, since $v_{1}$ is a weak solution of (\ref{1.2}), we have $J_{v_{1}}(v_{1}) \leq J_{v_{1}}(u_{1})$, i.e.,
\begin{align}\label{6.3}
    \int_{\{v_{1}>u_{1}\}} |\nabla_{0} v_{1}| + v_{1}|\nabla_{0} v_{1}| \leq \int_{\{v_{1}>u_{1}\}} |\nabla_{0} u_{1}| + u_{1}|\nabla_{0} v_{1}| .
\end{align}
Adding (\ref{6.2}) and (\ref{6.3}), we obtain
\begin{align}\label{6.4}
    \int_{\{v_{1}>u_{1}\}} (v_{1}-u_{1}) \left(|\nabla_{0} v_{1}| - |\nabla_{0} u_{1}|\right) + \varepsilon \int_{\{v_{1}>u_{1}\}} (v_{1}-u_{1})|\nabla_{0} u_{1}|  \leq 0 .
\end{align}
From (\ref{6.1}) and $u_{1} \leq w_{1}$, we know that $u_{1}$ is a weak supersolution of (\ref{1.2}), then $J_{u_{1}}(u_{1})\leq J_{u_{1}}(w_{1})$. Replace $w_{1}$ by $u_{1}+\left(v_{1}-s-u_{1}\right)_{+}$ with $s \geq 0$ yields
\begin{align*}
    \int_{\{v_{1}-s>u_{1}\}} |\nabla_{0} u_{1}| + u_{1}|\nabla_{0} u_{1}| \leq \int_{\{v_{1}-s>u_{1}\}} |\nabla_{0} v_{1}| + (v_{1}-s)|\nabla_{0} u_{1}| .
\end{align*}
Integrating over $s$ yields
\begin{align*}
    \int_{0}^{\infty} \int_{\{v_{1}-s>u_{1}\}} |\nabla_{0} u_{1}| + u_{1}|\nabla_{0} u_{1}| dx ds \leq \int_{0}^{\infty} \int_{\{v_{1}-s>u_{1}\}} |\nabla_{0} v_{1}| + (v_{1}-s)|\nabla_{0} u_{1}| dx ds ,
\end{align*}
i.e.,
\begin{align*}
    \int_{-\infty}^{+\infty} \int_{\mathbb{H}^{1}} \chi _{\{s>0\}} \chi _{\{v_{1}-s>u_{1}\}} (1+u_{1}-v_{1}+s)&|\nabla_{0} u_{1}|  dx ds \\
    & \leq  \int_{-\infty}^{+\infty} \int_{\mathbb{H}^{1}} \chi _{\{s>0\}} \chi _{\{v_{1}-s>u_{1}\}} |\nabla_{0} v_{1}| dx ds.
\end{align*}
Switching the order of integration, we have
\begin{align*}
    \int_{\mathbb{H}^{1}} |\nabla_{0} u_{1}| & \left(\int_{-\infty}^{+\infty} \chi _{\{s>0\}} \chi _{\{v_{1}-s>u_{1}\}} (1+u_{1}-v_{1}+s) ds\right) dx \\
    &\leq
    \int_{\mathbb{H}^{1}} |\nabla_{0} v_{1}| \left(\int_{-\infty}^{+\infty} \chi _{\{s>0\}} \chi _{\{v_{1}-s>u_{1}\}}  ds\right) dx.
\end{align*}
Then
\begin{align*}
    \int_{\mathbb{H}^{1}} - \frac{1}{2}\left(v_{1}-u_{1}\right)^{2} |\nabla_{0} u_{1}| dx\leq \int_{\mathbb{H}^{1}} (v_{1}-u_{1})(|\nabla_{0} v_{1}|-|\nabla_{0} u_{1}|) dx,
\end{align*}
which yields
\begin{align}\label{6.5}
    \int_{\{v_{1}>u_{1}\}} - \frac{1}{2}\left(v_{1}-u_{1}\right)^{2} |\nabla_{0} u_{1}| dx  \leq \int_{\{v_{1}>u_{1}\}} (v_{1}-u_{1})(|\nabla_{0} v_{1}|-|\nabla_{0} u_{1}|) dx .
\end{align}
Substituting (\ref{6.5}) into (\ref{6.4}), we have
\begin{align}\label{6.6}
    \int_{\{v_{1}>u_{1}\}} |\nabla_{0} u_{1}|(v_{1}-u_{1})\left[\varepsilon - \frac{1}{2} (v_{1}-u_{1})\right] dx \leq 0.
\end{align}

Case 1: if $\varepsilon - \frac{1}{2} (v_{1}-u_{1}) \geq 0$, then from (\ref{6.6}), we have $|\nabla_{0} u_{1}|=0$ a.e. on $\{v_{1}>u_{1}\}$. It follows from (\ref{6.4}) that $|\nabla_{0} v_{1}|=0$ a.e. on $\{v_{1}>u_{1}\}$. Hence, ($v_{1}-u_{1}$) is constant on each component of $\{v_{1}>u_{1}\}$, and $(v_{1}-u_{1})_{+}$ is a locally Lipschitz function. Since $\{v_{1}>u_{1}\} \subset \subset E_{t} \setminus F_{0}$ is precompact and $\Omega$ has no compact components, ($v_{1}-u_{1}$) should be zero. Hence, $\{ v_{1} >u_{1}\} = \emptyset$, i.e., $v_{1} \leq u_{1}$ on $E_{t}\setminus F_{0}$.

Case 2: if $\varepsilon - \frac{1}{2} (v_{1}-u_{1}) < 0$, i.e., $v_{1}>u_{1}+2\varepsilon$. Since $ \{v_{1}>u_{1}+2\varepsilon \}  \subset \{v_{1}>u_{1}\} \subset \subset E_{t}\setminus F_{0}$ is precompact, there exists a countable collection of non-empty sets $\{A_{k}\}_{k=1,\cdots,n}$ such that
\begin{align*}
    \{v_{1}>u_{1} + 2\varepsilon\}  = \bigcup ^{n}_{k=1}A_{k},\qquad
     A_{k}:=\{ 2k\varepsilon < v_{1} \leq  u_{1} + 2(k+1)\varepsilon\}.
\end{align*}
Hence, $v_{1} \leq u_{1} + 2(n+1)\varepsilon$. Let $v_{k} := v_{1} -2k \varepsilon$, then $v_{n} \leq u_{1} + 2\varepsilon$. Since each $v_{k}$ can be verified as a weak solution to (\ref{1.2}), it follows from Case 1 that $v_{n} \leq u_{1}$. This yields $\{v_{1} > u_{1} + 2\varepsilon\} = \bigcup_{k=1}^{n-1} A_{k}$. By iterating this procedure, we conclude that $v_{1} \leq u_{1}$ on $E_{t} \setminus F_{0}$.

In summary, for any $\varepsilon \in (0,1)$, we have $\min\{v ,t\} \leq \frac{u}{1-\varepsilon}$ on $E_{t} \setminus F_{0}$. Thus, $\min\{v ,t\} \leq u$ on $E_{t} \setminus F_{0}$. Since $u< t$ on $E_{t}\setminus F_{0}$, we obtain $v \leq u$ on $E_{t}\setminus F_{0}$. Consequently, $E_{t} \subseteq F_{t}$.
\hfill${\square}$

The uniqueness of the weak solution follows immediately from the comparison theorem. This completes the proof of Theorem \ref{thm1.1}.

\section{ Geometric Properties and application of the generalized evolution }

\noindent First, we identify suitable initial conditions that ensure the weak solution remains smooth for a short time. Next, we show that the $\mathbb{H}$-perimeter of the evolving surface increases exponentially along the generalized horizontal inverse mean curvature flow. Finally, as an application of the HIMCF, we provide a proof of the Heintze-Karcher type inequality.

\subsection{ Smoothness of the weak evolution  }
\
   \vglue-10pt
    \indent

Before proving Corollary \ref{cor1.2}, we introduce a property of classical solutions to the HIMCF.
\begin{lemma}\label{lem7.1}
    Let $E_{t} = \{x \in \mathbb{H}^{1} : u(x) < t\}$ and $M_{t} =\partial E_{t} $. Assume that $\{M_{t}\}_{a\leq t <b}$ is a family smooth classical solution of (\ref{1.2}) and $H_{0}|_{M_{t}}>0$. Then, for any $t \in [a,b)$, $E_{t}$ minimizes the functional (\ref{1.6}) on $E_{b} \setminus E_{a}$.
\end{lemma}
\noindent {\bf Proof}: Since the horizontal mean curvature $H_{0} |_{M_{t}} = |\nabla_{0} u| >0 $, the horizontal unit outer normal vector is given by $\nu_{0} = \frac{\nabla_{0} u}{|\nabla_{0} u|}$, which is a smooth vector field on $E_{b} \setminus E_{a}$. For any set $F \subset E_{b} \setminus E_{a}$ with locally finite $\mathbb{H}$-perimeter and any compact set $K \subset E_{b} \setminus E_{a}$ containing $F \Delta E_{t}$, we have
\begin{align*}
    J_{u}(E_{t}) &= P_{\mathbb{H}}(E_{t}, K) - \int_{E_{t} \cap K} |\nabla_{0} u| \\
    &= \int_{\partial E_{t} \cap K} \left\langle \nu_{\partial E_{t}}, \nu_{0}\right\rangle  - \int_{E_{t} \cap K} |\nabla_{0} u|= \int_{\partial F \cap K} \left\langle \nu_{\partial F}, \nu_{0}\right\rangle  - \int_{F \cap K} |\nabla_{0} u| \\
    & \leq  P_{\mathbb{H}}(F, K) - \int_{F\cap K} |\nabla_{0} u| =  J_{u}(F),
\end{align*}
where in the third equality we use Theorem 2.1 in \cite{ACV2007} and  $\nu_{0}$  is regarded as calibration.
\hfill${\square}$

\noindent {\bf Proof of Corollary \ref{cor1.2}}:
  From Lemma \ref{lem2.7}, the initial problem (\ref{0.1}) of $u(x)$ can be expressed as
\begin{equation}\label{7.1}
    \begin{cases}
        \frac{\partial}{\partial t} v = - |\nabla_{0} v|^{2} Tr\left[ \left(I - \frac{\nabla_{0} v \otimes \nabla_{0} v }{|\nabla_{0} v|^{2}}\right) \left(\nabla_{0}^{2} v\right)^{*}\right]^{-1},  \qquad x \in M_{t} \setminus \Sigma (M_{t}),  \\
        v(x,0) = 0,  \qquad x \in M_{0},
    \end{cases}
\end{equation}
where $v(x,t) = u(x) - t$, and both $\nabla_{0} v$ and $\left(\nabla_{0}^{2} v\right)^{*}$ are non-degenerate at regular points.

Let $Dv(x)$ and $D^{2}v(x)$ be the gradient and the Hessian of $v$ in $\mathbb{R}^{3}$. As in \cite{NFM2010}, the equation (\ref{7.1}) can be rewritten as
\begin{align}\label{7.2}
    \frac{\partial}{\partial t} v &= - |\alpha(x)Dv |^{2} Tr\left[ \left(I - \frac{\alpha(x) Dv \otimes \alpha(x) Dv }{|\alpha(x) Dv|^{2}}\right) \left(\alpha(x) D^{2}v \alpha^{T}(x) \right)\right]^{-1} \\
    &:=F(x, Dv, D^{2}v)  \notag
\end{align}
where $\alpha(x) = [X_{1}(x),X_{2}(x)]^{T}$ and $|\alpha(x)Dv |^{2} = \left\langle \alpha(x)Dv, \alpha(x)Dv\right\rangle_{\mathbb{R}^{3}} $. Combining $\Sigma(M_{0}) = \emptyset$, $H_{0}\big|_{M_{0}} > 0$ and the bottom of page 7 of \cite{NFM2010}, we can readily derive that $\frac{\partial F}{\partial (D_{ij}v)} > 0$, i.e., (\ref{7.2}) is a parabolic equation in the regular region.

Using the classical theory of parabolic equations in \cite{LSU1968}, we know that there exists a classical solution of (\ref{0.1}) for a short time. According to Lemma \ref{lem7.1}, there exists a smooth solution $E_{t}^{*}$ of the initial problem (\ref{3.11}) on an open neighborhood $U$ of $E_{0}$. From Proposition \ref{prop4.8}, we have $E_{0} = E_{0}^{'} = E_{0}^{+}$. Since $E_{t}$ is precompact, there exists a sequence $\{t_{i}\}$ such that
\begin{align*}
    \int_{\mathbb{H}^{1}} |\chi_{E_{t_{i}}} - \chi_{E_{0}}| d\mathcal{H}^{2} \to 0, \qquad \text{ as } t_{i} \to 0,
\end{align*}
which implies that $E_{t}$ is precompact on $U$ for a short time. Hence, using Theorem \ref{thm6.1}, we know that $E_{t}$ is coincident with $E_{t}^{*}$ for a short time.
\hfill${\square}$

\subsection{ Exponential expansion of $\mathbb{H}$-perimeter }
\
   \vglue-10pt
    \indent

\noindent {\bf Proof of Corollary \ref{cor1.3}}: It follows from Lemma \ref{lem2.18} that
\begin{align*}
    J_{u}(E_{t}) = P_{\mathbb{H}}(E_{t}) - \int_{E_{t}} |\nabla_{0} u| dx  =  P_{\mathbb{H}}(E_{t}) - \int_{0}^{t}  P_{\mathbb{H}}(E_{s}) ds .
\end{align*}
Since $E_{t}$ is a minimizer of the functional (\ref{1.6}) for any $t>0$, we have
\begin{align*}
    0 = \frac{d}{dt}J_{u}(E_{t}) = \frac{d}{dt} P_{\mathbb{H}}(E_{t}) - P_{\mathbb{H}}(E_{t}).
\end{align*}
Thus,
\begin{align}\label{7.3}
    P_{\mathbb{H}}(E_{t}) = e^{t} P_{\mathbb{H}}(E_{0}), \qquad  \text{ for any }  t > 0.
\end{align}
According to Proposition \ref{prop4.8} (4), the equality $P_{\mathbb{H}}(E_{t}) = P_{\mathbb{H}}(E^{+}_{t})$ holds for all $t \geq 0$, as long as $E_{t}^{+}$ is precompact and $E_{0}$ is a horizontal minimizing hull. Therefore
\begin{align}\label{7.4}
    P_{\mathbb{H}}(E_{t}) = P_{\mathbb{H}}(E^{+}_{t}) = e^{t} P_{\mathbb{H}}(E^{+}_{0}) = e^{t} P_{\mathbb{H}}(E_{0}), \qquad  \text{ for any }  t\geq 0.
\end{align}
This completes the proof.
\hfill${\square}$

\subsection{ Heintze-Karcher type inequality  }
\
   \vglue-10pt
    \indent

\noindent {\bf Proof of Theorem \ref{thm1.4}} : Let $ \lambda(t)=  e^{\frac{1}{3}t}$ and $G_{\lambda(t)} = \delta_{\lambda(t)}(E_{0}) = \{\delta_{\lambda(t)}(x), x\in E_{0} \} \subset \mathbb{H}^{1}$. Then, $\lambda(0) =1$ and $G_{\lambda(0)} =  \delta_{1}(E_{0}) = E_{0}$. From Lemma \ref{lem2.17}, we have
\begin{align}\label{7.5}
    P_{\mathbb{H}}(G_{\lambda(t)}) = e^{t}P_{\mathbb{H}}(E_{0}),\qquad  |G_{\lambda(t)}| = e^{\frac{4}{3} t} |E_{0}|, \qquad \forall t\in [0,\infty).
\end{align}
The subsequent proof is carried out on a finite time interval, which suffices to establish the theorem. It is clear that $G_{\lambda(t)}$ is an open and bounded set whenever $t \in [0, T_{0}]$ and $T_{0} < \infty$. In fact, $i$): it can be proved that $G_{\lambda(t)}$ is an open set. As $E_{0}$ is an open set, then $\forall x \in E_{0}$, there exists a Kor$\acute{ \rm a}$nyi ball $ \mathsf{B}_{r}(x) \subset E_{0}$, where $r > 0$. Let $y = \delta_{\lambda(t)}(x) \in G_{\lambda(t)}$. For any point $\tilde{y}  \in \mathsf{B}_{\lambda(t)r }(y)$, it holds that $\rho(y, \tilde{y}) <  \lambda(t)r $. Applying (\ref{2.2}), we have $\rho (\delta_{\lambda^{-1}}(y) ,\delta_{\lambda^{-1}}(\tilde{y}))  < r$. This implies $\rho (x ,\delta_{\lambda^{-1}}(\tilde{y}))  < r$, so $\delta_{\lambda^{-1}}(\tilde{y}) \in \mathsf{B}_{r}(x) \subset E_{0}$. Consequently, $\tilde{y}  \in \delta_{\lambda(t)}(E_{0}) = G_{\lambda(t)}$. Therefore, we have $\mathsf{B}_{\lambda(t) r}(y) \subset G_{\lambda(t)}$, which implies that $G_{\lambda(t)}$ is an open set.

$ii$): it can be demonstrated that $G_{\lambda(t)}$ is both bounded and has locally finite $\mathbb{H}$-perimeter. Since $E_{0}$ is bounded, there exist a point $x_{0} \in E_{0}$ and a finite radius $R > 0$ such that $E_{0} \subseteq \mathsf{B}_{R}(x_{0})$. Thus, for any $x \in E_{0}$, $\rho(x, x_{0}) \leq R$. Using (\ref{2.2}), we obtain
\begin{align*}
    \rho (\delta_{\lambda(t)}(x), \delta_{\lambda(t)}(x_{0})) \leq \lambda(t) R \leq e^{\frac{1}{3}T_{0}} R < \infty.
\end{align*}
Denote $y_{0} = \delta_{\lambda(t)}(x_{0})$, then $y_{0} \in G_{\lambda(t)}$ and $\forall y \in G_{\lambda(t)}$, $\rho (y,y_{0}) \leq \lambda(t) R < \infty$. Hence,
\begin{align*}
    G_{\lambda(t)} \subseteq  \mathsf{B}_{\lambda(t) R}(y_{0}) \subseteq \mathsf{B}_{e^{\frac{1}{3}T_{0}}  R}(y_{0}).
\end{align*}
Therefore, $G_{\lambda(t)}$ is bounded whenever $t \in [0,T_{0}]$. Moreover, from (\ref{7.5}), we have
\begin{align*}
    P_{\mathbb{H}}(G_{\lambda(t)}) = e^{t}P_{\mathbb{H}}(E_{0}) \leq e^{T_{0}}P_{\mathbb{H}}(E_{0}) < \infty.
\end{align*}

Let $E_{t}$ be a weak solution to (\ref{0.1}) with initial data $E_{0}$. Then, for any measurable set $F$ with locally finite $\mathbb{H}$-perimeter and $F\Delta E_{t} \subset \subset \mathbb{H}^{1} \setminus E_{0}$, there holds
\begin{align}\label{7.6}
    J_{u}(E_{t}) \leq J_{u}(F).
\end{align}
Moreover, from (\ref{7.3}), we have
\begin{align}\label{7.7}
  P_{\mathbb{H}}(E_{t})= e^{t}P_{\mathbb{H}}(E_{0}),\qquad \forall t \in (0, \infty).
\end{align}
Combining (\ref{7.5}) with (\ref{7.7}), we obtain
\begin{align}\label{7.8}
    P_{\mathbb{H}}(G_{\lambda(t)})= P_{\mathbb{H}}(E_{t}),\qquad \forall t \in [0, \infty).
\end{align}
Next, we prove $|E_{t}| \geq |G_{\lambda(t)}|$ holds for any $t \in [0,T_{0}]$.  Since $P_{\mathbb{H}}(G_{\lambda(t)}) < \infty$ and $G_{\lambda(t)} \Delta E_{t} \subset \subset \mathbb{H}^{1} $, we may choose $F = G_{\lambda(t)}$ in (\ref{7.6}) to obtain
\begin{align*}
     P_{\mathbb{H}}(E_{t}) - \int_{E_{t}} |\nabla_{0} u|dx \leq P_{\mathbb{H}}(G_{\lambda(t)}) - \int_{G_{\lambda(t)}} |\nabla_{0} u|dx  .
\end{align*}
Thus, for any $t \in [0,T_{0}]$,
\begin{align*}
    \int_{G_{\lambda(t)} \setminus E_{t}} |\nabla_{0} u|dx \leq P_{\mathbb{H}}(G_{\lambda(t)}) - P_{\mathbb{H}}(E_{t}) = 0,
\end{align*}
where the last equality utilises (\ref{7.8}). From the above inequality, either $|G_{\lambda(t)} \setminus E_{t}| = 0$ or $|G_{\lambda(t)} \setminus E_{t}| >  0$ and $|\nabla_{0} u (x) |= 0$ almost everywhere on $G_{\lambda} \setminus E_{t}$.

Case 1 : if $|G_{\lambda(t)} \setminus E_{t}| = 0$, then $|E_{t}| \geq |G_{\lambda(t)} |$.

Case 2 : if $|G_{\lambda(t)} \setminus E_{t}| > 0$, we denote
\begin{align*}
    A_{1} &= \{x\in (G_{\lambda(t)} \setminus E_{t}) \setminus \partial E_{t},  |\nabla_{0} u (x)| = 0\},\\
    A_{2} &= \{x\in (G_{\lambda(t)} \setminus E_{t}) \cap \partial E_{t},   |\nabla_{0} u (x)| = 0\},\\
    A_{3} &= \{x\in G_{\lambda(t)} \setminus E_{t} ,  |\nabla_{0} u (x)| \neq  0\}.
\end{align*}
Then, $G_{\lambda(t)} \setminus E_{t} = A_{1} \cup A_{2} \cup A_{3}$, $A_{1} \cap A_{2} = \emptyset$, $A_{2} \cap A_{3} = \emptyset$, and $A_{1} \cap A_{3} = \emptyset$. Since $|\nabla_{0} u (x) |= 0$ almost everywhere on $G_{\lambda(t)} \setminus E_{t}$, we have $|A_{3}| = 0$. Furthermore, $A_{2} \subseteq \Sigma (\partial E_{t})  = \{x \in \partial E_{t}, |\nabla_{0} u(x)| = 0\}$. By Remark \ref{rem2.6}, we have $P_{\mathbb{H}}(A_{2}) = 0$ and $|A_{2}| = 0$. Hence, $ 0 < |G_{\lambda(t)} \setminus E_{t}| = |A_{1}|$.

From Lemma \ref{lem2.15} and $(G_{\lambda(t)} \setminus E_{t}) \cap E_{t} = \emptyset $, it follows that
\begin{align*}
    P_{\mathbb{H}}(G_{\lambda(t)}) &=  P_{\mathbb{H}}( E_{t}) + P_{\mathbb{H}}(G_{\lambda(t)} \setminus E_{t}) = P_{\mathbb{H}}( E_{t}) + P_{\mathbb{H}}(A_{1} \cup A_{2} \cup A_{3}) \\
    &=P_{\mathbb{H}}( E_{t}) + P_{\mathbb{H}}(A_{1}) +P_{\mathbb{H}}(A_{2}) + P_{\mathbb{H}}(A_{3})\\
    &\geq P_{\mathbb{H}}( E_{t}) + P_{\mathbb{H}}(A_{1}).
\end{align*}
Thus, by (\ref{7.8}), we have $P_{\mathbb{H}}(A_{1}) \leq P_{\mathbb{H}}(G_{\lambda(t)}) - P_{\mathbb{H}}( E_{t}) = 0$. By the isoperimetric inequality in $\mathbb{H}^{1}$ (cf. Theorem 7.1 in \cite{CDPT2007}), there exists a constant $C>0$ such that $ |A_{1}|^{\frac{3}{4}} \leq C P_{\mathbb{H}}(A_{1}) = 0$. This contradicts $|A_{1}| > 0$.

In summary, $|G_{\lambda (t)} \setminus E_{t}| = 0 $ holds for all $t \in [0,T_{0}]$, that is, $ |E_{t}| \geq |G_{\lambda (t)}|$.  From (\ref{7.5}), we have $ |G_{\lambda(t)} | = e^{\frac{4}{3} t} |E_{0}|$. Thus, $|E_{t}| \geq  e^{\frac{4}{3} t} |E_{0}|$. Let $g(t) = e^{-\frac{4}{3} t}|E_{t}|$. From the inequality $e^{-\frac{4}{3} t}|E_{t}| \geq  |E_{0}|$, it follows that $g(t) \geq g(0)$ holds for all $t \in [0, T_{0}]$. Differentiating $g(t)$ yields
\begin{align*}
    \frac{d}{dt} g(t) = e^{-\frac{4}{3} t} \left(\frac{\partial }{\partial t} |E_{t}| - \frac{4}{3}|E_{t}|\right).
\end{align*}
We claim that $\frac{d}{dt} g \Big|_{t=0} \geq 0$. If $\frac{d}{dt} g \Big|_{t=0} < 0$, then there exists $\varepsilon >0$ such that $g(t) < g(0)$ for all $t \in (0, \varepsilon)$, which contradicts the fact that $g(t) \geq g(0)$ for all $t \in [0,T_{0}]$. Therefore
\begin{align}\label{7.9}
    \frac{d}{dt} g \Big|_{t=0} = \left(\frac{\partial }{\partial t} |E_{t}| \Big|_{t=0}- \frac{4}{3}|E_{0}|\right) \geq 0.
\end{align}
It is known that $|E_{t}| = \int_{\{0 < u(x) <t\}} dx = \int_{0}^{t} \int_{M_{s}} d\sigma ds  $, where $d\sigma$ denotes the surface measure on $M_{s}$. From Remark \ref{rem2.6} and (\ref{2.6})-(\ref{2.8}), we obtain
\begin{align*}
    |E_{t}|= \int_{0}^{t} \int_{ M_{s} \setminus \Sigma(M_{s})}d\sigma ds = \int_{0}^{t} \int_{M_{s} \setminus \Sigma(M_{s})} |\nabla_{0} u|^{-1} d\sigma_{\mathbb{H}} ds =  \int_{0}^{t} \int_{M_{s} \setminus \Sigma(M_{s})} H_{0}^{-1} d\sigma_{\mathbb{H}} ds.
\end{align*}
Consequently,
\begin{align}\label{7.10}
    \frac{\partial }{\partial t} |E_{t}|\Big|_{t=0} = \int_{M_{0} \setminus \Sigma (M_{0})} H_{0}^{-1} d\sigma_{\mathbb{H}} .
\end{align}
In view of (\ref{7.9}) and (\ref{7.10}), we establish (\ref{1.7}).

If ${M}_{0}$ is a horizontal constant mean curvature surface, i.e., $H_{0}$ is constant along the noncharacteristic locus, then from Corollary \ref{cor1.8} (whose proof does not rely on Theorem 1.4), we have
\begin{align*}
    \frac{4}{3}|E_{0}| = H_{0}^{-1} P_{\mathbb{H}}(E_{0}) = \int_{M_{0} \setminus \Sigma(M_{0})} H_{0}^{-1} d\sigma_{\mathbb{H}}.
\end{align*}
Therefore, the equality of (\ref{1.7}) holds for the horizontal constant mean curvature surface.
\hfill${\square}$

\section{ $\mathbb{H}$-perimeter preserving flow and its application}
\noindent In this section, we construct a $\mathbb{H}$-perimeter preserving flow using a Heisenberg dilation of the HIMCF. Moreover, the geometric properties of the $\mathbb{H}$-perimeter preserving flow are used to establish the Minkowski type formula (\ref{1.9}).

\subsection{ $\mathbb{H}$-perimeter preserving flow }
\
   \vglue-10pt
    \indent

Let $\widehat{E}_{t} = \delta_{\lambda}(E_{t}) = \{y \in \mathbb{H}^{1}, y=\delta_{\lambda}(x), x \in E_{t} \}$ and $\widehat{M}_{t} = \partial \widehat{E}_{t}$, where $E_{t}$ is the weak solution of the HIMCF (\ref{1.2}).

\begin{lemma}\label{lem8.1}
    Let $\widehat{V}(y)$ and $\widehat{\nu}_{0}(y)$ be the non-unit outer normal vector and the horizontal unit outer normal vector at $y \in \widehat{M}_{t} \setminus \Sigma(\widehat{M}_{t})$ respectively. Then
    \begin{align}\label{8.1}
        \left\langle  \widehat{\nu}_{0}(y), \widehat{V}(y) \right\rangle = \lambda^{3} \left\langle  \nu_{0}(x), V(x) \right\rangle .
    \end{align}
    where $V(x)$ and $\nu_{0}(x)$ are the non-unit outer normal vector and the horizontal unit outer normal vector at $x = \delta_{\lambda^{-1}}(y) \in M_{t} \setminus \Sigma(M_{t})$ respectively.
\end{lemma}
\noindent {\bf Proof}: It is convenient to introduce the following notations :
\begin{align*}
    p=\left\langle V(x), X_{1}\right\rangle ,\quad q=\left\langle V(x), X_{2}\right\rangle, \quad \omega =\left\langle V(x), T \right\rangle.
\end{align*}
Then,
\begin{align}\label{8.2}
    V(x) = pX_{1} + qX_{2}+ \omega T, \quad V_{0}(x) = pX_{1} + qX_{2}, \quad |V_{0}(x)| = \sqrt{p^{2}+ q^{2}}.
\end{align}
Similarly, let
\begin{align*}
    \widehat{p}=\left\langle \widehat{V}(y), X_{1}\right\rangle ,\quad \widehat{q}=\left\langle \widehat{V}(y), X_{2}\right\rangle, \quad \widehat{\omega} =\left\langle \widehat{V}(y), T \right\rangle.
\end{align*}
Then,
\begin{align}\label{8.3}
    \widehat{V}(y)=  \widehat{p}X_{1} + \widehat{q}X_{2} + \widehat{\omega}T,\quad \widehat{V}_{0}(y) = \widehat{p} X_{1} + \widehat{q}X_{2}, \quad |\widehat{V}_{0}(y)| = \sqrt{\widehat{p}^{2}+ \widehat{q}^{2}}.
\end{align}
Given an open set $\mathcal{O} \subset \mathbb{R}^{2}$, let $\vartheta : \mathcal{O} \to \mathbb{H}^{1}$ be a $C^{2}$ parametrization of $M_{t}$, then
\begin{align*}
    M_{t} = \left\{\vartheta(a,b) =x_{1}(a,b)X_{1}(\vartheta) + x_{2}(a,b)X_{2}(\vartheta) + x_{3}(a,b)T (\vartheta) \in \mathbb{H}^{1} | (a,b) \in \mathcal{O}\right\}.
\end{align*}
Similarly, let $\widehat{\vartheta}: \mathcal{O} \to \mathbb{H}^{1}$ be a $C^{2}$ parametrization of $\widehat{M}_{t}$, then
\begin{align*}
    \widehat{M}_{t} = \left\{\widehat{\vartheta}(a,b) =y_{1}(a,b)X_{1}(\widehat{\vartheta}) + y_{2}(a,b)X_{2}(\widehat{\vartheta}) + y_{3}(a,b)T (\widehat{\vartheta}) \in \mathbb{H}^{1} | (a,b) \in \mathcal{O}\right\}.
\end{align*}
From \cite[p.~361]{DGN2006}, we know that $V(x)= \vartheta_{a} \wedge \vartheta_{b}$ and
\begin{equation*}
    \begin{cases}
        & p =x_{2,a}x_{3,b} - x_{2,b}x_{3,a} - \frac{x_{2}}{2}\left(x_{1,a}x_{2,b}-x_{1,b}x_{2,a}\right) \\
        & q = x_{1,b}x_{3,a} - x_{1,a}x_{3,b} + \frac{x_{1}}{2}\left(x_{1,a}x_{2,b}-x_{1,b}x_{2,a}\right)\\
        & \omega = x_{1,a}x_{2,b}-x_{1,b}x_{2,a}
    \end{cases}
\end{equation*}
where $x_{i,a} = \frac{\partial x_{i}}{\partial a}$, $x_{i,b} = \frac{\partial x_{i}}{\partial b}$, $i=1,2,3$. Also, $\widehat{V}(y) = \widehat{\vartheta}_{a} \wedge \widehat{\vartheta}_{b}$ and
\begin{equation*}
    \begin{cases}
        & \widehat{p} =y_{2,a}y_{3,b} - y_{2,b}x_{3,a} - \frac{y_{2}}{2}\left(y_{1,a}y_{2,b}-y_{1,b}y_{2,a}\right) \\
        & \widehat{q} = y_{1,b}y_{3,a} - y_{1,a}y_{3,b} + \frac{y_{1}}{2}\left(y_{1,a}y_{2,b}-y_{1,b}y_{2,a}\right)\\
        & \widehat{\omega} = y_{1,a}y_{2,b}-y_{1,b}y_{2,a}
    \end{cases}
\end{equation*}
where $y_{i,a} = \frac{\partial y_{i}}{\partial a}$, $y_{i,b} = \frac{\partial y_{i}}{\partial b}$, $i=1,2,3$. Since $(y_{1}, y_{2}, y_{3}) = (\lambda x_{1}, \lambda x_{2}, \lambda^{2} x_{3})$, we obtain
\begin{equation*}
    \begin{cases}
        &y_{1,a} = \lambda x_{1,a} ; \qquad y_{1,b} = \lambda x_{1,b} \\
        &y_{2,a} = \lambda x_{2,a} ; \qquad y_{2,b} = \lambda x_{2,b}\\
        &y_{3,a} = \lambda^{2} x_{3,a} ; \qquad y_{3,b} = \lambda^{2} x_{3,b}
    \end{cases}
\end{equation*}
Thus
\begin{align}\label{8.4}
    \widehat{p} = \lambda^{3} p  ,\quad \widehat{q} = \lambda^{3} q, \quad \widehat{\omega} = \lambda^{2} \omega .
\end{align}
Using (\ref{2.3}), (\ref{8.2}) and (\ref{8.3}), we have
\begin{align*}
    \nu_{0}(x) = \frac{p}{\sqrt{p^{2}+q^{2}}}X_{1} + \frac{q}{\sqrt{p^{2}+q^{2}}}X_{2}, \quad  \widehat{\nu}_{0}(y) = \frac{\widehat{p}}{\sqrt{\widehat{p}^{2}+\widehat{q}^{2}}}X_{1} + \frac{\widehat{q}}{\sqrt{\widehat{p}^{2}+\widehat{q}^{2}}}X_{2}.
\end{align*}
Moreover
\begin{align}\label{8.5}
    \left\langle \nu_{0}(x), V(x)\right\rangle  = \sqrt{p^{2}+q^{2}}, \quad \left\langle \widehat{\nu}_{0}(y), \widehat{V}(y)\right\rangle =  \sqrt{\widehat{p}^{2}+\widehat{q}^{2}}.
\end{align}
Combining (\ref{8.4}) and (\ref{8.5}), we obtain (\ref{8.1}).
\hfill${\square}$

\begin{lemma}\label{lem8.2}
    Let $H_{0}(y)$ be the horizontal mean curvature at $y \in \widehat{M}_{t} \setminus \Sigma(\widehat{M}_{t})$. Then
    \begin{align}\label{8.6}
        H_{0} (y)  = \lambda^{-1}H_{0}(x)
    \end{align}
    where $H_{0}(x)$ is the horizontal mean curvature at $x = \delta_{\lambda^{-1}}(y) \in M_{t} \setminus \Sigma(M_{t})$ and $\lambda$ depends only on $t$.
\end{lemma}
\noindent {\bf Proof}: By the proof of Lemma \ref{lem8.1}, it follows that
\begin{align*}
    \nu_{0}(x) = \frac{p}{\sqrt{p^{2}+q^{2}}}X_{1}(x) + \frac{q}{\sqrt{p^{2}+q^{2}}}X_{2}(x), \quad
    \widehat{\nu}_{0}(y) = \frac{\widehat{p}}{\sqrt{\widehat{p}^{2}+\widehat{q}^{2}}}X_{1}(y) + \frac{\widehat{q}}{\sqrt{\widehat{p}^{2}+\widehat{q}^{2}}}X_{2}(y),
\end{align*}
and
\begin{align*}
    \widehat{p} = \lambda^{3} p  ,\quad \widehat{q} = \lambda^{3} q, \quad \widehat{\omega} = \lambda^{2} \omega .
\end{align*}
Thus
\begin{align*}
    \widehat{\nu}_{0}(y) = \frac{p}{\sqrt{p^{2}+q^{2}}}X_{1}(y) + \frac{q}{\sqrt{p^{2}+q^{2}}}X_{2}(y).
\end{align*}

From (\ref{2.4}), we have
\begin{align}\label{8.7}
    H_{0}(x) = \operatorname{div}_{0} (\nu_{0}(x)) = X_{1}\big|_{M_{t}}\left(\frac{p}{\sqrt{p^{2}+q^{2}}}\right) +  X_{2}\big|_{M_{t}}\left(\frac{q}{\sqrt{p^{2}+q^{2}}}\right),
\end{align}
and
\begin{align}\label{8.8}
    H_{0}(y) = \operatorname{div}_{0} (\widehat{\nu}_{0}(y)) = X_{1}\big|_{\widehat{M}_{t}}\left(\frac{p}{\sqrt{p^{2}+q^{2}}}\right) +  X_{2}\big|_{\widehat{M}_{t}}\left(\frac{q}{\sqrt{p^{2}+q^{2}}}\right).
\end{align}

Since $x = \delta_{\lambda^{-1}}(y)$, i.e., $(x_{1}, x_{2},x_{3}) = (\lambda^{-1}y_{1}, \lambda^{-1}y_{2}, \lambda^{-2}y_{3})$, we get
\begin{align*}
    \frac{\partial}{\partial y_{i}} = \frac{\partial x_{i}}{\partial y_{i}} \frac{\partial }{\partial x_{i}} = \lambda^{-1}\frac{\partial }{\partial x_{i}}, \quad i = 1,2; \qquad \frac{\partial}{\partial y_{3}} = \frac{\partial x_{3}}{\partial y_{3}} \frac{\partial }{\partial x_{3}} = \lambda^{-2}\frac{\partial }{\partial x_{3}}.
\end{align*}
From (\ref{2.1}), we have
\begin{align}
    X_{1}(y)\big|_{\widehat{M}_{t}} &= \frac{\partial }{\partial y_{1}} - \frac{y_{2}}{2}\frac{\partial }{\partial y_{3}} = \lambda^{-1}\frac{\partial }{\partial x_{1}} - \frac{\lambda x_{2}}{2}\left(\lambda^{-2}\frac{\partial }{\partial x_{3}}\right) = \lambda^{-1}X_{1}(x)\big|_{M_{t}},\label{8.9}\\
    X_{2}(y)\big|_{\widehat{M}_{t}} &= \frac{\partial }{\partial y_{2}} + \frac{y_{1}}{2}\frac{\partial }{\partial y_{3}} = \lambda^{-1}\frac{\partial }{\partial x_{2}} + \frac{\lambda x_{1}}{2}\left(\lambda^{-2}\frac{\partial }{\partial x_{3}}\right) = \lambda^{-1}X_{2}(x)\big|_{M_{t}} .\label{8.10}
\end{align}
Combining (\ref{8.7})-(\ref{8.10}), one can deduce that
\begin{align*}
    H_{0}(y)\big|_{\widehat{M}_{t} \setminus \Sigma(\widehat{M}_{t})} &=\lambda^{-1} X_{1}\big|_{M_{t}}\left(\frac{p}{\sqrt{p^{2}+q^{2}}}\right) +  \lambda^{-1}X_{2}\big|_{M_{t}}\left(\frac{q}{\sqrt{p^{2}+q^{2}}}\right) \\
    &= \lambda^{-1}H_{0}(x)\big|_{M_{t}\setminus \Sigma(M_{t})}.
\end{align*}
This completes the proof.
\hfill${\square}$

\begin{lemma}\label{lem8.3}
    Let $y=\delta_{\lambda(t)}(x)$ and $\lambda(t) = e^{-\frac{1}{3}t}$. Suppose that $x \in M_{t} \setminus \Sigma(M_{t})$ satisfy the equation (\ref{1.2}), then $y \in \widehat{M}_{t} \setminus \Sigma(\widehat{M}_{t})$ satisfies
    \begin{align}\label{8.11}
        \left\langle \frac{\partial y}{\partial t}, \widehat{\nu}\right\rangle = H_{0}^{-1}(y) \left\langle \widehat{\nu}_{0}, \widehat{\nu}\right\rangle  - \frac{1}{3}\left\langle y^{'}, \widehat{\nu} \right\rangle ,
    \end{align}
    where $y^{'} = (y_{1}, y_{2}, 2y_{3}) = y_{1} X_{1}(y) + y_{2} X_{2}(y) + 2y_{3} T$.
\end{lemma}
\noindent {\bf Proof}: Since $x \in M_{t} \setminus \Sigma(M_{t})$ satisfies (\ref{1.2}), we have
\begin{align}\label{8.12}
    \left\langle \frac{\partial x}{\partial t}, V(x)\right\rangle = H_{0}^{-1}\left\langle \nu_{0}(x), V(x)\right\rangle
\end{align}
where
\begin{align}\label{8.13}
    \frac{\partial x}{\partial t} = \frac{\partial x_{1}}{\partial t}X_{1} + \frac{\partial x_{2}}{\partial t}X_{2} +\left(\frac{\partial x_{3}}{\partial t}+\frac{x_{2}}{2}\frac{\partial x_{1}}{\partial t} -\frac{x_{1}}{2} \frac{\partial x_{2}}{\partial t}\right)T.
\end{align}
Combining (\ref{8.2}), (\ref{8.12}) and (\ref{8.13}) we get
\begin{align}\label{8.14}
    \left\langle \frac{\partial x}{\partial t}, V(x)\right\rangle =\frac{\partial x_{1}}{\partial t}p + \frac{\partial x_{2}}{\partial t}q +\left(\frac{\partial x_{3}}{\partial t}+\frac{x_{2}}{2}\frac{\partial x_{1}}{\partial t} -\frac{x_{1}}{2} \frac{\partial x_{2}}{\partial t}\right)\omega = H_{0}^{-1}\left\langle \nu_{0}(x), V(x)\right\rangle.
\end{align}
Moreover
\begin{align}\label{8.15}
    \frac{\partial y}{\partial t} = & \frac{\partial y_{1}}{\partial t}X_{1} + \frac{\partial y_{2}}{\partial t}X_{2} +\left(\frac{\partial y_{3}}{\partial t}+\frac{y_{2}}{2}\frac{\partial y_{1}}{\partial t} -\frac{y_{1}}{2} \frac{\partial y_{2}}{\partial t}\right)T \notag \\
     =&  \left(e^{-\frac{1}{3}t}\frac{\partial x_{1}}{\partial t} - \frac{y_{1}}{3}\right)X_{1} + \left(e^{-\frac{1}{3}t}\frac{\partial x_{2}}{\partial t} - \frac{y_{2}}{3}\right)X_{2} \\
    & + \left(e^{-\frac{2}{3}t}\frac{\partial x_{3}}{\partial t} - \frac{2y_{3}}{3} + e^{-\frac{2}{3}t} \frac{x_{2}}{2}\frac{\partial x_{1}}{\partial t} - e^{-\frac{2}{3}t} \frac{x_{1}}{2}\frac{\partial x_{2}}{\partial t}\right)T . \notag
\end{align}
Using (\ref{8.3}), (\ref{8.4}), (\ref{8.14}) and (\ref{8.15}), we obtain
\begin{align}\label{8.16}
    \left\langle \frac{\partial y}{\partial t}, \widehat{V}(y)\right\rangle =& \left[e^{-\frac{1}{3}t}\widehat{p}\frac{\partial x_{1}}{\partial t}+ e^{-\frac{1}{3}t}\widehat{q}\frac{\partial x_{2}}{\partial t} + e^{-\frac{2}{3}t}\widehat{\omega}\left(\frac{\partial x_{3}}{\partial t} +\frac{x_{2}}{2}\frac{\partial x_{1}}{\partial t} -\frac{x_{1}}{2} \frac{\partial x_{2}}{\partial t} \right)\right] \notag \\
    &- \frac{1}{3}\left\langle y^{'}, \widehat{V}(y)\right\rangle. \notag \\
    =& e^{-\frac{4}{3}t}\left[p\frac{\partial x_{1}}{\partial t}+ q\frac{\partial x_{2}}{\partial t} + \omega\left(\frac{\partial x_{3}}{\partial t} +\frac{x_{2}}{2}\frac{\partial x_{1}}{\partial t} -\frac{x_{1}}{2} \frac{\partial x_{2}}{\partial t} \right)\right] - \frac{1}{3}\left\langle y^{'}, \widehat{V}(y)\right\rangle. \notag \\
    = & e^{-\frac{4}{3} t } H_{0}^{-1}(x)\left\langle \nu_{0}(x), V(x)\right\rangle - \frac{1}{3}\left\langle y^{'}, \widehat{V}(y)\right\rangle
\end{align}
From (\ref{8.1}) and (\ref{8.6}), we have
\begin{align} \label{8.17}
    H_{0}^{-1}(y)\left\langle \widehat{\nu}_{0}(y), \widehat{V}(y)\right\rangle = e^{-\frac{4}{3} t } H_{0}^{-1}(x)\left\langle \nu_{0}(x), V(x)\right\rangle.
\end{align}
Combining (\ref{8.16}) with (\ref{8.17}), we obtain
\begin{align*}
    \left\langle \frac{\partial y}{\partial t}, \widehat{V}(y)\right\rangle = H_{0}^{-1}(y)\left\langle \widehat{\nu}_{0}(y), \widehat{V}(y)\right\rangle - \frac{1}{3}\left\langle y^{'}, \widehat{V}(y) \right\rangle.
\end{align*}
This implies that  (\ref{8.11}) holds.
\hfill${\square}$

Suppose there exists a function $h: \mathbb{H}^{1} \times [0, + \infty) \to \mathbb{R}$ such that
\begin{align*}
    \widehat{E}_{t} = \delta_{\lambda(t)}(E_{t}) =\left\{y \in \mathbb{H}^{1}, y=\delta_{\lambda(t)}(x), u(x)< t \right\}= \{y \in \mathbb{H}^{1}, h(y,t)< 0 \}.
\end{align*}
Then, using the level set representation, we obtain
\begin{align*}
    0  = \left\langle \frac{\partial h}{\partial y},\frac{\partial y}{\partial t}\right\rangle _{\mathbb{R}^{3}}+ \frac{\partial h}{\partial t}= \left\langle \nabla h, \frac{\partial y}{\partial t} \right\rangle + \frac{\partial h}{\partial t}.
\end{align*}
Thus
\begin{align*}
    \frac{\partial h}{\partial t}  = -\left\langle \widehat{V}, \frac{\partial y}{\partial t}\right\rangle  &= -\left(H_{0}^{-1}(y) \left\langle \widehat{V}, \widehat{\nu}_{0}\right\rangle  - \frac{1}{3}\left\langle y^{'}, \widehat{V} \right\rangle\right) \notag \\
    &= - \left(H_{0}^{-1}(y) \left\langle \nabla h, \frac{\nabla_{0} h}{| \nabla_{0} h|}\right\rangle  - \frac{1}{3}\left\langle y^{'}, \nabla h  \right\rangle\right) \notag \\
    & = - \left(\left[\operatorname{div}_{0} \left(\frac{\nabla_{0} h}{| \nabla_{0} h|}\right)\right]^{-1}| \nabla_{0} h|  - \frac{1}{3}\left\langle y^{'}, \nabla h  \right\rangle\right) .
\end{align*}

\begin{lemma}\label{lem8.4}
    The $\mathbb{H}$-perimeter of $\widehat{E}_{t}$ remains constant for all $t> 0$ along the flow
    \begin{align*}
        \left\langle \frac{\partial y}{\partial t}, \widehat{\nu}\right\rangle = H_{0}^{-1}(y) \left\langle \widehat{\nu}_{0}, \widehat{\nu}\right\rangle  - \frac{1}{3}\left\langle y^{'}, \widehat{\nu} \right\rangle.
    \end{align*}
    Moreover, this invariance extends to $t = 0$ provided that $E_{0}$ is a horizontal minimizing hull.
\end{lemma}
\noindent {\bf Proof}: Since $\widehat{E}_{t} = \delta_{\lambda(t)}(E_{t})$ and $\lambda(t) = e^{-\frac{1}{3}t}$, it follows from Lemma \ref{lem2.17} that
\begin{align}\label{8.18}
    P_{\mathbb{H}}(\widehat{E}_{t}) = e^{-t} P_{\mathbb{H}}(E_{t}).
\end{align}
Similarly, let $\widehat{E}^{+}_{t} = \delta_{e^{-\frac{1}{3}t}}({E}^{+}_{t}) $. Then,
\begin{align}\label{8.19}
    P_{\mathbb{H}}(\widehat{E}^{+}_{t}) = e^{-t} P_{\mathbb{H}}(E^{+}_{t}).
\end{align}
Combining (\ref{7.3}) and (\ref{8.18}), we obtain
\begin{align*}
    P_{\mathbb{H}}(\widehat{E}_{t}) = e^{-t} P_{\mathbb{H}}(E_{t}) = P_{\mathbb{H}}(E_{0}), \qquad t\in (0, \infty).
\end{align*}
Also, if $E_{0}$ is a horizontal minimizing hull, combining (\ref{7.4}) and (\ref{8.19}), we get
\begin{align*}
    P_{\mathbb{H}}(\widehat{E}_{t}) = P_{\mathbb{H}}(\widehat{E}^{+}_{t}) = e^{-t} P_{\mathbb{H}}(E^{+}_{t}) = P_{\mathbb{H}}(E^{+}_{0}) =P_{\mathbb{H}}(E_{0}), \qquad  t\in [0 ,\infty).
\end{align*}
Therefore, the rescaled flow (\ref{1.8}) is an $\mathbb{H}$-perimeter preserving flow.
\hfill${\square}$

\subsection{ Minkowski type formula }
\
   \vglue-10pt
    \indent

\noindent {\bf Proof of Theorem \ref{thm1.7}}: Combining the first variation formula (\ref{2.15}) for the $\mathbb{H}$-perimeter and the horizontal variation (\ref{1.8}), i.e.,
\begin{align*}
    \left\langle \frac{\partial y}{\partial t}, \widehat{\nu}\right\rangle = \left(H_{0}^{-1}(y)   - \frac{\left\langle y^{'}, \widehat{\nu} \right\rangle}{3\left\langle \widehat{\nu}_{0}, \widehat{\nu}\right\rangle} \right) \left\langle \widehat{\nu}_{0}, \widehat{\nu}\right\rangle,
\end{align*}
we have
\begin{align}\label{8.20}
    \frac{\partial }{\partial t} P_{\mathbb{H}}(\widehat{E}_{t}) &= \int_{\widehat{M}_{t} \setminus \Sigma(\widehat{M}_{t})} \left(H_{0}^{-1}(y)   - \frac{\left\langle y^{'}, \widehat{\nu} \right\rangle}{3\left\langle \widehat{\nu}_{0}, \widehat{\nu}\right\rangle} \right)H_{0}(y) d\sigma_{\mathbb{H}} \notag \\
    &= \int_{\widehat{M}_{t} \setminus \Sigma(\widehat{M}_{t})} d\sigma_{\mathbb{H}} - \frac{1}{3} \int_{\widehat{M}_{t} \setminus \Sigma(\widehat{M}_{t})}  H_{0} \left\langle y^{'}, \widehat{\nu} \right\rangle d\mathcal{H}^{2} \notag \\
    & = P_{\mathbb{H}}(\widehat{E}_{t}) - \frac{1}{3} \int_{\widehat{M}_{t} \setminus \Sigma(\widehat{M}_{t})}  H_{0} \left\langle y^{'}, \widehat{\nu} \right\rangle d\mathcal{H}^{2}.
\end{align}
From Lemma \ref{lem8.4}, $P_{\mathbb{H}}(\widehat{E}_{t})$ remains constant for all $t > 0$, thus
\begin{align}\label{8.21}
    P_{\mathbb{H}}(\widehat{E}_{t})= P_{\mathbb{H}}(E_{0}),\qquad \frac{\partial }{\partial t} P_{\mathbb{H}}(\widehat{E}_{t}) = 0.
\end{align}
Combining (\ref{8.20}) with (\ref{8.21}), we obtain
\begin{align*}
   3P_{\mathbb{H}}(E_{0}) = 3P_{\mathbb{H}}(\widehat{E}_{t}) = \int_{\widehat{M}_{t} \setminus \Sigma(\widehat{M}_{t})}  H_{0} \left\langle y^{'}, \widehat{\nu} \right\rangle d\mathcal{H}^{2}, \qquad \text{ for any } t \in (0, +\infty).
\end{align*}
The above equality yields
\begin{align}\label{8.22}
    3P_{\mathbb{H}}(E_{0}) =  \lim_{t \to 0^{+}} \int_{\widehat{M}_{t} \setminus \Sigma(\widehat{M}_{t})}  H_{0} \left\langle y^{'}, \widehat{\nu} \right\rangle d\mathcal{H}^{2}.
\end{align}

Next, we estimate the right-hand side of (\ref{8.22}). Fix $t \in (0, \varepsilon_{0} ]$ with $\varepsilon_{0} >0$ sufficiently small. Let $Q(y(t)) := \chi_{\widehat{M}_{t} \setminus \Sigma(\widehat{M}_{t})}  H_{0}(y) \left\langle y^{'}, \widehat{\nu}(y) \right\rangle $. From (\ref{2.6}), (\ref{2.7}) and (\ref{8.6}), we have
\begin{align*}
    H_{0}(y) = e^{\frac{1}{3}t}H_{0}(x) = e^{\frac{1}{3}t}|\nabla_{0}  u (x)|  \leq Ce^{\frac{1}{3}\varepsilon_{0}}, \qquad x \in M_{t}\setminus \Sigma(M_{t}) .
\end{align*}
where the last inequality follows from Proposition \ref{prop5.17} and $C >0$ is a constant. Since $y = (y_{1},y_{2},y_{3}) = \delta_{e^{-\frac{1}{3}t}}(x) = (e^{-\frac{1}{3}t}x_{1}, e^{-\frac{1}{3}t}x_{2}, e^{-\frac{2}{3}t}x_{3})$, we have
\begin{align*}
    \left\langle y^{'}, \widehat{\nu}(y) \right\rangle \leq \left\langle y^{'}, y^{'} \right\rangle ^{\frac{1}{2}} = \sqrt{y_{1}^{2} + y_{2}^{2} + 4y_{3}^{2}} \leq e^{\frac{1}{3}\varepsilon_{0}} \sqrt{x_{1}^{2} + x_{2}^{2} + 4e^{\frac{2}{3}\varepsilon_{0}}x_{3}^{2}}.
\end{align*}
Combining the above inequalities, we get
\begin{align*}
    |Q(y(t))| \leq  Ce^{\frac{2}{3}\varepsilon_{0}}   \sqrt{x_{1}^{2} + x_{2}^{2} + 4e^{\frac{2}{3}\varepsilon_{0}}x_{3}^{2}}.
\end{align*}
for any $ x \in M_{t}\setminus \Sigma(M_{t})$ and $t \in (0, \varepsilon_{0} ]$. Moreover, since $E_{0}$ is a bounded open set, there exist a Kor$\acute{ \rm a}$nyi ball $\mathsf{B}_{R_{0}}(o)$ with finite radius $R_{0} >0$ such that $E_{0} \subseteq \mathsf{B}_{R_{0}}(o)$. Let $R(t) = e^{\frac{1}{3}t} R_{0}$. Then, by Proposition \ref{prop3.6}, it follows that
\begin{align*}
    \mathsf{B}_{R(t)}(o) = \left\{x \in \mathbb{H}^{1}, \rho (x,o) = \left(\left(x_{1}^{2}+x_{2}^{2}\right)^{2}+16 x_{3}^{2}\right)^{1 / 4} < R(t)\right\}
\end{align*}
is a solution of the HIMCF with initial domain $\mathsf{B}_{R_{0}}(o)$. Since $E_{0} \subseteq \mathsf{B}_{R_{0}}(o)$, it follows from Theorem \ref{thm6.1} that $E_{t} \subseteq \mathsf{B}_{R(t)}(o)$. Thus, for all $t \in (0, \varepsilon_{0} ]$ and $x \in M_{t}\setminus \Sigma(M_{t})$,
\begin{align*}
    \left(\left(x_{1}^{2}+x_{2}^{2}\right)^{2}+16 x_{3}^{2}\right)^{1 / 4} \leq  e^{\frac{1}{3}t} R_{0} \leq e^{\frac{1}{3}\varepsilon_{0}} R_{0},
\end{align*}
which implies that $\sqrt{x_{1}^{2} + x_{2}^{2} + 4e^{\frac{2}{3}\varepsilon_{0}}x_{3}^{2}}$ is bounded. Consequently, $|Q(y(t))|$ is bounded for any $t \in (0, \varepsilon_{0}]$. By Lebesgue's dominated convergence theorem, we obtain
\begin{align}\label{8.23}
    \lim_{t \to 0^{+}} \int_{\mathbb{H}^{1}} Q(y(t)) d\mathcal{H}^{2} = \int_{\mathbb{H}^{1}} \lim_{t \to 0^{+}} Q(y(t)) d\mathcal{H}^{2} = \int_{\widehat{M}_{0} \setminus \Sigma(\widehat{M}_{0})}  H_{0} \left\langle y^{'}, \widehat{\nu} \right\rangle d\mathcal{H}^{2}.
\end{align}
Combining (\ref{8.22}) and (\ref{8.23}), we have
\begin{align*}
    3P_{\mathbb{H}}(E_{0}) = \int_{M_{0} \setminus \Sigma(M_{0})}  H_{0} \left\langle y^{'}, \widehat{\nu} \right\rangle d\mathcal{H}^{2}.
\end{align*}
This completes the proof.
\hfill${\square}$

\noindent {\bf Proof of Corollary \ref{cor1.8}}: For each $x\in M\setminus \Sigma(M)$, $H_{0}$ is a  positive constant. From (\ref{1.9}), we have
\begin{align}\label{8.24}
     3 P_{\mathbb{H}}(E) =\int_{ M\setminus \Sigma (M)} H_{0} \left\langle y^{'}, \widehat{\nu} \right\rangle  d\mathcal{H}^{2} = H_{0}  \int_{ M\setminus \Sigma (M)} \left\langle y^{'}, \widehat{\nu} \right\rangle d\mathcal{H}^{2} = H_{0}\int_{M } \left\langle y^{'}, \widehat{\nu} \right\rangle d\mathcal{H}^{2}.
\end{align}
where the last equality follows from the fact that $\mathcal{H}^{2}(\Sigma(M)) = 0$. Using the Riemannian divergence theorem, we have
\begin{align}\label{8.25}
    \int_{M} \left\langle y^{'} , \widehat{\nu}\right\rangle d\mathcal{H}^{2} = \int_{E} \operatorname{div} y^{'} dy_{1} \wedge  dy_{2} \wedge  dy_{3} = 4 |E|.
\end{align}
Combining (\ref{8.24}) and  (\ref{8.25}), the proof of Corollary \ref{cor1.8} is complete.
\hfill${\square}$

\section{ Appendix}
\noindent

\subsection{ Proof of Theorem \ref{thmS}}
\
   \vglue-10pt
    \indent

In order to prove Theorem \ref{thmS}, we shall need the following $L^{p}$-continuity of the Hardy-Littlewood maximal operator. Recall that the Hardy–Littlewood maximal operator with respect to $CC$-balls is defined by
\begin{align*}
    Mf(x) = \sup_{r>0}\frac{1}{|B(x,r)|} \int_{B(x,r)} |f(y)| dy.
\end{align*}
\begin{proposition}\label{prop9.1}
    Let $1 < p < \infty$, $f \in L^{p}(\mathbb{H}^{1})$ and $\lambda >0$. Then,
    \begin{align*}
        \int_{\mathbb{H}^{1}} |Mf(x)|^{p} dx \leq  C(c,Q)  \frac{p2^{p-1}}{p-1}\left\|f\right\|^{p}_{L^{p}(\mathbb{H}^{1})}
    \end{align*}
    where $C(c,Q)$ is a constant depends on doubling constant $c$ and homogeneous dimension $Q$.
\end{proposition}

\noindent {\bf Proof }: For each $R > 0$, we define
\begin{align*}
    M_{R}f(x) = \sup_{0<r\leq R}\frac{1}{|B(x,r)|} \int_{B(x,r)} |f(y)| dy, \quad  E_{\lambda}^{R} = \{x \in \mathbb{H}^{1}, M_{R}f(x) > \lambda\}.
\end{align*}
Then, $Mf(x)  = \lim_{R \to \infty } M_{R}f(x)$. Let $r \leq R$ and
\begin{align}\label{9.1}
    \lambda |B_{r}| < \int_{B_{r}} |f(y)| dy .
\end{align}
Let $\mathcal{B} $ be a family of $CC$-balls satisfies (\ref{9.1}). Then, the union of all $CC$-balls in $\mathcal{B} $ contains $E_{\lambda}^{R}$. It follows from the covering lemma (Lemma 1.1 in \cite{Bu1999}) that there exists a countable pairwise disjoint family $\{B(x_{i}, r_{i})\}_{i \in \mathbb{N}} \in  \mathcal{B}  $ such that $E_{\lambda}^{R} \subset \bigcup_{i \in \mathbb{N}} B(x_{i}, 5r_{i})  $. Moreover, we know that
\begin{align*}
    \lambda |B(x_{i},r_{i})| < \int_{B(x_{i},r_{i})} |f(y)|dy.
\end{align*}
Combining with Lemma \ref{lemD}, we get
\begin{align*}
    |E_{\lambda}^{R}| \leq \sum_{i}  |B(x_{i}, 5r_{i})| \leq c^{-1} 5^{Q} \sum_{i}  |B(x_{i}, r_{i})| \leq \frac{C(c,Q)}{\lambda} \int_{\bigcup_{i} B(x_{i}, r_{i})} |f(y)|dy
\end{align*}
where $C(c,Q)$ independent of $R$. Thus
\begin{align*}
   \lim_{R \to \infty}  |E_{\lambda}^{R}| \leq  \frac{C(c,Q)}{\lambda} \int_{\bigcup_{i} B(x_{i}, r_{i})} |f(y)|dy,
\end{align*}
and
\begin{align}\label{9.2}
    |\{x\in \mathbb{H}^{1} \vert  Mf(x) > \lambda\}| \leq \frac{C(c,Q)}{\lambda} \int_{\mathbb{H}^{1}} |f(y)|dy.
\end{align}
Let $f_{1} = f\cdot \chi_{\{|f| \leq \frac{\lambda}{2}\}}$ and $f_{2} = f\cdot \chi_{\{|f| >  \frac{\lambda}{2}\}} $, then $f = f_{1} + f_{2}$ and
\begin{align*}
    Mf(x) \leq Mf_{1}(x) + Mf_{2}(x), \quad  \{ Mf(x) > \lambda\} \subset \left(\left\{Mf_{1}(x) > \lambda/2\right\} \cup \left\{Mf_{2}(x) > \lambda/2\right\}\right)
\end{align*}
Since
\begin{align*}
    Mf_{1}(x) = \sup_{r>0}\fint_{B(x,r)} |f_{1}(y)|dy  = \sup_{r>0} \fint_{B(x,r)\cap \{ |f| \leq \frac{\lambda}{2}\}} |f(y)|dy \leq \lambda/2,
\end{align*}
we have $\{Mf_{1}(x) > \lambda /2\} = \emptyset  $ and $\{ Mf(x) > \lambda\} \subset \left\{Mf_{2}(x) > \lambda /2\right\}$. Using (\ref{9.2}), we obtain
\begin{align*}
    \int_{\mathbb{H}^{1}} |Mf(x)|^{p} dx &= p\int_{0}^{\infty} \lambda^{p-1} |\{Mf(x) > \lambda\}| d\lambda
    \leq p\int_{0}^{\infty} \lambda^{p-1} \left|\left\{Mf_{2}(x) > \lambda/2\right\}\right| d\lambda \notag \\
    &\leq p\int_{0}^{\infty} \lambda^{p-1} \left(\frac{C(c,Q)}{\lambda} \int_{\{|f| > \frac{\lambda}{2}\}} |f(y)|dy\right) d\lambda \leq \frac{p 2^{p-1}}{p-1} C(c,Q) \|f\|^{p}_{L^{p}(\mathbb{H}^{1})}
\end{align*}
This completes the proof of Proposition \ref{prop9.1}.
\hfill${\square}$
\begin{remark}\label{rem9.2}
    The above proof follows the standard argument for the $L^{p}$-continuity of the Hardy-Littlewood maximal operator on spaces of homogeneous type (see, e.g., \cite{C1976}). The details are included here in order to express the constant appearing in this result explicitly in terms of $p$.
\end{remark}

\begin{proposition}\cite{CGL1993}\label{prop9.3}
    Let $\Omega \subset \mathbb{H}^{1}$ be a $C^1$ domain and let $u \in C_0^1(\Omega)$. Then, for every $x \in \Omega$, we have
    \begin{align}\label{9.3}
        u(x)=\int_{\Omega} \left\langle \nabla_{0} \Gamma(x,y), \nabla_{0}u(y)\right\rangle  d y
    \end{align}
    where $\Gamma(x, y)$ denotes the (positive) fundamental solution of horizontal Laplacian
\end{proposition}

By \cite{NSW1985,SC1984}, there exist positive constants $C$ and $R_{0}$ such that for all $x \in U \subset \subset \mathbb{H}^{1}$ and $0 < d(x, y) < R_{0}$,
\begin{align}\label{9.4}
    \left|\nabla_{0} \Gamma(x, y)\right| & \leq C^{-1} \frac{d(x, y)}{|B(x, d(x, y))|}, \quad x \neq y
\end{align}
where $C$, $R_{0}$ depend on $\mathbb{H}^{1}$. For $u \in C_0^1\left(B_R\right)$ and $x \in B_R$, by (\ref{9.3}) and (\ref{9.4}), we have
\begin{align}\label{9.5}
    |u(x)| \leq \int_{B_R}\left|\nabla_{0}\Gamma(x, y)\right|\left|\nabla_{0} u(y)\right| d y \leq C \int_{B_R}\left|\nabla_{0} u(y)\right| \frac{d(x, y)}{|B(x, d(x, y))|} d y
\end{align}
Let $0<\beta \leq Q$. The fractional integral operator $I_\beta$ is defined by
\begin{align}\label{9.6}
    I_\beta(f)(x)=\int_{B_R}|f(y)| \frac{d(x, y)^\beta}{|B(x, d(x, y))|} d y .
\end{align}
Using (\ref{9.6}), we can rewrite (\ref{9.5}) as follows
\begin{align}\label{9.7}
    |u(x)| \leq C I_1\left(\left|\nabla_{0} u\right|\right)(x),\quad  \text { for } x \in B_R .
\end{align}
where $C$ depend on $\mathbb{H}^{1}$.

\begin{theorem}\label{thm9.4}
    Let $U \subset \mathbb{H}^{1}$ be a bounded open set. Assume $1 \leq p \leq \infty$ and $0<\beta \leq Q$. Then, $I_\beta$ maps $L^p\left(B_R\right)$ continuously into $L^q\left(B_R\right)$, with $0 \leq \frac{1}{p}-\frac{1}{q} \leq \frac{\beta}{Q}$. Moreover, there exist constants $\overline{C} >0$ and $R_0>0$ such that for any $x_0 \in U, B_R=B\left(x_0, R\right)$, with $R \leq R_0$, we have
    \begin{align*}
        \left( \fint_{B_R}\left|I_\beta(f)(x)\right|^q d x\right)^{\frac{1}{q}} \leq \overline{C} R^\beta\left(\fint_{B_R}|f(x)|^p d x\right)^{\frac{1}{p}},
    \end{align*}
    for every $f \in L^p\left(B_R\right)$, where the constant $\overline{C}$ is defined as (\ref{9.25}).  When $p=1$ one must have $0 \leq \frac{1}{p}-\frac{1}{q}<\frac{\beta}{Q}$.
\end{theorem}
\noindent {\bf Proof }: The proof is divided into two cases. We begin with the case $0 \leq \frac{1}{p} - \frac{1}{q} < \frac{\beta}{Q}$. Let $r \geq 1$ be defined by $1 -\frac{1}{r} = \frac{1}{p} - \frac{1}{q}$. For a fixed $x \in B_{R}$, let $h(x, y): y \longmapsto  \frac{d(x, y)^\beta}{|B(x, d(x, y))|}$, and $f \in L^p\left(B_R\right)$. Using H$\ddot{\rm o} $lder's inequality, we have
\begin{align*}
    \left|I_\beta(f)(x)\right|  = &\int_{B_{R}} |f(y)| h(x,y) dy = \int_{B_{R}} h(x, y)^{\frac{r}{q}} h(x, y)^{1-\frac{r}{q}}|f(y)|^{\frac{p}{q}}|f(y)|^{1-\frac{p}{q}} dy \\
    \leq &  \left(\int_{B_{R}} |f(y)|^{p}h(x,y)^{r} dy\right)^{\frac{1}{q}} \left(\int_{B_{R}} h(x,y)^{r} dy\right)^{\frac{q-r}{qr}} \left(\int_{B_{R}} |f(y)|^{p} dy\right)^{\frac{q-p}{qp}}.
\end{align*}
Thus,
\begin{align}\label{9.8}
    \int_{B_{R}} \left|I_\beta(f)(x)\right|^{q}  dx \leq  \sup_{x\in B_{R}} \|h(x,y)\|^{q}_{L^{r}(B_{R})} \left(\int_{B_{R}} |f(y)|^{p} dy\right)^{\frac{q}{p}}.
\end{align}
It remains to estimate $\sup _{x \in B_R}\|h(x, \cdot)\|_{L^r\left(B_R\right)}$. Since $x \in B_R$, we have $B_R \subset B(x, 2 R)$ and $B(x, R) \subset B_{2R}$. Then,
\begin{align}\label{9.9}
    \int_{B_R} h(x, y)^r d y & \leq \int_{B(x, 2 R)} \frac{d(x, y)^{\beta r}}{|B(x, d(x, y))|^r} d y =\sum_{k=0}^{\infty} \int_{2^{-k} R<d(x, y)<2^{-k+1} R} \frac{d(x, y)^{\beta r}}{|B(x, d(x, y))|^r} d y  \notag \\
    &\leq \sum_{k=0}^{\infty} \frac{(2^{-k+1}R)^{\beta r}}{|B(x, 2^{-k}R)|^{r}} |B(x, 2^{-k+1}R)|
\end{align}
By (\ref{2.d}), we have $|B(x, 2^{-k}R)| \geq c 2^{-Q}|B(x, 2^{-k+1}R)|$, $|B(x, 2^{-k}R)| \geq c 2^{-kQ}|B(x,R)|$ and
\begin{align}\label{9.10}
    |B(x, R)| \geq C(c,Q)|B_{2R}| \geq C(c,Q)|B_{R}|.
\end{align}
Thus, (\ref{9.9}) implies that
\begin{align*}
    \int_{B_R} h(x, y)^r d y &\leq C(c,Q)  \frac{2^{\beta r} R^{\beta r}}{|B_{R}|^{r-1}}\sum_{k=0}^{\infty} \left(2^{\beta r - Q(r-1)} \right)^{-k}\\
    &\leq C(c,Q) \frac{2^{\beta r -Q(r-1) }}{2^{\beta r - Q(r-1)} -1}\frac{2^{\beta r }R^{\beta r}}{|B_{R}|^{r-1}}\leq \frac{C(c,Q)4^{\beta r }}{2^{\beta r - Q(r-1)} -1}\frac{R^{\beta r}}{|B_{R}|^{r-1}}
\end{align*}
where we used $\sum_{k=0}^{\infty} \left(2^{\beta r - Q(r-1)} \right)^{-k} = \frac{2^{\beta r - Q(r-1)}}{2^{\beta r - Q(r-1)} -1}$. Hence,
\begin{align}\label{9.11}
    \sup _{x \in B_R}\|h(x, \cdot)\|_{L^r\left(B_R\right)}^r \leq \frac{C(c,Q) 4^{\beta r}}{2^{\beta r -Q(r-1)} -1} \frac{R^{\beta r}}{\left|B_R\right|^{r-1}}.
\end{align}
Combining (\ref{9.8}) and (\ref{9.11}), we obtain
\begin{align*}
    \left(\int_{B_{R}} \left|I_\beta(f)(x)\right|^{q}  dx \right)^{\frac{1}{q}} \leq \frac{C(c,Q)4^{\beta }}{(2^{\beta r - Q(r-1)} -1)^{\frac{1}{r}}} \frac{R^{\beta}}{|B_{R}|^{1-\frac{1}{r}}}\left(\int_{B_{R}} |f(y)|^{p} dy\right)^{\frac{1}{p}}.
\end{align*}
Since $|B_{R}|^{1-\frac{1}{r}} = |B_{R}|^{\frac{1}{p}-\frac{1}{q}}$ and $4^{\beta} \leq 4^{Q}$, we have
\begin{align}\label{9.12}
     \left(\fint_{B_{R}} \left|I_\beta(f)(x)\right|^{q}  dx \right)^{\frac{1}{q}} \leq \frac{C(c,Q)R^{\beta}}{\left(2^{\frac{\beta pq-Q(q-p)}{pq-q+p}} -1\right)^{\frac{pq-q+p}{pq}}}\left(\fint_{B_{R}} |f(y)|^{p} dy\right)^{\frac{1}{p}}.
\end{align}

Next, consider the case $\frac{1}{p}-\frac{1}{q} = \frac{\beta}{Q}$, which implies $1< p < \frac{Q}{\beta}$. For any $x \in B_{R}$ and $0 < \varepsilon < R$, let
\begin{align*}
     I^{1}_\beta(f)(x) = \int_{B(x,\varepsilon)}  |f(y)| \frac{d(x, y)^\beta}{|B(x, d(x, y))|} d y,\quad
     I^{2}_\beta(f)(x) = \int_{B^{c}(x,\varepsilon)\cap B_{R}}  |f(y)| \frac{d(x, y)^\beta}{|B(x, d(x, y))|} d y .
\end{align*}
Then, $I_\beta(f)(x) = I^{1}_\beta(f)(x) + I^{2}_\beta(f)(x)$. We first estimate $I^{1}_\beta(f)(x)$.
\begin{align}\label{9.13}
    I^{1}_\beta(f)(x)&=\sum_{k=0}^{\infty} \int_{2^{-(k+1)} \varepsilon <d(x, y)<2^{-k} \varepsilon} |f(y)| \frac{d(x, y)^{\beta }}{|B(x, d(x, y))|} d y \notag \\
    &\leq  \sum_{k=0}^{\infty}\frac{(2^{-k} \varepsilon)^{\beta}}{|B(x,2^{-(k+1)} \varepsilon)|} \int_{B(x,2^{-k} \varepsilon)} |f(y)| dy \notag \\
    &\leq C(c,Q)\varepsilon^{\beta} \sum_{k=0}^{\infty} (2^{-\beta} )^{k} \left(\frac{1}{|B(x,2^{-k} \varepsilon)|}\int_{B(x,2^{-k} \varepsilon)} |f(y)| dy\right) \notag \\
    &\leq  C(c,Q) \frac{2^{\beta}}{2^{\beta} -1}\varepsilon^{\beta}Mf(x)
\end{align}
where we used $|B(x,2^{-(k+1)} \varepsilon)| \geq c 2^{-Q} |B(x,2^{-k} \varepsilon)| $ and $\sum_{k=0}^{\infty} (2^{-\beta})^{k}  = \frac{2^{\beta}}{2^{\beta} -1}$.

It remains to establish an upper bound for $I^{2}_\beta(f)(x)$. Set $f \equiv 0$ on $B_{R}^{c}$. Observe that $B_{R} \subset B(x,2R)$ for every $x \in B_{R}$. Applying H$\ddot{\rm o} $lder's inequality with exponent $p$, we obtain
\begin{align}\label{9.14}
    I^{2}_\beta(f)(x) \leq \left(\int_{B_{R}} |f(y)|^{p} dy\right)^{\frac{1}{p}} \left(\int_{B^{c}(x,\varepsilon)\cap B(x,2R)} \left(\frac{d(x, y)^\beta}{|B(x, d(x, y))|}\right)^{\frac{p}{p-1}} d y\right)^{\frac{p-1}{p}}.
\end{align}
Choose $k_{0} \in \mathbb{N}$ such that $2^{k_{0}}\varepsilon \leq 2R < 2^{k_{0}+1}\varepsilon$. Then, $ B^{c}(x,\varepsilon) \cap B(x,2R) \subseteq \bigcup_{k=0}^{k_{0}} \{ 2^{k} \varepsilon< d(x,y) < 2^{k+1} \varepsilon\} $ and
\begin{align}\label{9.15}
    \int_{B^{c}(x,\varepsilon)\cap B(x,2R)} \left(\frac{d(x, y)^\beta}{|B(x, d(x, y))|}\right)^{\frac{p}{p-1}} d y &\leq  \sum_{k=0}^{k_{0}} \int_{\{ 2^{k} \varepsilon< d(x,y) < 2^{k+1} \varepsilon\} } \left(\frac{d(x, y)^\beta}{|B(x, d(x, y))|}\right)^{\frac{p}{p-1}} d y \notag \\
    &\leq \sum_{k=0}^{k_{0}}\left(\frac{(2^{k+1} \varepsilon)^{\beta}}{ |B(x,2^{k} \varepsilon)|}\right)^{\frac{p}{p-1}}|B(x,2^{k+1} \varepsilon )|.
\end{align}
Obviously, $2^{k+1}\varepsilon \leq 2^{k-k_{0}+2}R $ and $ 2^{k-k_{0}}R  < 2^{k}\varepsilon $. By (\ref{2.d}), we have
\begin{align*}
    |B(x,2^{k} \varepsilon)| \geq c2^{-Q}|B(x,2^{k+1} \varepsilon )|, \quad
    |B(x,2^{k}\varepsilon) |\geq |B(x,2^{k-k_{0}}R) |\geq c2^{(k-k_{0})Q}|B(x,R)|
\end{align*}
Hence, (\ref{9.15}) yields that
\begin{align}\label{9.16}
    \int_{B^{c}(x,\varepsilon)\cap B(x,2R)}& \left(\frac{d(x, y)^\beta}{|B(x, d(x, y))|}\right)^{\frac{p}{p-1}} d y \notag \\
    & \leq C(c,Q) \sum_{k=0}^{k_{0}}(2^{k-k_{0}+2}R)^{\frac{\beta p}{p-1}} |B(x,2^{k} \varepsilon)|^{-\frac{1}{p-1}} \notag \\
    & \leq  C(c,Q) \sum_{k=0}^{k_{0}}(2^{k-k_{0}+2}R)^{\frac{\beta p}{p-1}} \left(c2^{(k-k_{0})Q}|B(x,R)|\right)^{-\frac{1}{p-1}} \notag \\
    &= C(c,Q) c^{-\frac{1}{p-1}}|B(x,R)|^{-\frac{1}{p-1}} 4^{\frac{\beta p} {p-1}} R^{\frac{\beta p} {p-1}} 2^{-k_{0}(\frac{\beta p -Q} {p-1} )} \sum_{k=0}^{k_{0}} (2^{\frac{\beta p -Q} {p-1}})^{k}.
\end{align}
From $1< p < \frac{Q}{\beta}$, we have $\beta p -Q < 0$, and
\begin{align}\label{9.17}
    \sum_{k=0}^{k_{0}} (2^{\frac{\beta p -Q} {p-1}})^{k} = \frac{1-(2^{\frac{\beta p-Q}{p-1}})^{k_{0}+1}}{1-2^{\frac{\beta p-Q}{p-1}}} =\frac{2^{\frac{Q-\beta p}{p-1}}}{2^{\frac{Q-\beta p}{p-1}} -1} \left(1-2^{(\frac{\beta p-Q}{p-1})(k_{0}+1)}\right) < \frac{2^{\frac{Q-\beta p}{p-1}}}{2^{\frac{Q-\beta p}{p-1}} -1}.
\end{align}
Putting (\ref{9.10}), (\ref{9.16}) and (\ref{9.17}) together we get
\begin{align}\label{9.18}
    \left(\int_{B^{c}(x,\varepsilon)\cap B(x,2R)} \right. &\left.  \left(\frac{d(x, y)^\beta}{|B(x, d(x, y))|}\right)^{\frac{p}{p-1}} d y\right)^{\frac{p-1}{p}} \\
    & \leq C(c,Q)|B_{R}|^{-\frac{1}{p}} R^{\beta} 2^{-k_{0}(\frac{\beta p -Q} {p})}\left(\frac{2^{\frac{Q-\beta p}{p-1}}}{2^{\frac{Q-\beta p}{p-1}} -1}\right)^{\frac{p-1}{p}} .\notag
\end{align}
Since $2^{k_{0}} \leq \frac{2R}{\varepsilon}$,  we have
\begin{align}\label{9.19}
    2^{-k_{0}(\frac{\beta p -Q} {p})} \leq \left(\frac{2R}{\varepsilon}\right)^{\frac{Q-\beta p}{p}} = 2^{\frac{Q-\beta p}{p}} \varepsilon^{\frac{\beta p- Q}{p}}R^{\frac{Q-\beta p}{p}} \leq 2^{Q}\varepsilon^{\frac{\beta p- Q}{p}}R^{\frac{Q-\beta p}{p}}.
\end{align}
Combining (\ref{9.14}), (\ref{9.18}) with (\ref{9.19}), we obtain
\begin{align}\label{9.20}
    I^{2}_\beta(f)(x) &\leq C(c,Q) \left(\frac{2^{\frac{Q-\beta p}{p-1}}}{2^{\frac{Q-\beta p}{p-1}} -1}\right)^{\frac{p-1}{p}}R^{\frac{Q}{p}}|B_{R}|^{-\frac{1}{p}} \|f\|_{L^{p}(B_{R})} \varepsilon^{\frac{\beta p- Q}{p}}.
\end{align}
Adding (\ref{9.13}) and (\ref{9.20}) we get
\begin{align}\label{9.21}
    I_\beta(f)(x) \leq & C(c,Q) \frac{2^{\beta}}{2^{\beta} -1}Mf(x)\varepsilon^{\beta} \\
    &+C(c,Q) \left(\frac{2^{\frac{Q-\beta p}{p-1}}}{2^{\frac{Q-\beta p}{p-1}} -1}\right)^{\frac{p-1}{p}}R^{\frac{Q}{p}}|B_{R}|^{-\frac{1}{p}} \|f\|_{L^{p}(B_{R})} \varepsilon^{\frac{\beta p- Q}{p}} .\notag
\end{align}
Let $A := C(c,Q) \frac{2^{\beta}}{2^{\beta} -1}Mf(x) $ and $B:= C(c,Q) \left(\frac{2^{\frac{Q-\beta p}{p-1}}}{2^{\frac{Q-\beta p}{p-1}} -1}\right)^{\frac{p-1}{p}}R^{\frac{Q}{p}}|B_{R}|^{-\frac{1}{p}} \|f\|_{L^{p}(B_{R})}$. Minimizing the right-hand side of (\ref{9.21}) with respect to $\varepsilon$, we find that $\varepsilon$ is given by
\begin{align*}
    \varepsilon = \left(\frac{(Q-\beta p)}{\beta p } \frac{B}{A}\right)^{\frac{p}{Q}}.
\end{align*}
Substituting the above equality into (\ref{9.21}) yields
\begin{align}\label{9.22}
    I_\beta(f)(x) \leq C(c,Q)C(\beta,p,Q)|B_{R}|^{-\frac{\beta}{Q}} R^{\beta} |Mf(x)|^{\frac{Q-\beta p}{Q}} \|f\|^{\frac{\beta p}{Q}}_{L^{p}(B_{R})}
\end{align}
where
\begin{align}\label{9.23}
    C(\beta,p, Q) := \left(\frac{Q - \beta p}{\beta p}\right)^{\frac{\beta p}{Q}}\frac{Q}{Q-\beta p} \left(\frac{2^{\frac{Q-\beta p}{p-1}}}{2^{\frac{Q-\beta p}{p-1}} -1}\right)^{\frac{\beta(p-1)}{Q}}\left(\frac{2^{\beta}}{2^{\beta}-1}\right)^{\frac{Q-\beta p}{Q}}.
\end{align}
From $\frac{1}{p} - \frac{1}{q} = \frac{\beta}{Q}$, Proposition \ref{prop9.1} and (\ref{9.22}), we have
\begin{align}\label{9.24}
    \left(\fint_{B_{R}} |I_\beta(f)(x)|^{q}dx\right)^{\frac{1}{q}}  &\leq C(c,Q)C(\beta,p,Q)R^{\beta} |B_{R}|^{-\frac{1}{p}}\|f\|_{L^{p}(B_{R})}^{\frac{\beta p}{Q}} \left(\int_{B_{R}}|Mf(x)|^{p}  dx\right)^{\frac{1}{q}}  \notag \\
    & \leq C(c, Q)C(\beta,p, Q) \left(\frac{p2^{p-1}}{p-1}\right)^{\frac{1}{q}}R^{\beta}\left( \fint_{B_{R}}|f(y)|^{p} dy\right)^{\frac{1}{p}}.
\end{align}
Denote
\begin{eqnarray}\label{9.25}
    \overline{C} =
    \begin{cases}
        \frac{C(c,Q)}{\left(2^{\frac{\beta pq - Q(q-p)}{pq -q+p}} -1\right)^{\frac{pq -q +p}{pq}}} , \quad 0\leq \frac{1}{p}-\frac{1}{q} <\frac{\beta}{Q}\\
        C(c, Q)C(\beta,p, Q) \left(\frac{p2^{p-1}}{p-1}\right)^{\frac{1}{q}} ,\quad \frac{1}{p}-\frac{1}{q} = \frac{\beta}{Q}
    \end{cases}
\end{eqnarray}
where $C(\beta,p, Q)$ is given by (\ref{9.23}). Combining (\ref{9.12}) and (\ref{9.24}), the proof of Theorem \ref{thm9.4} is complete.
\hfill${\square}$

\noindent {\bf Proof of Theorem \ref{thmS}}: Let $1 \leq \kappa \leq \frac{Q}{Q-p}$ be defined by $q= \kappa p$. For $\beta =1$ and $\frac{1}{p}-\frac{1}{q} <\frac{1}{Q}$, from (\ref{9.12}), we have
\begin{align}\label{9.26}
    \left(\fint_{B_{R}} |I_{1}(f)(x)|^{\kappa p} dx\right)^{\frac{1}{\kappa p}} \leq  \frac{C(c,Q) R}{\left(2^{\frac{Q-p}{\kappa(p-1)+1}(\frac{Q}{Q-p} -\kappa)}-1\right)^{\frac{\kappa(p-1)+1}{\kappa p}}}  \left(\fint_{B_{R}} |f(y)|^{p} dy\right)^{\frac{1}{p}}.
\end{align}

For $\beta =1$ and $\frac{1}{p}-\frac{1}{q} = \frac{1}{Q}$, from (\ref{9.23}) and (\ref{9.24}), we have
\begin{align}\label{9.27}
    \left(\fint_{B_{R}} |I_{1}(f)(x)|^{\kappa p} dx \right)^{\frac{1}{q}}
    \leq  \frac{C(c, Q)}{Q-p}\left(1 + \frac{1}{2^{\frac{Q-p}{p-1}} -1 }\right)^{\frac{Q-1}{Q}} \left(\frac{p2^{p-1}}{p-1}\right)^{\frac{Q-p}{pQ}}R\left( \fint_{B_{R}}|f(y)|^{p} dy\right)^{\frac{1}{p}}
\end{align}
where we used
\begin{align*}
    C(1,p,Q) &= \left(\frac{Q-p}{p}\right)^{\frac{p}{Q} }\frac{Q}{Q-p}\left(\frac{2^{\frac{Q-p}{p-1}}}{2^{\frac{Q- p}{p-1}} -1}\right)^{\frac{p-1}{Q}} 2^{\frac{Q-p}{Q}} < \frac{2Q^{2}}{Q-p}\left(1 + \frac{1}{2^{\frac{Q-p}{p-1}} -1 }\right)^{\frac{Q-1}{Q}}.
\end{align*}
By (\ref{9.7}), we have
\begin{align}\label{9.28}
    \left(\fint_{B_{R} } |u(x)|^{\kappa p}\right)^{\frac{1}{\kappa p}} &\leq C \left(\fint_{B_{R}}| I_1\left(\left|\nabla_{0} u\right|\right)(x)|^{\kappa p} dx \right)^{\frac{1}{\kappa p }} .
\end{align}
Combining (\ref{9.26})-(\ref{9.28}), we arrive at
\begin{align*}
    \left(\fint_{B_{R} } |u(x)|^{p\kappa}\right)^{\frac{1}{\kappa p}}  \leq \mathcal{C}_{s} R \left(\fint_{B_{R}}|\nabla_{0} u|^{p} dy\right)^{\frac{1}{p}}.
\end{align*}
This is completes the proof of Theorem \ref{thmS}.
\hfill${\square}$

\subsection{ Proofs of Lemmas \ref{lem5.8}-\ref{lem5.11}}
\
   \vglue-10pt
    \indent

\noindent {\bf Proof Lemma \ref{lem5.8} }: Let $\varphi = \eta^{\beta+2} |Tw_{p}|^{\beta+1} $. From the regular result (\ref{5.10}), we have $\varphi \in HW^{1,2}_{0}(D)$. Using $\varphi$ as a test function in the weak formula of  (\ref{5.46}) we get
\begin{align*}
    & (\beta+1) \int_D \sum_{i=1}^2\left(\sum_{j=1}^2 \partial_{z_j} a_{i} X_j T w_{p}\right)|T w_{p}|^\beta X_i T w_{p} \eta^{\beta+2} d x \\
    = & -(\beta+2) \int_D \sum_{i=1}^2\left(\sum_{j=1}^2 \partial_{z_j} a_{i} X_j T w_{p}\right) X_i \eta \eta^{\beta+1}|T w_{p}|^{\beta+1} d x .
\end{align*}
Using (\ref{5.42}) and (\ref{5.43}), we have
\begin{align*}
    (p-1) \int_D \eta^{\beta+2}\Phi^{\frac{p-2}{2}}|T w_{p}|^\beta\left|\nabla_{0} T w_{p}\right|^2 d x
    \leq   \frac{\beta+1}{\beta+2} \int_D\Phi^{\frac{p-2}{2}}\left|\nabla_{0} T w_{p}\right|\left|\nabla_{0} \eta\right| \eta^{\beta+1}|T w_{p}|^{\beta+1} d x .
\end{align*}
Applying H$\ddot{\rm o} $lder's inequality with exponent 2 to the right-hand side of the above inequality, we complete the proof of Lemma \ref{lem5.8}.
\hfill${\square}$

\noindent {\bf Proof Lemma \ref{lem5.9}}: Using $\varphi = \eta^{2} \Phi^{\frac{\beta}{2}} X_{1}w_{p}$ as a test function in the weak formulation of (\ref{5.44}), we have
\begin{align*}
    I_1= & \int_D \eta^2 \sum_{i, j=1}^2 \partial_{z_j} a_i X_j X_1 w_{p} X_i\left(\Phi^{\frac{\beta}{2}} X_1 w_{p}\right) d x \\
    = & -2 \int_D \sum_{i, j=1}^2 \partial_{z_j} a_i X_j X_1 w_{p} X_i \eta \eta \Phi^{\frac{\beta}{2}} X_1 w_{p} d x  -\int_D \sum_{i=1}^2 \partial_{z_2} a_i T w_{p} X_i\left(\eta^2\Phi^{\frac{\beta}{2}} X_1 w_{p}\right) d x  \\
    &-\int_D a_2 T\left(\eta^2\Phi^{\frac{\beta}{2}} X_1 w_{p}\right) d x  =  I_2+I I+I I I
\end{align*}
We estimate the terms $I_2$, $I I$, and $I I I$, respectively. By (\ref{5.43}), we obtain
\begin{align}\label{9.29}
    I_2  \leq 2\int_D  \eta\left|\nabla_{0} \eta\right| \Phi^{\frac{p-1+\beta}{2}} \left|\nabla^{2}_{0} w_{p}\right| d x
\end{align}
Now calculating derivatives and using (\ref{5.43}), we get
\begin{align}\label{9.30}
    II= & 2 \int_D \sum_{i=1}^2 \partial_{z_2} a_i T w_{p} X_i \eta \eta\Phi^{\frac{\beta}{2}} X_1 w_{p} d x  + \int_D \eta^2 \sum_{i=1}^2 \partial_{z_2} a_i T w_{p}\Phi^{\frac{\beta}{2}} X_i X_1 w_{p} d x \notag \\
    & +\beta \int_D \eta^2 \sum_{i=1}^2 \partial_{z_2} a_i T w_{p}\Phi^{\frac{\beta-2}{2}}\left\langle X_i \nabla_{0} w_{p}, \nabla_{0} w_{p}\right\rangle X_1 w_{p} d x  \notag \\
    \leq & 2 \int_{D} \eta\left|\nabla_{0} \eta\right|\Phi^{\frac{p-1+\beta}{2}}|T w_{p}| d x  + (\beta +1 )\int_D \eta^2\Phi^{\frac{p-2+\beta}{2}}|T w_{p}|\left|\nabla^{2}_{0} w_{p}\right| d x
\end{align}
The term $III$ is estimated as follows.
\begin{align*}
    III=&  2 \int_D a_2 \eta T \eta\Phi^{\frac{\beta}{2}} X_1 w_{p} d x +\beta \int_D \eta^2 a_2\Phi^{\frac{\beta-2}{2}}\left\langle T \nabla_{0} w_{p}, \nabla_{0} w_{p}\right\rangle X_1 w_{p} d x  \\
    &+ \int_D \eta^2 a_2\Phi^{\frac{\beta}{2}} T X_1 w_{p} d x \\
    =&  III_1+III_2+III_3.
\end{align*}
Using (\ref{5.41}), we get
\begin{align}
    III_1 \leq 2 \int_D \eta |T\eta| \Phi^{\frac{p+\beta}{2}}  d x.
\end{align}
To estimate $III_{2}$ and $III_{3}$, we integrate by parts and apply (\ref{5.41}) and (\ref{5.43}),
\begin{align}
    III_2= & \beta \int_D \eta^2 a_2\Phi^{\frac{\beta-2}{2}} \sum_{k=1}^2 X_k T w_{p} X_k w_{p} X_1 w_{p} d x \notag \\
    &=  -\beta \int_D \sum_{k=1}^2 X_k\left(\eta^2 a_2\Phi^{\frac{\beta-2}{2}} X_k w_{p} X_1 w_{p}\right) T w_{p} d x \notag \\
    \leq & 2\beta \int_D \eta\left|\nabla_{0} \eta\right| \Phi^{\frac{p-1+\beta}{2}} |T w_{p}| d x  +\beta(\beta+1) \int_D \eta^2\Phi^{\frac{p-2+\beta}{2}}|T w_{p}| \left|\nabla_{0}^2 w_{p}\right| d x  \label{9.32} \\
    III_3= & \int_D \eta^2 a_2\Phi^{\frac{\beta}{2}} X_1 T w_{p} d x =  -\int_D X_1\left(\eta^2 a_2\Phi^{\frac{\beta}{2}}\right) T w_{p} d x \notag \\
    \leq & 2 \int_D \eta|\nabla_{0} \eta| \Phi^{\frac{p-1+\beta}{2}} |Tw_{p}| dx +(\beta +1) \int_D \eta^{2} \Phi^{\frac{p-2+\beta}{2}}|Tw_{p}| |\nabla^{2}_{0} w_{p}|^{2} dx \label{9.33}.
\end{align}
Adding (\ref{9.29})-(\ref{9.33}), we get
\begin{align*}
   I_{2}+ II + III \leq& 2 \int_D \eta|\nabla_{0} \eta| \Phi^{\frac{p-1+\beta}{2}} |\nabla_{0}^{2} w_{p}| dx +  2(\beta +2) \int_D \eta|\nabla_{0} \eta| \Phi^{\frac{p-1+\beta}{2}} |Tw_{p}| dx\\
    &+ 2 \int_D \eta|T \eta| \Phi^{\frac{p+\beta}{2}} dx  + (\beta +1)(\beta +2) \int_D \eta^{2} \Phi^{\frac{p-2+\beta}{2}} |\nabla^{2}_{0} w_{p}||Tw_{p}| dx.
\end{align*}
Using Young's inequality with exponent 2 introducing a parameter $\varepsilon>0$ that will be chosen later, and $|T w_{p}| \leq 2\left|\nabla_{0}^2 w_{p}\right|$, we get
\begin{align}\label{9.34}
    I_1= & \int_D \eta^2 \sum_{i, j=1}^2 \partial_{z_j} a_i X_j X_1 w_{p} X_i\left(\Phi^{\frac{\beta}{2}} X_1 w_{p}\right) d x \notag \\
    \leq& \varepsilon(\beta^{2}+11\beta +20) \int_D \eta^2\Phi^{\frac{p-2+\beta}{2}}\left|\nabla_{0}^2 w_{p}\right|^2 d x  +\frac{(\beta+1)(\beta +2)}{4 \varepsilon} \int_D \eta^2\Phi^{\frac{p-2+\beta}{2}}|T w_{p}|^2 d x  \\
    & +\left(\frac{\beta+3}{2 \varepsilon} \left\|\nabla_{0} \eta\right\|^2_{L^{\infty}(D)} + 2 \left\|\eta T\eta\right\|_{L^{\infty}(D)}\right)\int_D\Phi^{\frac{p+\beta}{2}} d x \notag
\end{align}
Similarly, using $\varphi=\eta^2\Phi^{\frac{\beta}{2}} X_2 w_{p}$ as a test function in the weak formulation of (\ref{5.45}) we get
\begin{align}\label{9.35}
    J_1= & \int_D \eta^2 \sum_{i, j=1}^2 \partial_{z_j} a_i X_j X_2 w_{p} X_i\left(\Phi^{\frac{\beta}{2}} X_2 w_{p}\right) d x \notag \\
    \leq& \varepsilon(\beta^{2}+11\beta +20) \int_D \eta^2\Phi^{\frac{p-2+\beta}{2}}\left|\nabla_{0}^2 w_{p}\right|^2 d x  +\frac{(\beta+1)(\beta +2)}{4 \varepsilon} \int_D \eta^2\Phi^{\frac{p-2+\beta}{2}}|T w_{p}|^2 d x  \\
    & +\left(\frac{\beta+3}{2 \varepsilon} \left\|\nabla_{0} \eta\right\|^2_{L^{\infty}(D)} + 2 \left\|\eta T\eta\right\|_{L^{\infty}(D)}\right)\int_D\Phi^{\frac{p+\beta}{2}} d x \notag.
\end{align}
Now we will estimate $I_1$ and $J_1$. For $j=1,2$, we have
\begin{align*}
    &X_{i}\left(\Phi^{\frac{\beta}{2}}X_{j}w_{p}\right)= \frac{\beta}{2}\Phi^{\frac{\beta -2}{2}}X_{i}(|\nabla_{0} w_{p}|^{2})X_{j} w_{p} + \Phi^{\frac{\beta}{2}}X_{i}X_{j}w_{p}, \\
    &X_j X_1 w_{p} X_1 w_{p}+X_j X_2 w_{p} X_2 w_{p}=\frac{1}{2} X_j\left(\left|\nabla_{0} w_{p}\right|^2\right).
\end{align*}
From the above equalities and  (\ref{5.42}), we get
\begin{align*}
    I_1 \geq &  (p-1)\int_D \eta^2\Phi^{\frac{\beta-2+\beta}{2}}\left|\nabla_{0} X_1 w_{p}\right|^2 d x  \\
    &+\frac{\beta}{2} \int_D \eta^2 \sum_{i, j=1}^2 \partial_{z_j} a_i X_j X_1 w_{p} X_1 w_{p} X_i\left(\left|\nabla_{0} w_{p}\right|^2\right)\Phi^{\frac{\beta-2}{2}} d x
\end{align*}
and
\begin{align*}
    J_1 \geq  (p-1)\int_D \eta^2\Phi^{\frac{\beta-2+a}{2}}\left|\nabla_{0} X_2 w_{p}\right|^2 d x  +\frac{\beta}{2} \int_D \eta^2 \sum_{i, j=1}^2 \partial_{z_j} a_iX_j X_2 w_{p} X_2 w_{p} X_i\left(\left|\nabla_{0} w_{p}\right|^2\right)\Phi^{\frac{a-2}{2}} d x
\end{align*}
Thus,
\begin{align}\label{9.36}
    I_1+J_1
    \geq & (p-1)\int_D \eta^2\Phi^{\frac{p-2+\beta}{2}}\left|\nabla_{0}^2 w_{p}\right|^2 d x +\frac{\beta (p-1)}{4} \int_D \eta^2\Phi^{\frac{p-2}{2}}\left|\nabla_{0}\left(\left|\nabla_{0} w_{p}\right|^2\right)\right|^2\Phi^{\frac{\beta-2}{2}} d x \notag \\
    \geq & (p-1)(\beta +1)\int_D \eta^2\Phi^{\frac{p-2+\beta}{2}}\left|\nabla_{0}^2 w_{p}\right|^2 d x .
\end{align}
Combining (\ref{9.34})-(\ref{9.36}), there exist absolute constants $C_{5}$, $C_{6}$ and $C_{7}$ such that
\begin{align}\label{9.37}
    \int_D  \eta^2 \Phi^{\frac{p-2+\beta}{2}}\left|\nabla_{0}^2 w_{p}\right|^2 d x \leq& \frac{C_{5}(\beta +1)\varepsilon }{p-1} \int_D \eta^2\Phi^{\frac{p-2+\beta}{2}}\left|\nabla_{0}^2 w_{p}\right|^2 d x  \\
    &+ \frac{C_{6}(\beta +1)}{2(p-1) \varepsilon} \int_D \eta^2\Phi^{\frac{p-2+\beta}{2}}|T w_{p}|^2 d x  \notag \\
    &+ \frac{C_{7}}{p-1}\left(\left\|\nabla_{0} \eta\right\|^2_{L^{\infty}(D)}+\left\| \eta T\eta\right\|_{L^{\infty}(D)}\right) \int_D\Phi^{\frac{p+\beta}{2}} d x  \notag
\end{align}
By choosing $\varepsilon=\frac{p-1}{2C_{5}(\beta +1)}$, the first integral on the right-hand side can be absorbed into the left-hand side of (\ref{9.37}). This completes the proof of Lemma \ref{lem5.9}.
\hfill${\square}$

\noindent {\bf Proof Lemma \ref{lem5.10}}: From the equation (4.57) in \cite{R2015} and integration by parts, we have
\begin{align*}
    \int_D \sum_{i=1}^2 X_1\left(a_i\right) X_i \varphi d x-\int_D T\left(a_2\right) \varphi d x=0 \quad \text { for all } \varphi \in H W_0^{1,2}(D).
\end{align*}
Let $\varphi=\eta^{\beta+2}|T w_{p}|^\beta X_1 w_{p}$ as a test function, we have
\begin{align}\label{9.38}
    I_1  = &\int_D \sum_{i=1}^2 X_1\left(a_i\right) X_i X_1 w_{p} \eta^{\beta+2}|T w_{p}|^\beta d x \notag \\
    =& -(\beta+2) \int_D \sum_{i=1}^2 X_1\left(a_i\right) X_i \eta \eta^{\beta+1}|T w_{p}|^\beta X_1 w_{p} d x\\
    & -\beta \int_D \sum_{i=1}^2 X_1\left(a_i\right) X_i T w_{p}|T w_{p}|^{\beta-1} \operatorname{sign}(T w_{p}) X_1 w_{p} \eta^{\beta+2} d x \notag \\
    & +\int_D T\left(a_2\right) \eta^{\beta+2}|T w_{p}|^\beta X_1 w_{p} d x=I_2+I_3+II \notag.
\end{align}
Using $X_2 X_1=X_1 X_2-T$, we get
\begin{align*}
    I_1 & =\int_D \eta^{\beta+2} \sum_{i=1}^2 X_1\left(a_i\right) X_1 X_i w_{p} |T w_{p}|^\beta d x-\int_D \eta^{\beta+2} X_1\left(a_2\right) T w_{p} |T w_{p}|^\beta d x =I_{1,1}-I_{1,2}.
\end{align*}
The equation (\ref{9.38}) is equivalent to $I_{1,1}=I_{1,2}+I_2+I_3+I I$. By (\ref{5.42}), there holds
\begin{align}\label{9.39}
    I_{1,1} \geq (p-1) \int_D \eta^{\beta+2}\Phi^{\frac{p-2}{2}}\left|\nabla_{0} X_1 w_{p}\right|^2|T w_{p}|^\beta d x.
\end{align}
For $I_{1,2}$, the regularity result (\ref{5.10}) ensures that integration by parts is valid, and hence
\begin{align}\label{9.40}
    I_{1,2} &= \int_{D} a_{2}X_{1}\left(\eta^{\beta+2}Tw_{p}|Tw_{p}|^{\beta}\right) \notag \\ &\leq  (\beta+2)  \int_D \eta^{\beta+1}\Phi^{\frac{p-1}{2}}\left|\nabla_{0} \eta\right||T w_{p}|^{\beta+1} d x +(\beta+1)  \int_D \eta^{\beta+2}\Phi^{\frac{p-1}{2}}|T w_{p}|^\beta\left|\nabla_{0} T w_{p}\right| d x
\end{align}
where we used (\ref{5.41}). Analogously, using (\ref{5.41}) and (\ref{5.43}) to get
\begin{align*}
    I_2 &\leq (\beta +2) \int_{D} \eta^{\beta+1}\Phi^{\frac{p-1}{2}}|\nabla_{0} \eta| |Tw_{p}|^{\beta} |\nabla^{2}_{0} w_{p}|dx, \\
    I_3 &\leq  \beta \int_{D} \eta^{\beta+2}\Phi^{\frac{p-1}{2}} |\nabla_{0} (Tw_{p})||Tw_{p}|^{\beta -1} |\nabla^{2}_{0} w_{p}|,dx\\
    I I &\leq \int_{D} \eta^{\beta+2}\Phi^{\frac{p-1}{2}} |\nabla_{0} (Tw_{p})||Tw_{p}|^{\beta} dx.
\end{align*}
Adding (\ref{9.40}) and the above three inequalities, we get
\begin{align*}
    I_{1,1}\leq& (\beta +2) \int_{D} \eta^{\beta+1}\Phi^{\frac{p-1}{2}}\left|\nabla_{0}\eta\right||T w_{p}|^{\beta+1} d x  + (\beta +2) \int_{D} \eta^{\beta+1}\Phi^{\frac{p-1}{2}}\left|\nabla_{0}\eta\right||\nabla_{0}^{2} w_{p}| |T w_{p}|^{\beta} d x  \\
    & + (\beta +2) \int_{D} \eta^{\beta+2}\Phi^{\frac{p-1}{2}}|\nabla_{0} (Tw_{p})| |T w_{p}|^{\beta} d x \\
    &+ \beta \int_{D} \eta^{\beta+2}\Phi^{\frac{p-1}{2}}|\nabla_{0} (Tw_{p})||\nabla_{0}^{2} w_{p}| |T w_{p}|^{\beta-1} d x
\end{align*}
By $|T w_{p}| \leq 2\left|\nabla_{0}^2 w_{p}\right|$ and Young's inequality with a parameter $\varepsilon>0$ to be chosen later, we get
\begin{align}\label{9.41}
    I_{1,1}
    \leq & 5(\beta +2)\varepsilon \int_{D} \eta^{\beta+2}\Phi^{\frac{p-2}{2}}|T w_{p}|^{\beta} |\nabla^{2}_{0} w_{p}|^{2}d x \\
    &+ \frac{2\varepsilon(\beta +1)}{\left\|\nabla_{0} \eta\right\|^{2}_{L^{\infty}(D)}}  \int_{D} \eta^{\beta+4}\Phi^{\frac{p-2}{2}}|T w_{p}|^{\beta} |\nabla_{0} (w_{p})|^{2}d x \notag \\
    &+\left(3(\beta +2) + \frac{\beta}{4}\right)\frac{\left\|\nabla_{0} \eta\right\|^{2}_{L^{\infty}(D)}}{\varepsilon}\int_{D} \eta^{\beta}\Phi^{\frac{p}{2}}|T w_{p}|^{\beta-2} |\nabla^{2}_{0} w_{p}|^{2}d x \notag
\end{align}
From Lemma \ref{lem5.8}, there exists an absolute constant $C$ such that
\begin{align}\label{9.42}
    \int_{D}  \eta^{\beta+4}\Phi^{\frac{p-2}{2}} |Tw_{p}|^{\beta} |\nabla_{0}(Tw_{p})|^{2}&\leq  \frac{(\beta+4)^{2}\left\|\nabla_{0} \eta\right\|^{2}_{L^{\infty}(D)}}{(\beta+1)^{2}(p-1)^{2}} \int_{D}\eta^{\beta}\Phi^{\frac{p-2}{2}} |Tw_{p}|^{\beta +2}  \notag \\
    &\leq  \frac{C\left\|\nabla_{0} \eta\right\|^{2}_{L^{\infty}(D)}}{(p-1)^{2}} \int_{D}\eta^{\beta}\Phi^{\frac{p-2}{2}} |Tw_{p}|^{\beta +2}
\end{align}
Putting (\ref{9.39}), (\ref{9.41}) and (\ref{9.42}) together we get
\begin{align}\label{9.43}
    (p-1)\int_D \eta^{\beta+2} & \Phi^{\frac{p-2}{2}}\left|X_1 \nabla_{0} w_{p}\right|^2|T w_{p}|^\beta d x \\
    & \leq \varepsilon \left(\frac{C(\beta +1)}{(p-1)^{2}} + 5(\beta +2)\right)  \int_D \eta^{\beta+2}\Phi^{\frac{p-2}{2}}|T w_{p}|^\beta\left|\nabla_{0}^2 w_{p}\right|^2 d x  \notag \\
    & +\left(3(\beta +2) + \frac{\beta}{4}\right)\frac{1}{\varepsilon}\left\|\nabla_{0} \eta\right\|_{L^{\infty}(D)}^2 \int_D \eta^\beta\Phi^{\frac{p}{2}}|T w_{p}|^{\beta-2}\left|\nabla_{0}^2 w_{p}\right|^2 d x \notag
\end{align}
Similarly, let $\varphi=\eta^{\beta+2}|T w_{p}|^\beta X_2 w_{p}$ as a test function in the weak formulation of (\ref{5.45}), we obtain
\begin{align}\label{9.44}
    (p-1) \int_D \eta^{\beta+2} & \Phi^{\frac{p-2}{2}}\left|X_2 \nabla_{0} w_{p}\right|^2|T w_{p}|^\beta d x  \\
     \leq & \varepsilon \left(\frac{C(\beta +1)}{(p-1)^{2}} + 5(\beta +2)\right)  \int_D \eta^{\beta+2}\Phi^{\frac{p-2}{2}}|T w_{p}|^\beta\left|\nabla_{0}^2 w_{p}\right|^2 d x \notag \\
    & +\left(3(\beta +2) + \frac{\beta}{4}\right)\frac{1}{\varepsilon}\left\|\nabla_{0} \eta\right\|_{L^{\infty}(D)}^2 \int_D \eta^\beta\Phi^{\frac{p}{2}}|T w_{p}|^{\beta-2}\left|\nabla_{0}^2 w_{p}\right|^2 d x\notag
\end{align}
Adding (\ref{9.43}) and (\ref{9.44}), we have
\begin{align*}
    \int_D \eta^{\beta+2}&\Phi^{\frac{p-2}{2}}|T w_{p}|^\beta \left|\nabla_{0}^2 w_{p}\right|^2 d x \\
    \leq & \left(\frac{C(\beta +1)}{(p-1)^{3}} + \frac{10(\beta +2)}{p-1}\right) \varepsilon \int_D \eta^{\beta+2}\Phi^{\frac{p-2}{2}}|T w_{p}|^\beta\left|\nabla_{0}^2 w_{p}\right|^2 d x \\
    & +\left(\frac{6(\beta +2)}{p-1} + \frac{\beta}{2(p-1)}\right)\frac{1}{\varepsilon}\left\|\nabla_{0} \eta\right\|_{L^{\infty}(D)}^2 \int_D \eta^\beta\Phi^{\frac{p}{2}}|T w_{p}|^{\beta-2}\left|\nabla_{0}^2 w_{p}\right|^2 d x\\
    \leq & \frac{C_{8}(\beta +1)\varepsilon }{(p-1)^{3}} \int_D \eta^{\beta+2}\Phi^{\frac{p-2}{2}}|T w_{p}|^\beta\left|\nabla_{0}^2 w_{p}\right|^2 d x \\
    & + \frac{C_{9}(\beta +1)}{\varepsilon(p-1)}\left\|\nabla_{0} \eta\right\|_{L^{\infty}(D)}^2 \int_D \eta^\beta\Phi^{\frac{p}{2}}|T w_{p}|^{\beta-2}\left|\nabla_{0}^2 w_{p}\right|^2 d x
\end{align*}
where $C_{8}$ and $C_{9}$ are absolute constants. Choose $\varepsilon=\frac{(p-1)^{3}}{2C_{8}(\beta +1)} $, we obtain
\begin{align}\label{9.45}
    \int_D \eta^{\beta+2} \Phi^{\frac{p-2}{2}}|T w_{p}|^\beta &\left|\nabla_{0}^2 w_{p}\right|^2 d x \notag \\
    &\leq 2C_{8}C_{9}\frac{(\beta +1)^{2}}{(p-1)^{4}}\left\|\nabla_{0} \eta\right\|_{L^{\infty}(D)}^2 \int_D \eta^\beta\Phi^{\frac{p}{2}}|T w_{p}|^{\beta-2}\left|\nabla_{0}^2 w_{p}\right|^2 d x .
\end{align}
Using H$\ddot{\rm o} $lder's inequality with exponent $\frac{\beta}{\beta-2}$, we have
\begin{align*}
    \int_D &\eta^\beta\Phi^{\frac{p}{2}}|T w_{p}|^{\beta-2}\left|\nabla_{0}^2 w_{p}\right|^2 d x
    \\
    &\leq \left(\int_D \eta^{\beta+2}\Phi^{\frac{p-2}{2}}|T w_{p}|^\beta\left|\nabla_{0}^2 w_{p}\right|^2 d x\right)^{\frac{\beta-2}{\beta}} \left(\int_D \eta^2\Phi^{\frac{p-2+\beta}{2}}\left|\nabla_{0}^2 w_{p}\right|^2 d x\right)^{\frac{2}{\beta}}.
\end{align*}
Combining the above inequality with (\ref{9.45}), the proof of Lemma \ref{lem5.10} is complete.
\hfill${\square}$

\noindent {\bf Proof Lemma \ref{lem5.11}}: Using $|T w_{p}|^{2} \leq 4|\nabla_{0}^{2} w_{p}|^{2}$ and H$\ddot{\rm o} $lder’s inequality with exponent $\frac{\beta +2}{2}$, we obtain
\begin{align*}
    \int_{D} \eta^{2}\Phi^{\frac{p-2+\beta}{2}} |Tw_{p}|^{2}
    \leq &\left(\int_{D} \eta^{\beta + 2} \Phi^{\frac{p-2}{2}} |Tw_{p}|^{\beta +2}\right)^{\frac{2}{\beta +2}} \left(\int_{D}\Phi^{\frac{p+\beta}{2}}  \right)^{\frac{\beta}{\beta +2}} \\
    \leq & \left( 4 \int_{D} \eta^{\beta + 2} \Phi^{\frac{p-2}{2}} |Tw_{p}|^{\beta}|\nabla^{2}_{0} w_{p}|^{2}\right)^{\frac{2}{\beta +2}} \left(\int_{D}\Phi^{\frac{p+\beta}{2}}  \right)^{\frac{\beta}{\beta +2}}.
\end{align*}
By Lemma \ref{lem5.10}, we get
\begin{align*}
    \int_{D} \eta^{2}\Phi^{\frac{p-2+\beta}{2}} |Tw_{p}|^{2}
    \leq &\left[4 \left(\frac{C(\beta +1)^{2}}{(p-1)^{4}}\right)^{\frac{\beta}{2}}\left\|\nabla_{0} \eta\right\|^{\beta}_{L^{\infty}(D)} \int_{D} \eta^{2}\Phi^{\frac{p-2+\beta}{2}}  |\nabla^{2}_{0}w_{p}|^{2}\right]^{\frac{2}{\beta +2}} \left(\int_{D}\Phi^{\frac{p+\beta}{2}}  \right)^{\frac{\beta}{\beta +2}}\notag \\
    \leq &C\left(\frac{(\beta +1)^{2}}{(p-1)^{4}}\right)^{\frac{\beta}{\beta +2}}\left\|\nabla_{0} \eta\right\|^{\frac{2\beta}{\beta +2}}_{L^{\infty}(D)} \left(\int_{D} \eta^{2}\Phi^{\frac{p-2+\beta}{2}}  |\nabla^{2}_{0}w_{p}|^{2}\right)^{\frac{2}{\beta +2}} \left(\int_{D}\Phi^{\frac{p+\beta}{2}}  \right)^{\frac{\beta}{\beta +2}}.
\end{align*}
Using Young's inequality with a parameter $\varepsilon>0$ to be chosen later, we have
\begin{align*}
    \int_{D} \eta^{2}\Phi^{\frac{p-2+\beta}{2}}  |Tw_{p}|^{2}
    \leq \frac{ 2\varepsilon }{\beta +2}\int_{D}  \eta^{2}\Phi^{\frac{p-2+\beta}{2}} |\nabla^{2}_{0}w_{p}|^{2} dx + \frac{C \beta}{\beta +2} \varepsilon^{-\frac{2}{\beta}} \frac{(\beta +1)^{2}}{(p-1)^{4}}\left\|\nabla_{0} \eta\right\|^{2}_{L^{\infty}(D)}\int_{D}  \Phi^{\frac{p+\beta}{2}} dx
\end{align*}
Let $A := \int_{D} \eta^{2}  \Phi^{\frac{p-2+\beta}{2}} |\nabla^{2}_{0}w_{p}|^{2} dx $ and $B:= \int_{D}  \Phi^{\frac{p+\beta}{2}} dx$. By Lemma \ref{lem5.9} and the above inequality, we get
\begin{align*}
    A \leq& C_{3} \frac{\beta +1}{(p-1)^{2}}K_{\eta}B  + C_{4}\frac{(\beta +1)^{2}}{(p-1)^{2}}  \left[\frac{ 2\varepsilon }{\beta +2} A + \frac{C \beta}{\beta +2} \varepsilon^{-\frac{2}{\beta}} \frac{(\beta +1)^{2}}{(p-1)^{4}}\left\|\nabla_{0} \eta\right\|^{2}_{L^{\infty}(D)}B \right] \\
    \leq & \frac{C_{10}(\beta +1)\varepsilon}{(p-1)^{2}} A + \left(\frac{C_{11}(\beta +1)}{(p-1)^{2}} + \frac{C_{12}(\beta +1)^{4}}{(p-1)^{6}}\varepsilon^{-\frac{2}{\beta}} \right) K_{\eta} B
\end{align*}
where $C_{10}$, $C_{11}$ and $C_{12}$ are absolute constants. Choose $\varepsilon = \frac{(p-1)^{2}}{2C_{10}(\beta +1)} $, we have
\begin{align*}
    A \leq 2 \left(\frac{C_{11}(\beta +1)}{(p-1)^{2}}+ \frac{C_{12}(\beta +1)^{4}}{(p-1)^{6}} \left(\frac{2C_{10}(\beta +1)}{(p-1)^{2}}\right)^{\frac{2}{\beta}} \right)K_{\eta}  B
\end{align*}
From $\beta \geq 2$ and $1<p<2$, there exists an absolute constant $C$ such that
\begin{align*}
    A \leq C K_{\eta}\frac{(\beta +1)^{5}}{(p-1)^{8}} B.
\end{align*}
This completes the proof of Lemma \ref{lem5.11}.
\hfill${\square}$

\subsection{ Proofs of Theorems \ref{thm5.13} and \ref{thm5.14}}
\
   \vglue-10pt
    \indent

\noindent {\bf Proof of Theorem \ref{thm5.13}}: Fix a $CC$-ball $B_{R} \subset \Omega$. Since $C^{\infty}\left(B_{R}\right) \cap H W^{1, p}\left(B_{R}\right)$ is dense in $H W^{1, p}\left(B_{R}\right)$, there exists a sequence of regular functions $\psi_{\varepsilon} \in C^{\infty}\left(B_{R}\right) \cap H W^{1, p}\left(B_{R}\right)$ that converges to $w_{p}$ in $H W^{1, p}\left(B_{R}\right)$. Take an Euclidean ball $B_{\sigma r}^E \subset B_{R}$: it is satisfies condition (\ref{5.9}). Consider the Dirichlet problem
\begin{equation}\label{9.47}
    \begin{cases}
        \operatorname{div}_{0}\left(\left(\delta^2+\left|\nabla_{0} v\right|^2\right)^{\frac{p-2}{2}} \nabla_{0} v\right)=0, \qquad  \text { in } B_{\sigma r}^E \\
    v-\psi_{\varepsilon} \in H W_0^{1, p}\left(B_{\sigma r}^E\right)
    \end{cases}
\end{equation}
Let $w_{p}^{\varepsilon}$ be the weak solution of (\ref{9.47}). Taking $\varphi_{\varepsilon}=w_{p}^{\varepsilon}-\psi_{\varepsilon}$ as a test function in the weak formulation of (\ref{5.6}), we get
\begin{align*}
    \int_{B_{\sigma r}^E} \left(\delta^2+\left|\nabla_{0} w_{p}^{\varepsilon}\right|^2\right)^{\frac{p-2}{2}}\left|\nabla_{0} w_{p}^{\varepsilon}\right|^2 dx  &=\int_{B_{\sigma r}^E}\left(\delta^2+\left|\nabla_{0} w_{p}^{\varepsilon}\right|^2\right)^{\frac{p-2}{2}}\left\langle\nabla_{0} w_{p}^{\varepsilon}, \nabla_{0} \psi_{\varepsilon}\right\rangle dx \notag \\
    &\leq \int_{B_{\sigma r}^E}\left(\delta^2+\left|\nabla_{0} w_{p}^{\varepsilon}\right|^2\right)^{\frac{p-1}{2}}\left|\nabla_{0} \psi_{\varepsilon}\right| d x .
\end{align*}
Using the above inequality and Young's inequality with exponent $\frac{p}{p-1}$, we have
\begin{align*}
    \int_{B_{\sigma r}^E}\left(\delta^2+\left|\nabla_{0} w_{p}^{\varepsilon}\right|^2\right)^{\frac{p}{2}} dx = & \int_{B_{\sigma r}^E}\left(\delta^2+\left|\nabla_{0} w_{p}^{\varepsilon}\right|^2\right)^{\frac{p-2}{2}} \delta^2 d x +\int_{B_{\sigma r}^E}\left(\delta^2+\left|\nabla_{0} w_{p}^{\varepsilon}\right|^2\right)^{\frac{p-2}{2}}\left|\nabla_{0} w_{p}^{\varepsilon}\right|^2 \\
    \leq & \int_{B_{\sigma r}^E} \delta^{p} dx +\int_{B_{\sigma r}^E}\left(\delta^2+\left|\nabla_{0} w_{p}^{\varepsilon}\right|^2\right)^{\frac{p-1}{2}}\left|\nabla_{0} \psi_{\varepsilon}\right| dx \\
    \leq & \int_{B_{\sigma r}^E} \delta^p dx +\frac{p-1}{p} \int_{B_{\sigma r}^E}\left(\delta^2+ \left|\nabla_{0} w_{p}^{\varepsilon}\right|^2\right)^{\frac{p}{2}} dx+\frac{1}{p}\int_{B_{\sigma r}^E}\left|\nabla_{0} \psi_{\varepsilon}\right|^p dx.
\end{align*}
Thus,
\begin{align}\label{9.48}
    \int_{B_{\sigma r}^E}\left(\delta^2+\left|\nabla_{0} w_{p}^{\varepsilon}\right|^2\right)^{\frac{p}{2}} dx \leq  \int_{B_{\sigma r}^E}  p\delta^p  + \left|\nabla_{0} \psi_{\varepsilon}\right|^p dx.
\end{align}
From $|\nabla_{0} \psi_{\varepsilon}|^{p} \leq 2^{p-1}\left(|\nabla_{0} \psi_{\varepsilon} - \nabla_{0} w_{p} |^{p} + |\nabla_{0} w_{p}|\right) $, we have
\begin{align*}
    p\delta^p  + \left|\nabla_{0} \psi_{\varepsilon}\right|^p &\leq p \delta^{p} + 2^{p-1}\left(|\nabla_{0} \psi_{\varepsilon} - \nabla_{0} w_{p} |^{p} + |\nabla_{0} w_{p}|^{p}\right) \\
    &\leq 2|\nabla_{0} \psi_{\varepsilon} - \nabla_{0} w_{p} |^{p}  + 2 \left(\delta^{2} + |\nabla_{0} w_{p}|^{2}\right)^{\frac{p}{2}} .
\end{align*}
Hence, (\ref{9.48}) implies that
\begin{align*}
    \int_{B_{\sigma r}^E}\left(\delta^2+\left|\nabla_{0} w_{p}^{\varepsilon}\right|^2\right)^{\frac{p}{2}} dx \leq 2 \int_{B_{\sigma r}^E}\left(\delta^2+\left|\nabla_{0} w_{p}\right|^2\right)^{\frac{p}{2}}+ 2 \int_{B_{\sigma r}^E}\left|\nabla_{0} w_{p}-\nabla_{0} \psi_{\varepsilon}\right|^p dx .
\end{align*}
Since $\psi_{\varepsilon} \in C^{\infty}\left(B_{R}\right) \cap H W^{1, p}\left(B_{R}\right)$ and $\psi_{\varepsilon} \to w_{p} $ in $H W^{1, p}\left(B_{R}\right)$ as $\varepsilon \to 0$, we have
\begin{align}\label{9.49}
    \int_{B_{\sigma r}^E}\left(\delta^2+\left|\nabla_{0} w_{p}^{\varepsilon}\right|^2\right)^{\frac{p}{2}} dx \leq 2 \int_{B_{\sigma r}^E}\left(\delta^2+\left|\nabla_{0} w_{p}\right|^2\right)^{\frac{p}{2}} + o(\varepsilon)
\end{align}
where $o(\varepsilon) \to 0$ as $\varepsilon \to 0$. This implies $\left\| \nabla_{0} w_{p}^{\varepsilon}\right\|_{L^{p}(B_{\sigma r}^E)}$ is uniformity bounded with respect to $\varepsilon$. Using (\ref{5.s}) with $\kappa =1$ for $w_{p}^{\varepsilon}$ in $B_{\sigma r}^E \subset B_{R}$, it follows that $\left\|w_{p}^{\varepsilon}\right\|_{HW^{1,p}(B_{\sigma r}^E)}$ is uniformity bounded with respect to $\varepsilon$. Then, we can extract a weakly convergent subsequence of $\{w_{p}^{\varepsilon}\} $ in $H W^{1, p}\left(B_{\sigma r}^E\right)$ , and there exists a function $\bar{w}_{p} \in H W^{1, p}\left(B_{\sigma r}^E\right)$ such that $w_{p}^{\varepsilon_{j}} \rightharpoonup \bar{w}_{p}$ in $H W^{1, p}\left(B_{\sigma r}^E\right)$. Since $w_{p}^{\varepsilon}-\psi_{\varepsilon} \in H W_0^{1, p}\left(B_{\sigma r}^E\right)$ and this space is closed under weak convergence, we obtain $\bar{w}_{p}-w_{p} \in H W_0^{1, p}\left(B_{\sigma r}^E\right)$.

Using Lemma 3.6 in \cite{R2015}, we get that $\nabla_{0} w_{p}^{\varepsilon} \rightarrow \nabla_{0} \bar{w}_{p}$ in $L^p\left(B_{\sigma r}^E\right)$ and $\bar{w}_{p}$ is a weak solution of (\ref{9.47}) in $B_{\sigma r}^E$. Combining  $\bar{w}_{p}-w_{p} \in H W_0^{1, p}\left(B_{\sigma r}^E\right)$ with the uniqueness theorem the Dirichlet problem (\ref{9.47}) (cf. Theorem 3.3 in \cite{R2015}), we get $\bar{w}_{p}=w_{p}$.

Let $B_{2 \lambda R} \subset B_{\sigma r}^E$ with $\lambda >0$. Then, from Theorem \ref{thm5.12}, we have
\begin{align}\label{9.50}
    \left\|\nabla_{0} w_{p}^{\varepsilon}\right\|_{L^{\infty}\left(B_{\lambda R}\right)} \leq \mathcal{C}^{l}(c,Q)(p-1)^{-\frac{4Q}{p}} \left(\fint_{B_{2 \lambda R}}\left(\delta^2+\left|\nabla_{0} w_{p}^{\varepsilon}\right|^2\right)^{\frac{p}{2}} dx\right)^{\frac{1}{p}}.
\end{align}
Since $\nabla_{0} w_{p}^{\varepsilon} \rightarrow \nabla_{0} w_{p}$ in $L^p\left(B_{\sigma r}^E\right)$, we can pass to a subsequence $\{w_{p}^{\varepsilon_{j}}\}$ such that $\nabla_{0} w_{p}^{\varepsilon_{j}} \rightarrow \nabla_{0} w_{p}$ pointwisely a.e. in $B_{\sigma r}^E$. We can conclude from (\ref{9.49}) and (\ref{9.50}) that
\begin{align*}
    \left\|\nabla_{0} w_{p}\right\|_{L^{\infty}(B_{ \lambda R})} \leq \mathcal{C}^{l}(c,Q) (p-1)^{-\frac{4Q}{p}}\left(\fint_{B_{R}}  \left(\delta^{2} + |\nabla_{0} w_{p}|^{2}\right)^{\frac{p}{2}} \right)^{\frac{1}{p}}.
\end{align*}
By a covering argument, the proof of Theorem \ref{thm5.13} is complete.
\hfill${\square}$

\noindent {\bf Proof  of Theorem \ref{thm5.14}}: Let $w_{p}^{\delta}$ be the unique weak solution of the Dirichlet problem
\begin{equation*}
    \begin{cases}
        \operatorname{div}_{0}\left(\left(\delta^2+\left|\nabla_{0} w_{p}^{\delta}\right|\right)^{\frac{p-2}{2}} \nabla_{0} w_{p}^{\delta}\right)=0,\qquad  \text { in } \Omega \\
        w_{p}^{\delta}-w_{p} \in H W_0^{1, p}(\Omega)
    \end{cases}
\end{equation*}
where $\delta >0$ and $w_{p} \in H W_0^{1, p}(\Omega)$ is a weak solution of (\ref{5.4a})-(\ref{5.4b}). It is obvious that $w_{p}$ is a minimum of the functional $\mathcal{L}_{w_{p},0} $ in $ \mathcal{A}_{1}= \left\{v \in HW^{1,p}(\Omega) , v = 1 \text{ a.e.  on }  \partial \Omega \right\}$. In addition, $w_{p}^{\delta}$ is the minimum of the functional $\mathcal{L}_{w_{p},\delta}$. Thus
\begin{align*}
    \mathcal{L}_{w_{p},0}(w_{p}) \leq \mathcal{L}_{w_{p},0}\left(w_{p}^{\delta}\right) \leq \mathcal{L}_{w_{p}, \delta}\left(w_{p}^{\delta}\right) \leq \mathcal{L}_{w_{p}, \delta}(w_{p})
\end{align*}
and subtracting $\mathcal{L}_{w_{p},0}(w_{p})$, we obtain
\begin{align*}
    0 \leq \mathcal{L}_{w_{p},0}\left(w_{p}^{\delta}\right)  - \mathcal{L}_{w_{p},0}(w_{p}) \leq \mathcal{L}_{w_{p}, \delta}(w_{p})- \mathcal{L}_{w_{p},0}(w_{p}),
\end{align*}
i.e.,
\begin{align*}
    0 \leq \int_{\Omega}\left|\nabla_{0} w_{p}^{\delta}\right|^p dx-\int_{\Omega}\left|\nabla_{0} w_{p}\right|^p dx \leq \int_{\Omega}\left(\delta^2+\left|\nabla_{0} w_{p}\right|^2\right)^{\frac{p}{2}} dx-\int_{\Omega}\left|\nabla_{0} w_{p}\right|^p dx .
\end{align*}
By Lebesgue's dominated convergence theorem, the last term in the above
inequality tends to zero, thus $\left\|\nabla_{0} w_{p}^{\delta}\right\|_{L^p(\Omega)}\to \left\|\nabla_{0} w_{p}\right\|_{L^p(\Omega)}$ as $\delta \to 0$. This means that $w_{p}^{\delta}$ is a minimizing sequence for the functional $\mathcal{L}_{w_{p},0}$. Hence, for any $\varphi  \in HW^{1,p}(\Omega)$, we have
\begin{align*}
    0 = \int_{\Omega} |\nabla_{0} w_{p}^{\delta}|^{p-2} \left\langle \nabla_{0} w_{p}^{\delta}, \nabla_{0} \varphi \right\rangle  ,\quad 0 = \int_{\Omega} |\nabla_{0} w_{p}|^{p-2} \left\langle \nabla_{0} w_{p}, \nabla_{0} \varphi \right\rangle ,
\end{align*}
which implies  $0 =  \int_{\Omega} \left\langle |\nabla_{0} w_{p}^{\delta}|^{p-2}\nabla_{0} w_{p}^{\delta} -  |\nabla_{0} w_{p}|^{p-2}\nabla_{0} w_{p}, \nabla_{0} \varphi  \right\rangle $. Let $\varphi = w_{p}^{\delta} - w_{p}$, we get
\begin{align*}
    0 =  \int_{\Omega} \left\langle |\nabla_{0} w_{p}^{\delta}|^{p-2}\nabla_{0} w_{p}^{\delta} -  |\nabla_{0} w_{p}|^{p-2}\nabla_{0} w_{p}, \nabla_{0}w_{p}^{\delta} - \nabla_{0} w_{p}  \right\rangle .
\end{align*}
Applying Lemma 3.2 and Lemma 3.3 in \cite{R2015} to the above equality, we obtain
\begin{align*}
    0 \geq  \frac{(p-1)2^{\frac{2-p}{2}}}{p}\int_{\Omega} \left(|\nabla_{0} w_{p}^{\delta}|^{2} + |\nabla_{0} w_{p}|^{2}\right)^{\frac{p-2}{2}} |\nabla_{0}w_{p}^{\delta} - \nabla_{0} w_{p} |^{2}.
\end{align*}
Thus, $|\nabla_{0}w_{p}^{\delta} - \nabla_{0} w_{p} |^{2} = 0$ a.e. in $\Omega$, i.e., $\left\|\nabla_{0}w_{p}^{\delta} - \nabla_{0} w_{p}\right\|_{L^{2}(\Omega)} = 0$. Using (\ref{5.s}) with $\kappa =1$ and $p=2$ for $(w_{p}^{\delta} - w_{p})$, we obtain
\begin{align*}
    \left\|w_{p}^{\delta} - w_{p}\right\|_{L^{2}(B_{R})} \leq C(c,Q) R \left\|\nabla_{0} (w_{p}^{\delta} -w_{p}) \right\|_{L^{2}(B_{R})} = 0
\end{align*}
Therefore, $w_{p}^{\delta}  = w_{p} $ a.e. in $B_{R}$.

Consider a minimizing sequence of $w_{p}^{\delta_{j}} \in \mathcal{A}_{1}$. Let $\lambda = \inf_{w_{p}^{\delta_{j}} \in \mathcal{A}_{1}} \mathcal{L}_{w_{p},0}$, and observe that $0 \leq \lambda \leq  \mathcal{L}_{w_{p},0}(w_{p}) < \infty$. Thus, $\mathcal{L}_{w_{p},0}(w_{p}^{\delta_{j}}) \leq  \lambda + \frac{1}{j}$ which implies
\begin{align*}
    \int_{B_{R}} |\nabla_{0}w_{p}^{\delta_{j}}|^{p} \leq  p \mathcal{L}_{w_{p},0}(w_{p}^{\delta_{j}}) \leq  p \left(\lambda + \frac{1}{j}\right), \qquad \forall j.
\end{align*}
This yields that $\left\|\nabla_{0}w_{p}^{\delta_{j}}\right\|_{L^{p}(B_{R})}$ is uniformly bounded with respect to $\delta_{j}$. Using (\ref{5.s})  with $\kappa =1$ for $w_{p}^{\delta_{j}}$, we have $\left\|w_{p}^{\delta_{j}}\right\|_{L^{p}(B_{R})} \leq C(c,Q)R \left\|\nabla_{0}w_{p}^{\delta_{j}}\right\|_{L^{p}(B_{R})}$. Thus, $\left\|w_{p}^{\delta_{j}}\right\|_{HW^{1,p}(B_{R})}$ is uniformly bounded with respect to $\delta_{j}$. Since $HW^{1,p}(B_{R})$ is reflexive, there exist a subsequence, still denoted by $w_{p}^{\delta_{j}}$, and a function $\hat{w}_{p} \in HW^{1,p}(B_{R})$ such that $w_{p}^{\delta_{j}}  \rightharpoonup \hat{w}_{p}$ in $HW^{1,p}(B_{R}) $. Moreover, $\hat{w}_{p} \in \mathcal{A}_{1}$. By the weakly lower semicontinuity of the functional $\mathcal{L}_{w_{p},0}$, it follows that $\mathcal{L}_{w_{p},0}(\hat{w}_{p}) \leq \liminf_{j \to \infty}\mathcal{L}_{w_{p},0}(w_{p}^{\delta_{j}}) \leq \lim_{j \to \infty} \left(\lambda + \frac{1}{j}\right) = \lambda$. Hence, $\hat{w}_{p}$ is a minimizer of the functional $\mathcal{L}_{w_{p},0}$.

By the uniqueness of solutions to the Dirichlet problem (\ref{5.4a})-(\ref{5.4b}), we get $\hat{w}_{p}=w_{p}$. Consequently, $w_{p}^{\delta_{j}}  \rightharpoonup  w_{p}$ in $HW^{1,p}(B_{R})$.  Combining with $\left\|\nabla_{0} w_{p}^{\delta_{j}}\right\|_{L^p(B_{R})} \to \left\|\nabla_{0} w_{p}\right\|_{L^p(B_{R})}$, we have $\nabla_{0} w_{p}^{\delta_{j}} \to  \nabla_{0} w_{p} $ strongly in  $L^{p}(B_{R}) $. Thus
\begin{align}\label{9.51}
    \left\|\nabla_{0}(w_{p}^{\delta_{j}} - w_{p})\right\|_{ L^{p}(B_{R})} \to 0 ,\qquad \text{ as } \delta_{j} \to 0.
\end{align}
Using (\ref{5.s}) with $\kappa =1$ for $(w_{p}^{\delta_{j}} - w_{p})$, we have
\begin{align*}
    \left\|w_{p}^{\delta_{j}} - w_{p}\right\|_{ L^{p}(B_{R})} \leq  C(c,Q) R \left\|\nabla_{0}(w_{p}^{\delta_{j}} - w_{p})\right\|_{ L^{p}(B_{R})} \to 0, \qquad \text{ as } \delta_{j} \to 0.
\end{align*}
Thus, $w_{p}^{\delta_{j}} \to  w_{p}$ in $HW^{1,p}(B_{R})$. In addition, by (\ref{9.51}), we can extract a subsequence such that $\nabla_{0} w_{p}^{\delta_{j}}$ converges to $ \nabla_{0} w_{p}$ for a.e. $x\in B_{R}$. Hence,
\begin{align*}
    |\nabla_{0} w_{p}^{\delta_{j}}|  \to  |\nabla_{0} w_{p}| \, \text{ and } \, |\nabla_{0} w_{p}(x)| \leq \liminf_{\delta_{j} \to 0}\|\nabla_{0} w_{p}^{\delta_{j}}\|_{L^{\infty}(B_{R})}, \qquad \text{ for } a.e. x \in B_{R}.
\end{align*}
By Lebesgue's dominated convergence theorem and Theorem \ref{thm5.13},  we obtain
\begin{align*}
    \|\nabla_{0} w_{p}\|_{L^{\infty}(B_{R})} &\leq  \liminf_{\delta_{j} \to 0} \|\nabla_{0} w_{p}^{\delta_{j}}\|_{L^{\infty}(B_{R})}\\
    &\leq \mathcal{C}^{l}(c,Q) (p-1)^{-\frac{4Q}{p}} \lim_{\delta_{j} \to 0} \left(\fint_{B_{R}}  \left(\delta^{2}_{j} + |\nabla_{0} w_{p}^{\delta_{j}}|^{2}\right)^{\frac{p}{2}}\right)^{\frac{1}{p}} \\
    &= \mathcal{C}^{l}(c,Q) (p-1)^{-\frac{4Q}{p}}\left(\fint_{B_{2R}}  |\nabla_{0} w_{p}|^{p} \right)^{\frac{1}{p}}.
\end{align*}
This completes the proof of Theorem \ref{thm5.14}.
\hfill${\square}$



\end{document}